%% file: main.tex
\documentclass[final,3p]{elsarticle}
\makeatletter 
\def\ps@pprintTitle{%
 \def\@oddfoot{\footnotesize\itshape
       Published in \ifx\@journal\@empty Elsevier
       \else\@journal\fi\hfill March 1, 2018}}
\makeatother

\graphicspath{{./graphics/}}
\input{header}
\input{header_Article}

\journal{Computer methods in applied mathematics and engineering}

\begin{document}

\begin{frontmatter}

\title{Isogeometric Analysis of Acoustic Scattering using Infinite Elements}

\author[ntnu]{Jon Vegard Ven{\aa}s\texorpdfstring{\corref{cor}}{}}
\ead{Jon.Venas@ntnu.no}

\author[ntnu]{Trond Kvamsdal}
\ead{Trond.Kvamsdal@ntnu.no}

\author[ffi]{Trond Jenserud}
\ead{Trond.Jenserud@ffi.no}

\address[ntnu]{Department of Mathematical Sciences, Norwegian University of Science and Technology, Trondheim, Norway}
\address[ffi]{Department of Marine Systems, Norwegian Defence Research Establishment, Horten, Norway}
\cortext[cor]{Corresponding author.}

\begin{abstract}
\input{\contents/abstract}

\end{abstract}


\end{frontmatter}

\input{\contents/introduction}
\input{\contents/exteriorHelmholtz}
\input{\contents/coupledFluidStructure}
\input{\contents/resultsAndDiscussion}
\input{\contents/conclusion}
\appendix
\input{\contents/infiniteElementsBilinearForm}
\input{\contents/radialIntegrals}
\input{\contents/betssi_description}
\input{\contents/transformNURBStoBspline}

%\clearpage
\section*{References}
\bibliographystyle{TK_CM} % same as elsarticle-num without the display of url and doi (the title contains the url info. If doi is given the link in the title should perhaps use this link instead?)
%\nocite{*}
\bibliography{references}
%\bibliography{references}


%\begin{thebibliography}{00}
%%% Text of bibliographic item
%
%\bibitem{}
%
%\end{thebibliography}
\end{document}

%% file: header.tex
\usepackage{amsmath,amsfonts,amssymb} % for maths
\usepackage{mathtools} % contains amsthm
\usepackage{upgreek} % to get greek letters in upright mode, like \upphi
\usepackage{easybmat} % to easily create large block matrices
\usepackage{fixmath} % provides upercase greek: \upxi,\upeta,...
\usepackage[font=footnotesize,list=true]{subcaption} % Needed for subfigures
\usepackage{booktabs} % needed for \toprule, \bottomrule,... in tables
\usepackage[exponent-product=\cdot]{siunitx} % for nice SI unit formatting
\renewcommand{\ang}[1]{{{#1}^\circ}}
\usepackage{bm} % to get nice symbols in boldface mode
\usepackage{xcolor}
\usepackage{microtype}
\usepackage{lmodern}

\newcommand{\graphicsFolder}{graphics}

\renewcommand{\vec}[1]{\mathbold #1} % ISO80000-2 format
\newcommand{\diff}{\mathrm{d}} % For differentials
\newcommand{\idiff}{\, \mathrm{d}} % For differentials in integrals.
\newcommand{\zerovec}{\mathbf{0}}
\newcommand{\transpose}{\top}

\newcommand{\R}{\mathbb{R}}

\newcommand{\PI}{\uppi}
\newcommand{\euler}{\mathrm{e}}
\newcommand{\imag}{\mathrm{i}}
\newcommand{\bigoh}{\mathcal{O}}
\providecommand*{\pderiv}[3][]{\frac{\partial^{#1}#2}{\partial #3^{#1}}}
\providecommand*{\deriv}[3][]{\frac{\diff^{#1}#2}{\diff #3^{#1}}}

\renewcommand{\geq}{\geqslant}
\renewcommand{\leq}{\leqslant}

\DeclareMathSymbol{\GAMMA}{\mathalpha}{operators}{0}

\DeclareMathOperator{\TS}{TS}

\renewcommand\Re{\operatorname{Re}}

\newcommand{\calV}{\mathcal{V}}
\newcommand{\calS}{\mathcal{S}}
\newcommand{\calI}{\mathcal{I}}
\newcommand{\calF}{\mathcal{F}}
\newcommand{\calH}{\mathcal{H}}
\renewcommand{\Xi}{{\vec{t}_1}}
\newcommand{\Eta}{{\vec{t}_2}}
\newcommand{\Zeta}{{\vec{t}_3}}

\newcommand{\energyNorm}[2]{%
  {\left\vert\kern-0.25ex\left\vert\kern-0.25ex\left\vert #1 
    \right\vert\kern-0.25ex\right\vert\kern-0.25ex\right\vert}_{#2}
}

\let\originalleft\left
\let\originalright\right
\renewcommand{\left}{\mathopen{}\mathclose\bgroup\originalleft}
\renewcommand{\right}{\aftergroup\egroup\originalright}

\usepackage{xspace}

%% file: header_Article.tex
\usepackage{hyperref} % needed for links
\hypersetup{hidelinks,colorlinks=true,citecolor=blue,linkcolor=blue,urlcolor=blue} % Color the links
\usepackage[noabbrev]{cleveref} % To have easy referencing for figures, equations, ... noabbrev option can be used for full name
\crefname{equation}{Eq.}{Eqs.}
\usepackage{longtable}
\usepackage{enumitem}

%% file: contents/abstract.tex
Isogeometric analysis (IGA) has proven to be an improvement on the classical finite element method (FEM) in several fields, including structural mechanics and fluid dynamics. In this paper, the performance of IGA coupled with the infinite element method (IEM) for some acoustic scattering problems is investigated. In particular, the simple problem of acoustic scattering by a rigid sphere, and the scattering of acoustic waves by an elastic spherical shell with fluid domains both inside and outside, representing a full acoustic-structure interaction (ASI) problem. Finally, a mock shell and a simplified submarine benchmark are investigated. The numerical examples include comparisons between IGA and the FEM. Our main finding is that the usage of IGA significantly increases the accuracy compared to the usage of $C^0$ FEM due to increased inter-element continuity of the spline basis functions.

%% file: contents/introduction.tex
\section{Introduction}
Acoustic scattering is the physical phenomena of how sound interacts with objects and medium fluctuations. When an acoustic wave hits a rigid object, it is totally reflected, and the object is left in a quiescent state. In the case of an elastic object, part of the sound is transmitted into the object, which is set into motion and starts radiating sound. This leads to a coupled acoustic-structure interaction (ASI) problem.
Applications include underwater acoustics~\cite{Gilroy2013bib} and noise propagation in air~\cite{Bouillard1999eea}. Inverse problems are also of interest, such as shape optimization of membranes~\cite{Manh2011iso} and the problem of designing submarines with low scattering strength. Assuming harmonic time dependency, the fluid and solid media can be modeled using the scalar and vector Helmholtz equations, respectively. The vector Helmholtz equation can be used to model electromagnetic waves~\cite{Manh2012iaa}, such that the work presented herein can also be used for electromagnetic scattering.

Herein, the acoustic scattering characterized by sound waves reflected by man-made elastic objects will be addressed. Shape optimization for optimal acoustic scattering on man-made objects, e.g. antennas, submarines etc., is a typical problem facing design engineers.

Isogeometric analysis (IGA) is basically an extension of the finite element method (FEM) using non-uniform rational B-splines (NURBS) as basis functions not only representing the solution space, but also the geometry. Being introduced in 2005 by Hughes et al.~\cite{Hughes2005iac}, followed by the book~\cite{Cottrell2009iat} in 2009, IGA tries to bridge the gap between finite element analysis (FEA) and computer aided design (CAD) tools. The important feature of IGA is that it uses the same basis as CAD software for describing the given geometry, and thus exact representation of the model is possible.

\begin{figure}
	\centering
	\includegraphics[scale=1]{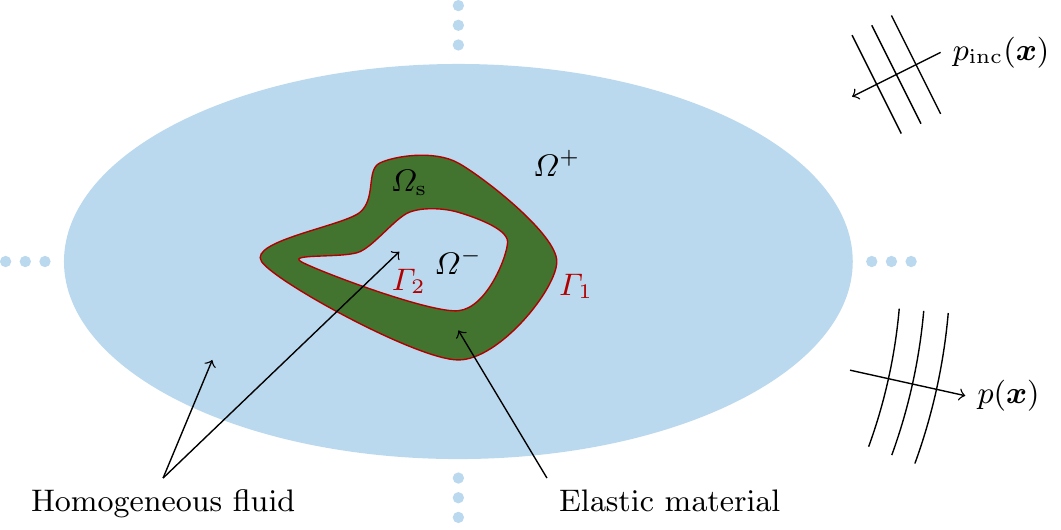}
	\caption{Illustration of the physical problem. A plane incident wave, $p_{\mathrm{inc}}(\vec{x})$, is scattered by the scatterer, $\Omega_{\mathrm{s}}$, in an unbounded domain, $\Omega^+$, resulting in the scattered wave, $p(\vec{x})$. The scatterer, which is bounded by the boundaries $\Gamma_1$ and $\Gamma_2$, envelops a fluid domain, $\Omega^-$.}
	\label{Fig2:physicalProblem}
\end{figure}
The physical problem is illustrated in \Cref{Fig2:physicalProblem} where the incoming sound waves, $p_{\rm inc}$, originate from a point source far from this object, such that the (spherical) sound waves are quite accurately approximated by plane waves when the waves reaches the proximity of the object. For rigid objects of irregular shape, the incoming wave may be reflected multiple times before leaving the object. When the object is elastic a coupled ASI problem results. The goal is then to calculate the scattered wave $p$ at an arbitrary far field point. Finally, to use the FEM or IGA the domain must be finite. A fictitious boundary is thus introduced, which must be implemented in such a way that outgoing waves reaching this boundary are absorbed.

The problem at hand is time dependent. However, harmonic time dependency will be assumed, such that all time dependent functions may be written as $\breve{F}=\breve{F}(\vec{x},t) = F(\vec{x})\euler^{-\imag\omega t}$ where $\omega$ is the angular frequency and $\imag = \sqrt{-1}$ the imaginary unit. This enables us to model the pressure $p$ in the fluid with the Helmholtz equation given by
\begin{equation}\label{Eq2:HelmholtzEquationIntro}
	\nabla^2 p + k^2 p = 0
\end{equation}
with the wave number $k=\frac{\omega}{c_{\mathrm{f}}}$ (where $c_{\mathrm{f}}$ is the wave speed in the fluid). Other important quantities include the frequency $f=\frac{\omega}{2\PI}$ and the wavelength $\lambda = \frac{2\PI}{k}$.

\begin{figure}
	\centering
	\includegraphics[width=0.3\textwidth]{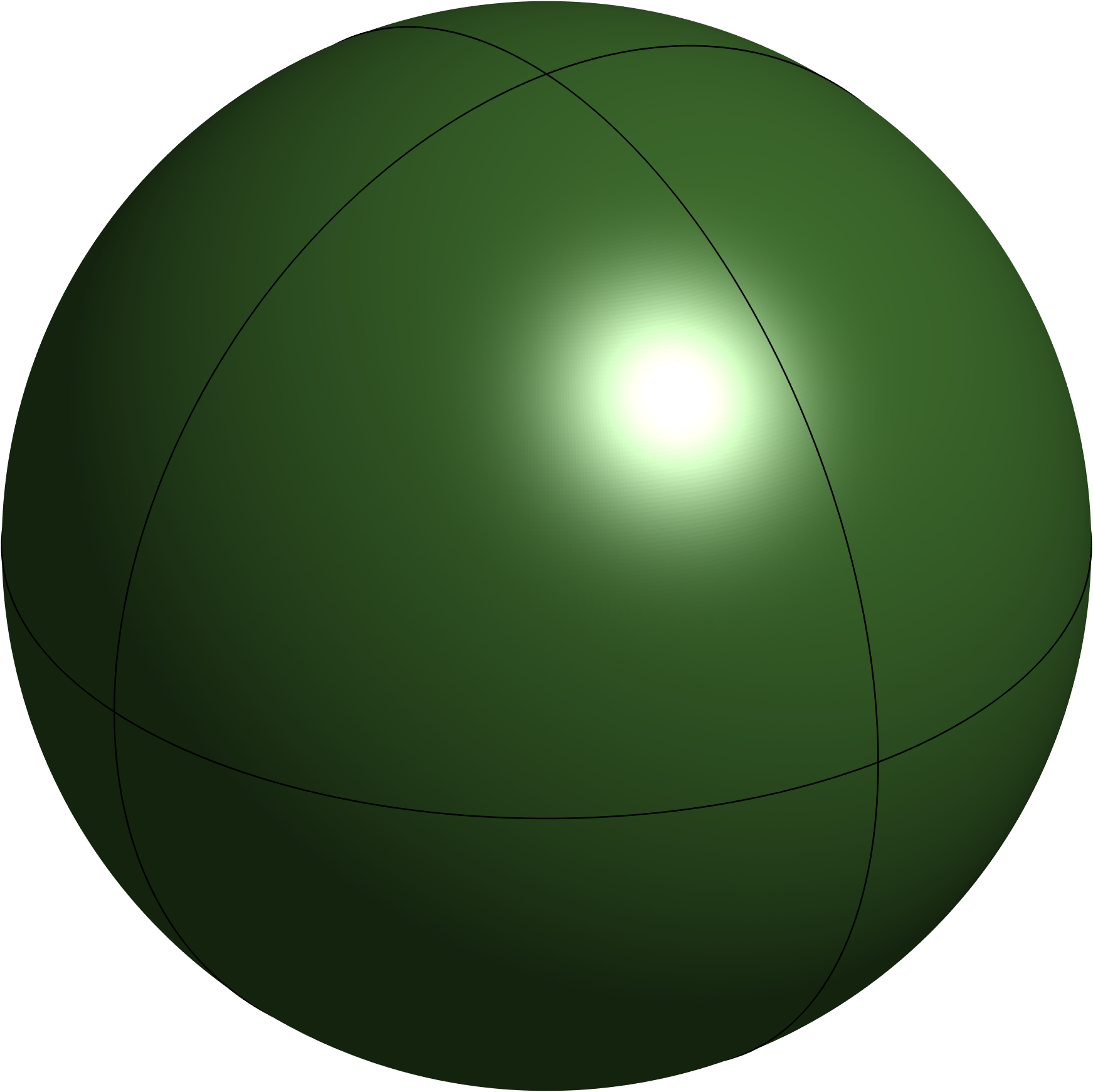}
	~
	\includegraphics[width=0.25\textwidth]{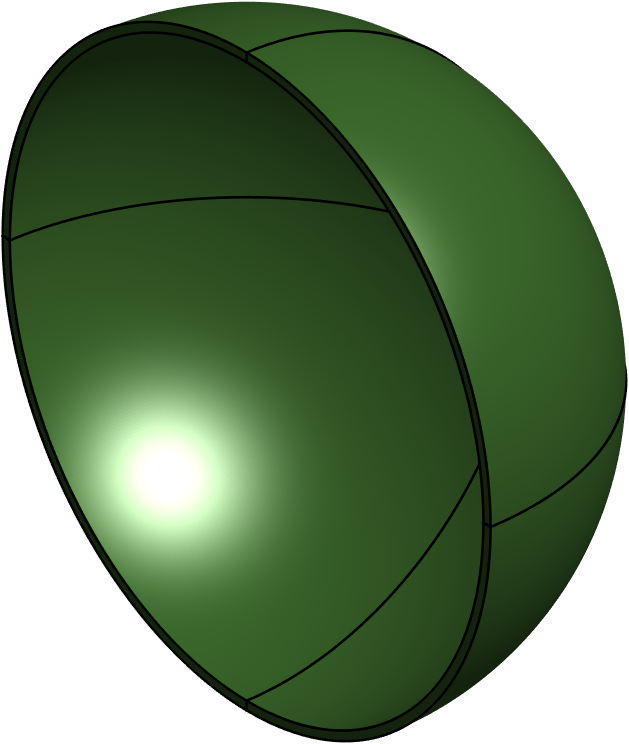}
	\caption{Exact geometry of a spherical shell using 8 elements.}
	\label{Fig2:SphericalShell}
\end{figure}

The geometry of the elastic object may be quite complex, but is typically exactly represented using NURBS. This fact is one of the motivating factor for using IGA, as it uses the same functions as basis functions for analysis. The spherical shell depicted in \Cref{Fig2:SphericalShell} is an example of a geometry that has an exact representation using NURBS, but is outside the space of standard (Lagrangian) FEM geometries.

It has been shown that the continuity of the basis functions plays an important role for the accuracy of solving elliptical problems (for instance the Helmholtz equation), see~\cite{BeiraodaVeiga2011sef} and~\cite{BeiraodaVeiga2014mao}. This motivates the use of IGA even further, as IGA enables control of the continuity of the basis function up to $C^{p-1}$ (in contrast with the $C^0$-continuity restriction in FEA). IGA has proven to be promising in a host of areas related to the problem at hand, which yields further motivation in the use of IGA. For instance, in~\cite{nortoft2015iao} the method was shown to be suited for the more complex scenario of sound propagation through laminar flow.

In addition to IGA, the so-called infinite element method (IEM) has been chosen to handle the boundary conditions at the artificial boundary. Typically, the boundary element method (BEM) \cite{Sauter2011bem,Schanz2007bea} has been used for this purpose. However, for higher frequencies and complex geometries, BEM becomes computationally expensive (although improvement in performance has been done in the recent decades \cite{Liu2012raa}). The main motivation for the infinite element method is computational efficiency as reported by Burnett~\cite{Burnett1994atd} and Gerdes and Demkowicz~\cite{Gerdes1996so3}.

Before starting on the full ASI problem, it is important to establish good results for the IEM. This method only applies for the outer fluid, and it would thus be natural to first investigate the scattering problem on rigid objects (that is, no acoustic-structure interaction occurs). An introduction to the IEM is presented in \Cref{Sec2:exteriorHelmholtz}. The extension to ASI problems (presented in~\Cref{Sec2:coupledFluidStruct}) naturally follows from the implementation of rigid scattering using IEM. In \Cref{Sec2:resultsDisc} the results obtained for both rigid and elastic scattering on a spherical shell is presented. Results for rigid scattering from a mock shell are included to investigate condition numbers. Moreover, results for a simplified submarine is presented to illustrate the performance of the implementation on complex geometries. Finally, conclusions and suggested future work can be found in \Cref{Sec2:conclusions}. 

%% file: contents/exteriorHelmholtz.tex
\section{Exterior Helmholtz problems}
\label{Sec2:exteriorHelmholtz}
\begin{figure}
	\centering
	\includegraphics[scale=1]{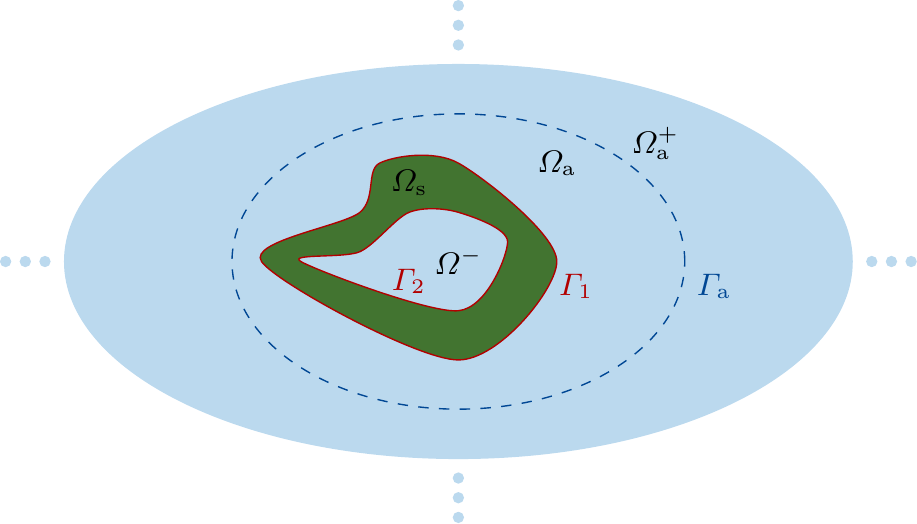}
	\caption[Illustration of artificial boundary]{An artificial boundary $\Gamma_{\mathrm{a}}$ is introduced such that the exterior domain $\Omega^+$ is decomposed by the two domains $\Omega_{\mathrm{a}}$ (which is bounded by $\Gamma_1$ and $\Gamma_{\mathrm{a}}$) and $\Omega_{\mathrm{a}}^+$. Thus, $\Omega^+ = \Omega_{\mathrm{a}} \cup \Omega_{\mathrm{a}}^+$.}
	\label{Fig2:artificialBoundary}
\end{figure}

Scattering problems involve \textit{unbounded exterior domains}, $\Omega^+$. A common method for solving such problems with the FEM is to introduce an artificial boundary that encloses the scatterer. On the artificial boundary some sort of absorbing boundary condition (ABC) is prescribed. The problem is then reduced to a finite domain, and both the elastic scatterer and the bounded domain between the scatterer and the artificial boundary can be discretized with finite elements. Several methods exist for handling the exterior Helmholtz problem (on unbounded domain), including
\begin{itemize}
	\item the perfectly matched layer (PML) method after B{\'e}renger~\cite{Berenger1994apm,Berenger1996pml}
	\item the boundary element method~\cite{Sauter2011bem,Schanz2007bea,Marburg2008cao,Chandler_Wilde2012nab}
	\item Dirichlet to Neumann-operators (DtN-operators)~\cite{Givoli2013nmf}
	\item local differential ABC operators~\cite{Shirron1995soe,Bayliss1982bcf,Hagstrom1998afo,Tezaur2001tdf}
	\item the infinite element method.~\cite{Bettess1977ie,Bettess1977dar}
\end{itemize}
Herein, the infinite element method is chosen. For the IEM, the unbounded domain $\Omega^+$ is partitioned into two domains by the artificial boundary $\Gamma_{\mathrm{a}}$; $\Omega_{\mathrm{a}}$ and $\Omega_{\mathrm{a}}^+$ (see \Cref{Fig2:artificialBoundary}). 
These domains are discretized by finite and infinite elements, respectively. A convergence analysis of a coupled FEM-IEM can be found in \cite{Demkowicz2001aoa}.

The exterior Helmholtz problem is given by
\begin{alignat}{3}
	\nabla^2 p + k^2 p &= 0 	&&\text{in}\quad \Omega^+,\label{Eq2:HelmholtzEqn}\\
	\partial_n p &= g						&&\text{on}\quad \Gamma_1,\label{Eq2:HelmholtzEqnNeumannCond}\\
	\pderiv{p}{r}-\imag k p &= o\left(r^{-1}\right)\quad &&\text{with}\quad r=|\vec{x}|\label{Eq2:sommerfeldCond}
\end{alignat}
where the Sommerfeld condition~\cite{Sommerfeld1949pde} in \Cref{Eq2:sommerfeldCond} restricts the field in the limit $r\to\infty$ uniformly in $\hat{\vec{x}}=\frac{\vec{x}}{r}$, such that no waves originate from infinity. The Neumann condition given by the function $g$ will in the case of rigid scattering be given by the incident wave $p_{\mathrm{inc}}$. Zero displacement of the fluid normal on the scatterer (rigid scattering) implies that $\partial_n(p+p_{\mathrm{inc}}) = 0$ where $\partial_n$ denotes the partial derivative in the normal direction on the surface $\Gamma_1$ (pointing ``out'' from $\Omega^+$), which implies that
\begin{equation}
	g = -\pderiv{p_{\mathrm{inc}}}{n}.
\end{equation}
Plane incident waves (with amplitude $P_{\mathrm{inc}}$) traveling in the direction $\vec{d}_{\mathrm{s}}$ can be written as
\begin{equation}\label{Eq2:p_inc}
	p_{\mathrm{inc}} = P_{\mathrm{inc}}\euler^{\imag k\vec{d}_{\mathrm{s}}\cdot\vec{x}}.
\end{equation}
The normal derivative on the surface of any smooth geometry may then be computed by
\begin{align}
	\pderiv{p_{\mathrm{inc}}}{n} &= \vec{n}\cdot\nabla p_{\mathrm{inc}} = \imag k\vec{d}_{\mathrm{s}}\cdot\vec{n} p_{\mathrm{inc}}.
\end{align}

\subsection{Weak formulation for the Helmholtz equation}
In order to choose the correct solution space in the infinite element method, the asymptotic behavior of the scattered pressure $p$ at large radii\footnote{Here, $r$ is referred to as the radius even though it does not necessarily represent the radius in spherical coordinates.} $r$ must be examined. In~\cite{Wilcox1956aet}, Wilcox shows that the scalar pressure field $p(\vec{x})$ satisfying the Helmholtz equation and the Sommerfeld radiation conditions can be written in the form\footnote{In some appropriate coordinate system $(r,\vartheta,\varphi)$ with the ``radial variable'', $r$, extending to infinity. Typically some degeneration of the ellipsoidal (in 3D) coordinate system.}
\begin{equation}\label{Eq2:wilcoxExpansion}
	p(\vec{x}) = \frac{\euler^{\imag kr}}{r}\sum_{n=0}^\infty \frac{p_n(\vartheta,\varphi)}{r^n},
\end{equation}
which implies that $|p|=\bigoh(r^{-1})$ asymptotically for large $r$. Considering a function which represents this asymptotic property
\begin{equation}
	\Psi(r)=\frac{\euler^{\imag kr}}{r},
\end{equation}
one can observe that the $L^2$ Hermitian inner product does not exist. Indeed, if $\Gamma_1$ is the unit sphere then
\begin{equation*}
	(\Psi,\Psi)_{L^2} = \int_{\Omega^+} \frac{\euler^{\imag kr}}{r}\frac{\euler^{-\imag kr}}{r} \idiff \Omega= 4\PI\int_1^\infty \frac{1}{r^2}r^2\idiff r,
\end{equation*}
which is not finite. The solution to the problem is to introduce weighted norms by defining the inner product
\begin{equation}
	(p,q)_w = \int_{\Omega^+} wp\bar{q}\idiff\Omega,\qquad\text{with}\quad w= \frac{1}{r^2}.
\end{equation}
The following norm may then be induced
\begin{equation}
	\| p\|_{1,w} = \sqrt{(p,p)_w + (\nabla p,\nabla p)_w}
\end{equation}
such that the trial functions satisfy $\| p\|_{1,w}<\infty$. The integrals
\begin{equation}
	\int_{\Omega^+} p\bar{q}\idiff\Omega\quad\text{and}\quad\int_{\Omega^+}\nabla p\cdot\nabla\bar{q}\idiff\Omega
\end{equation}
are well defined if the test functions $q$ are such that
\begin{equation}
	(q,q)_{w^*}<\infty\quad\text{and}\quad (\nabla q,\nabla q)_{w^*}<\infty
\end{equation}
with the inner product
\begin{equation}
	(p,q)_{w^*} = \int_{\Omega^+} w^*p\bar{q}\idiff\Omega,\qquad\text{with}\quad w^*= r^2,
\end{equation}
and the corresponding norm
\begin{equation}
	\| p\|_{1,w^*} = \sqrt{(p,p)_{w^*} + (\nabla p,\nabla p)_{w^*}}.
\end{equation}
Define now the following \textit{weighted Sobolev spaces} for the trial- and test spaces
\begin{equation}
	H_w^1(\Omega^+) = \{p\,:\,\| p\|_{1,w} < \infty\}\quad\text{and}\quad H_{w^*}^1(\Omega^+) =  \{q\,:\,\| q\|_{1,w^*} < \infty\},
\end{equation}
respectively. These definitions will not ensure that all trial function satisfy the Sommerfeld condition. Leis solved this problem in~\cite{Leis1986ibv} by modifying the trial space to be
\begin{equation}
	H_w^{1+}(\Omega^+) = \{p\,:\,\| p\|_{1,w}^+ < \infty\}
\end{equation}
where
\begin{equation}
	\| p\|_{1,w}^+ = \sqrt{\| p\|_{1,w}^2 + \int_{\Omega^+}\left|\pderiv{p}{r} - \imag kp\right|^2\idiff\Omega}.
\end{equation}
For a more detailed discussion of the functional analysis involved in these spaces refer to the book by Ihlenburg~\cite[pp. 41-43]{Ihlenburg1998fea}.

The weak form of the Helmholtz equation may now be found by multiplying \Cref{Eq2:HelmholtzEqn} with a test function and integration over the domain
\begin{equation*}
	\int_{\Omega^+} \left[q\nabla^2 p + k^2 qp\right]\idiff \Omega = 0.
\end{equation*}
Using Greens first identity this can be written as
\begin{equation*}
	-\int_{\Omega^+} \nabla q\cdot\nabla p\idiff\Omega + \int_{\partial\Omega^+} q\nabla p\cdot\vec{n}\idiff\Gamma + k^2\int_{\Omega^+}  qp\idiff \Omega = 0.
\end{equation*}
Thus,
\begin{equation}\label{Eq2:weakformulationHelmholtz}
	\int_{\Omega^+} \nabla q\cdot\nabla p\idiff\Omega -  k^2\int_{\Omega^+}qp\idiff\Omega = \int_{\partial\Omega^+} q g\idiff\Gamma.
\end{equation}
The weak formulation then becomes: 
\begin{equation}
	\text{Find} \quad p\in H_w^{1+}(\Omega^+)\quad\text{such that}\quad B(q,p) = L(q),\qquad \forall q\in H_{w^*}^1(\Omega^+),
\end{equation}
where the bilinear form is given by
\begin{equation*}
	B(q,p) = \int_{\Omega^+} \left[\nabla q\cdot\nabla p-  k^2qp\right]\idiff\Omega
\end{equation*}
and the corresponding linear form is given by
\begin{equation*}
	L(q) = \int_{\Gamma_1} qg\idiff\Gamma.
\end{equation*}

\subsection{Infinite elements}
\label{se:infElems}
In the following, a derivation of the weak formulation for infinite elements using a prolate spheroidal coordinate system is presented (cf.~\cite{Burnett1994atd}). The IEM is typically presented with four infinite element formulations:
\begin{itemize}
	\item Petrov--Galerkin conjugated (PGC)
	\item Petrov--Galerkin unconjugated (PGU)
	\item Bubnov--Galerkin conjugated (BGC)
	\item Bubnov--Galerkin unconjugated (BGU)
\end{itemize}
The Petrov--Galerkin formulations are based on the weighted Sobolev spaces after Leis~\cite{Leis1986ibv}. It turns out that it is possible to create Bubnov--Galerkin formulations as well when the integration in the weak formulation is understood in the sense of the Cauchy principal value (consider~\cite{Burnett1994atd} and~\cite{Gerdes1998tcv} for details). These spaces differ compared to the Petrov--Galerkin counterpart in that the test space and trial space are equal. The difference between the conjugated formulations and the unconjugated formulations is simply conjugations of the test functions in the weak formulation. The accuracy of these formulations has been assessed in the overview in~\cite{Astley2000ief}.

The idea of the IEM is to partition the unbounded domain $\Omega^+$ into $\Omega_{\mathrm{a}}$ and $\Omega_{\mathrm{a}}^+$ separated by an artificial boundary $\Gamma_{\mathrm{a}}$ (cf. \Cref{Fig2:artificialBoundary}). These two domains can then be discretized with finite elements and infinite elements, respectively. The boundary of the scatterer is assumed to be parameterized using 3D NURBS surface patches, such that the domain $\Omega_{\mathrm{a}}$ can be parameterized using 3D NURBS volume patches. Denote by $\calV_h(\Omega_{\mathrm{a}})$, the space spanned by these trivariate NURBS-basis functions. As the 3D NURBS volume representation of $\Omega_{\mathrm{a}}$ reduces to a NURBS surface parametrization at $\Gamma_{\mathrm{a}}$, a natural partition of $\Gamma_{\mathrm{a}}$ into surface elements arises. Denote by $\calV_h(\Gamma_{\mathrm{a}})$, the space spanned by the resulting bivariate basis functions. Consider now the following basis of the \textit{radial shape functions} which is motivated by the Wilcox expansion in \Cref{Eq2:wilcoxExpansion}
\begin{equation}
	\calI_{N, w}^+ = \mathrm{span}\left(\left\{\frac{\euler^{\imag k r}}{r^n}\right\}_{n=1,\dots,N}\right).
\end{equation}
Moreover, define corresponding spaces for the test-space
\begin{equation}
\calI_{N, w^*}^+ = \begin{cases}
	\calI_{N, w}^+  & \text{ for Bubnov--Galerkin formulations}\\
	\mathrm{span}\left(\left\{\frac{\euler^{\imag k r}}{r^{n+2}}\right\}_{n=1,\dots,N}\right) & \text{ for Petrov--Galerkin formulations}.
	\end{cases}
\end{equation}
The trial- and test spaces for the infinite elements can then be defined by
\begin{align}
	\calI_{h, w}^+ &= \calV_h(\Gamma_{\mathrm{a}})\times \calI_{N, w}^+,\\
	\calI_{h, w^*}^+ &= \calV_h(\Gamma_{\mathrm{a}})\times \calI_{N, w^*}^+,
\end{align}
respectively. Finally, the trial- and test spaces for the coupled FEM-IEM can be written as
\begin{align}
	\calF^+_{h,w} &= \left\{p\in H_w^{1+}(\Omega^+);\, p\big|_{\Omega_{\mathrm{a}}}\in \calV_h(\Omega_{\mathrm{a}})\quad\text{and}\quad p\big|_{\Omega_{\mathrm{a}}^+}\in \calI_{h, w}^+\right\},\\
	\calF^+_{h,w^*} &= \left\{q\in H_{w^*}^{1+}(\Omega^+);\, q\big|_{\Omega_{\mathrm{a}}}\in \calV_h(\Omega_{\mathrm{a}})\quad\text{and}\quad q\big|_{\Omega_{\mathrm{a}}^+}\in \calI_{h, w^*}^+\right\},
\end{align}
respectively. Note that $\calF^+_{h,w^*}=\calF^+_{h,w}$ for Bubnov--Galerkin formulations.

\begin{figure}
	\centering
	\includegraphics[scale=1]{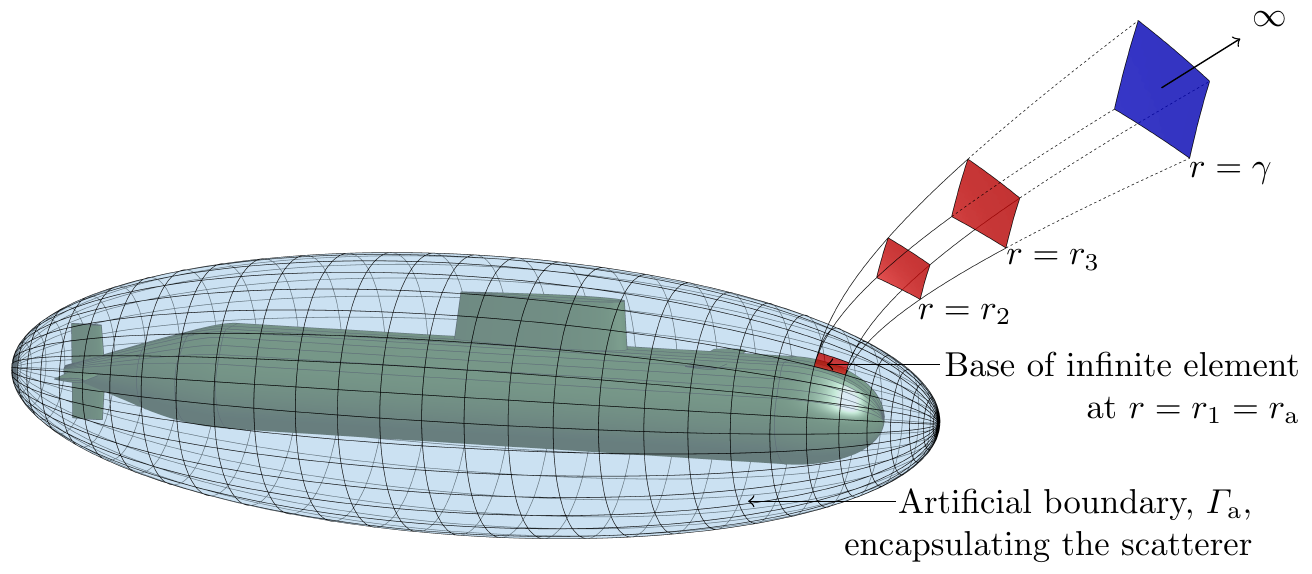}
	\caption{Illustration of an infinite element (with $N=3$) where the radial shape functions have the Kronecker delta property at radii $r_1=r_{\mathrm{a}}$, $r_2=\frac54 r_{\mathrm{a}}$ and $r_3=\frac64 r_{\mathrm{a}}$. The (green) scatterer inside $\Gamma_{\mathrm{a}}$ is the BeTSSi submarine which originates from the BeTSSi workshops~\cite{Gilroy2013bib}. Note that the volume elements discretizing the domain $\Omega_{\mathrm{a}}$ (bounded by the scatterer and the artificial boundary) are not shown here.}
	\label{Fig2:model3_in_waterInf}
\end{figure}

For the unconjugated formulations the Galerkin formulations now takes the form: 
\begin{equation}\label{Eq2:GalerkinFormulationHelmholtz}
	\text{Find}\quad p_h\in \calF^+_{h,w}\quad \text{such that} \quad B_{\mathrm{uc}}(q_h, p_h) = L(q_h),\qquad \forall q_h\in \calF^+_{h,w^*}
\end{equation}
where the bilinear form and linear form are respectively given by
\begin{align}\label{Eq2:infElemntBilinearForm}
	B_{\mathrm{uc}}(q,p) &= \lim_{\gamma\to\infty}\left(\int_{\Omega^\gamma} \left[\nabla q\cdot \nabla p - k^2 qp\right]\idiff\Omega - \int_{S^\gamma} q\partial_n p\idiff\Gamma\right),\\
	L(q) &= \int_{\Gamma_1} qg\idiff\Gamma.\nonumber
\end{align}
Here, $S^\gamma$ is the surface at $r = \gamma$ (and $\Omega^\gamma$ is the domain bounded by $\Gamma_1$ and $S^\gamma$, such that $\lim_{\gamma\to\infty}\Omega^\gamma = \Omega^+$) and the full domain can then be recovered by letting $\gamma\to\infty$ (see \Cref{Fig2:model3_in_waterInf}). Recall that $\partial_n p = \frac{\partial p}{\partial n}= \vec{n}\cdot\nabla p$ where $\vec{n}$ is pointing ``out'' of $\Omega^\gamma$. In the conjugated formulations the test functions $q_h$ are conjugated.

Let $r_{\mathrm{a}}$ be the radius in the prolate spheroidal coordinate system at the artificial boundary $\Gamma_{\mathrm{a}}$. Moreover, let the radial shape functions $\phi$ be defined by
\begin{equation}
	\phi_m(r) = \euler^{\imag k (r-r_{\mathrm{a}})}Q_m\left(\frac{r_{\mathrm{a}}}{r}\right),\quad m = 1,\dots,N
\end{equation}
where 
\begin{equation}\label{Eq2:Q_m}
	Q_m(x) = \sum_{\tilde{m}=1}^N D_{m\tilde{m}} x^{\tilde{m}}
\end{equation}
is a set of polynomial functions defined on the half open interval $(0,1]$. To obtain optimal sparsity of the global matrix, one should choose the polynomials such that $Q_m(1) = \delta_{m1}$, with the Kronecker delta function defined by
\begin{equation}\label{Eq2:Kronecker}
	\delta_{ij} = \begin{cases}
		1 & \text{if}\quad i = j\\
		0 & \text{if}\quad i \neq j
		\end{cases}
\end{equation}
which implies that $\phi_m(r_{\mathrm{a}}) = \delta_{m1}$. In~\cite{Burnett1994atd} Burnett includes the restrictions $\phi_m(r_n) = \delta_{mn}$ with radii $r_n$, $n=1,\dots,N$ (see \Cref{Fig2:model3_in_waterInf}).
Alternatively, one could use the shifted Chebyshev polynomials as done by Shirron and Dey in~\cite{Shirron2002aie}. These polynomials are defined by the three-term recurrence relation
\begin{equation}
	\tilde{T}_{m+1}(x) = 2(2x-1)\tilde{T}_m(x) - \tilde{T}_{m-1}(x)
\end{equation}
for $m\geq 1$ starting with
\begin{equation}
	\tilde{T}_0(x) = 1\quad\text{and}\quad\tilde{T}_1(x) = 2x-1.
\end{equation}
Let
\begin{equation}
	Q_m(x) = \begin{cases}x\left(\tilde{T}_{m-1}(x)-1\right)& m>1\\
				x & m=1.
				\end{cases}
\end{equation}
Then the coefficients $D_{m\tilde{m}}$ in \Cref{Eq2:Q_m} can be collected in the matrix (for $N \leq 6$)
\begin{equation*}
	\vec{D} = \begin{bmatrix}
		1 & 0 & 0 & 0 & 0 & 0\\
		-2 & 2 & 0 & 0 & 0 & 0\\	
		0 & -8 & 8 & 0 & 0 & 0\\
		-2 & 18 & -48 & 32 & 0 & 0\\
		0 & -32 & 160 & -256 & 128 & 0\\
		-2 & 50 & -400 & 1120 & -1280 & 512
	\end{bmatrix}.
\end{equation*}
For the Petrov--Galerkin formulations, a second set of shape functions (for the test space) must be created, namely
\begin{equation}
	\psi_n(r) = \euler^{\imag k (r-r_{\mathrm{a}})}\tilde{Q}_n\left(\frac{r_{\mathrm{a}}}{r}\right),\quad n = 1,\dots,N
\end{equation}
using
\begin{equation}
	\tilde{Q}_n(x) = \sum_{\tilde{n}=1}^N \tilde{D}_{n\tilde{n}} x^{\tilde{n}+2}
\end{equation}
where it is natural to choose $\tilde{D}_{n\tilde{n}}=D_{n\tilde{n}}$. The Bubnov--Galerkin formulations use the same shape functions for the test space, i.e., $\psi_n = \phi_n$.

Alternatively, the polynomials $Q$ can be based upon the Bernstein basis of order $\check{p}=N-1$ by
\begin{equation}
	Q_m(x) = xb_{p-m+1,\check{p}}(x)\qquad m=1,\dots N
\end{equation}
where
\begin{equation}
	b_{i,\check{p}}(x) = \binom{n}{i}(1-x)^{\check{p}-i}x^i = \sum_{j=0}^{\check{p}-i}(-1)^j\binom{\check{p}}{i}\binom{\check{p}-i}{j} x^{i+j},\qquad i=0,\dots,\check{p}.
\end{equation}

For completeness, note that the coefficients for the radial shape functions used by Burnett~\cite{Burnett1994atd} (for the Bubnov--Galerkin formulations) can be found by solving $\vec{D}\vec{B} = \vec{E}$ where
\begin{equation*}
	\vec{B} = \begin{bmatrix}
		x_1 & x_2 & \dots & x_N\\
		x_1^2 & x_2^2 & \dots & x_N^2\\
		\vdots & \vdots & \ddots & \vdots\\
		x_1^N & x_2^N & \dots & x_N^N\\
	\end{bmatrix},\quad \vec{E} = \begin{bmatrix}
		1 	& 	&  & \\
		   	& \euler^{\imag k (r_{\mathrm{a}}-r_2)} & & \\
			&  		& \ddots 	&\\
			&  		& 			& \euler^{\imag k (r_{\mathrm{a}}-r_N)}
	\end{bmatrix}, \quad x_n = \frac{r_{\mathrm{a}}}{r_n}.
\end{equation*}
The coefficients $D_{m\tilde{m}}$ are thus given by $\vec{D} = \vec{E}\vec{B}^{-1}$. For Petrov--Galerkin formulations, the coefficients $\tilde{D}_{n\tilde{n}}$ are found in the same way, but now with the matrix
\begin{equation*}
	\tilde{\vec{B}} = \begin{bmatrix}
		x_1^3 & x_2^3 & \dots & x_N^3\\
		x_1^4 & x_2^4 & \dots & x_N^4\\
		\vdots & \vdots & \ddots & \vdots\\
		x_1^{N+2} & x_2^{N+2} & \dots & x_N^{N+2}\\
	\end{bmatrix}
\end{equation*}
instead of $\vec{B}$. So with the notation presented, these basis functions are based on the Lagrange polynomials with polynomial order\footnote{The usage of a check sign above the polynomial order $p$ is to avoid ambiguity between the polynomial order and the scattered pressure.} $\check{p}=N-1$
\begin{equation}
	l_n(x) = \prod_{\substack{0\leq n\leq \check{p}\\ n\neq m}} \frac{x-x_n}{x_m-x_n},
\end{equation}
since the polynomials $Q_m$ can be written as
\begin{equation*}
	Q_m(x) = \euler^{\imag k(r_{\mathrm{a}}-r_m)} \frac{r_m}{r_{\mathrm{a}}} xl_m(x)
\end{equation*}
such that
\begin{equation*}
	\phi_m(r) = \euler^{\imag k(r-r_m)}\frac{r_m}{r}l_m\left(\frac{r_{\mathrm{a}}}{r}\right).
\end{equation*}
The radial shape functions in the test space for the Petrov--Galerkin formulations take the form
\begin{equation*}
	\psi_n(r) = \euler^{\imag k(r-r_n)}\left(\frac{r_n}{r}\right)^3 l_n\left(\frac{r_{\mathrm{a}}}{r}\right).
\end{equation*}
As all these sets of basis functions span the same space, they should only affect the conditioning of the system. Note that the sets of basis functions are identical for $N=1$.

The trial- and test functions now take the form 
\begin{equation}\label{Eq2:p_h}
	p_h(\vec{x}) = \begin{dcases}
		\sum_{J\in\vec{\kappa}_{\mathrm{a}}} \sum_{m=1}^N d_{m,J}\phi_m(r) R_J(\xi,\eta,\zeta)\big|_{\Gamma_{\mathrm{a}}} & \vec{x}\in\Omega_{\mathrm{a}}^+\\
		\sum_{J\in\vec{\kappa}} d_{1,J} R_J(\xi,\eta,\zeta) & \vec{x}\in\Omega_{\mathrm{a}}
	\end{dcases}
\end{equation}
and
\begin{equation}\label{Eq2:q_h}
	q_h(\vec{x}) = \begin{dcases}
		\sum_{I\in\vec{\kappa}_{\mathrm{a}}} \sum_{n=1}^N c_{n,I}\psi_n(r) R_I(\xi,\eta,\zeta)\big|_{\Gamma_{\mathrm{a}}} & \vec{x}\in\Omega_{\mathrm{a}}^+\\
		\sum_{I\in\vec{\kappa}} c_{1,I} R_I(\xi,\eta,\zeta) & \vec{x}\in\Omega_{\mathrm{a}},
	\end{dcases}
\end{equation}
respectively. Here, $\vec{\kappa}$ is the collection of the global indices of the NURBS basis functions and $\vec{\kappa}_{\mathrm{a}}$ the corresponding indices of the non-zero NURBS function at the surface $\Gamma_{\mathrm{a}}$. Moreover, $R_I(\xi,\eta,\zeta)$ is the set of NURBS basis functions. The system of equations will now be obtained by inserting the functions in \Cref{Eq2:p_h} and \Cref{Eq2:q_h} into the bilinear form (or sesquilinear form for the BGC and PGC formulations, i.e. the bilinear form with conjugated test functions). 

Before the insertion, it is advantageous to split the bilinear form in \Cref{Eq2:infElemntBilinearForm} as
\begin{equation}\label{Eq2:SplittingOfB_uc}
	B_{\mathrm{uc}}(q,p) = B_{\mathrm{a}}(q,p) + B_{\mathrm{uc},\mathrm{a}}^+(q,p)
\end{equation}
where
\begin{align}
	B_{\mathrm{a}}(q,p) &= \int_{\Omega_{\mathrm{a}}} \left[\nabla q\cdot \nabla p - k^2 qp\right]\idiff\Omega\nonumber\\
	B_{\mathrm{uc},\mathrm{a}}^+(q,p) &= \lim_{\gamma\to\infty}\left(\int_{\Omega_{\mathrm{a}}^\gamma} \left[\nabla q\cdot \nabla p - k^2 qp\right]\idiff\Omega - \int_{S^\gamma} q\partial_n p\idiff\Gamma\right).\label{Eq2:B_uc_a}	
\end{align}
Insertion of \Cref{Eq2:p_h} and \Cref{Eq2:q_h} into \Cref{Eq2:GalerkinFormulationHelmholtz} (using the splitting in \Cref{Eq2:SplittingOfB_uc}) results in the following system of equations
\begin{equation}
	(\vec{A}_{\mathrm{a}} + \vec{A}_{\mathrm{uc},\mathrm{a}}^+)\vec{d} = \vec{F}
\end{equation}
with components
\begin{alignat*}{3}
	& &&\vec{A}_{\mathrm{a}}[I,J] = B_{\mathrm{a}}(R_I,R_J)\qquad I, && J = 1,\dots,|\boldsymbol\kappa|\\
	& &&\vec{F}[I] = L(R_I)\qquad && I = 1,\dots,|\boldsymbol\kappa|\\
	& &&\vec{d}[J] = d_{1,J}\qquad && J = 1,\dots,|\boldsymbol\kappa|
\end{alignat*}
and
\begin{align*}
	\vec{A}_{\mathrm{uc},\mathrm{a}}^+[\tilde{I},\tilde{J}] &= B_{\mathrm{uc},\mathrm{a}}^+(R_I\psi_n,R_J\phi_m)\\
	\vec{d}[\tilde{J}] &= d_{m,J}
\end{align*}
where $I = \boldsymbol\kappa_{\mathrm{a}}[\tilde{i}]$ and $J = \boldsymbol\kappa_{\mathrm{a}}[\tilde{j}]$ for $\tilde{i},\tilde{j}=1,\dots,|\boldsymbol\kappa_{\mathrm{a}}|$ and $m,n=1,\dots,N$, and
\begin{align*}
	\tilde{I} &= \begin{cases}\boldsymbol\kappa_{\mathrm{a}}[\tilde{i}] & n = 1\\
	|\boldsymbol\kappa|+(n-2)|\boldsymbol\kappa_{\mathrm{a}}|+\tilde{i} & n>1\end{cases}\\
	\tilde{J} &= \begin{cases}\boldsymbol\kappa_{\mathrm{a}}[\tilde{j}] & m = 1\\
	|\boldsymbol\kappa|+(m-2)|\boldsymbol\kappa_{\mathrm{a}}|+\tilde{j} & m>1.\end{cases}
\end{align*}
Note that $\vec{A}_{\mathrm{a}}$ and $\vec{F}$ are independent of the IEM and that there are $|\boldsymbol\kappa| +|\boldsymbol\kappa_{\mathrm{a}}|(N-1)$ linear equations. The matrices are assembled as in the classical FEM. That is, instead of looping through the indices, one loops through the elements. A formula for $B_{\mathrm{uc},\mathrm{a}}^+(R_I\psi_n,R_J\phi_m)$ for the Petrov Galerkin formulation is derived in \Cref{Sec2:AppendixDerivationOfBilinearForm} and the final bilinear form is given in \Cref{Eq2:finalBilinearFormB_uc_a}. The final formulas for the other three formulations are also added in this appendix.

\subsection{Far field pattern}
The problem is solved inside an artificial boundary, computing the so-called near field. However, the far field is also often of interest. To solve this issue, one uses the integral solution given by\footnote{For the conjugated formulations one may also compute the far field using the radial shape functions in the infinite elements, but for the unconjugated formulations it is mentioned in~\cite[p. 137]{Shirron1998aco} that the expansion does not converge in the far field, such that it must be computed by other means.} (cf.~\cite[Theorem 2.21]{Chandler_Wilde2012nab})
\begin{equation}\label{Eq2:KirchhoffIntegral}
	p(\vec{x}) = \int_{\Gamma_1}\left[ p(\vec{y})\pderiv{\Phi_k(\vec{x},\vec{y})}{n(\vec{y})} - \Phi_k(\vec{x},\vec{y})\pderiv{p(\vec{y})}{n(\vec{y})}\right]\idiff \Gamma(\vec{y})
\end{equation}
where $\vec{y}$ is a point on the surface $\Gamma_1$, $n$ lies on $\Gamma_1$ pointing ``into'' $\Omega^+$ at $\vec{y}$ and $\Phi_k$ is the free space Green's function for the Helmholtz equation in \Cref{Eq2:HelmholtzEqn} given (in 3D) by
\begin{equation}\label{Eq2:FreeSpaceGrensFunction}
	\Phi_k(\vec{x},\vec{y}) = \frac{\euler^{\imag kR}}{4\PI R},\quad\text{where}\quad R = |\vec{x} - \vec{y}|.
\end{equation}
The derivative of both Green's function and the numerical solution for the pressure is therefore needed
\begin{equation}
	\pderiv{\Phi_k(\vec{x},\vec{y})}{n(\vec{y})} = \frac{\Phi_k(\vec{x},\vec{y})}{R}(\imag kR-1)\pderiv{R}{n(\vec{y})},\quad\text{where}\quad\pderiv{R}{n(\vec{y})} = -\frac{(\vec{x}-\vec{y})\cdot\vec{n}(\vec{y})}{R}.
\end{equation}
Note that for sound-hard scattering (where $\partial_n(p+p_{\mathrm{inc}}) = 0$) the values for $\partial_n p$ are known at the boundary $\Gamma_1$ (given by \Cref{Eq2:HelmholtzEqnNeumannCond}). To use the exact expression for the derivative seems to give better results, and is for this reason used in the sound-hard scattering cases when computing the field outside the artificial boundary.

The \textit{far field pattern} for the scattered pressure $p$, is now defined by
\begin{equation}\label{Eq2:farfield}
	p_0(\hat{\vec{x}}) =  \lim_{r\to\infty} r \euler^{-\imag k r}p(r\hat{\vec{x}}),
\end{equation}
with $r = |\vec{x}|$ and $\hat{\vec{x}} = \vec{x}/|\vec{x}|$. Using the limits
\begin{equation}\label{Eq2:Phi_k_limits}
	\lim_{r\to\infty} r\euler^{-\imag k r}\Phi_k(r\hat{\vec{x}},\vec{y}) = \frac{1}{4\PI}\euler^{-\imag k \hat{\vec{x}}\cdot\vec{y}}\quad\text{and}\quad 
	\lim_{r\to\infty} r\euler^{-\imag k r}\pderiv{\Phi_k(r\hat{\vec{x}},\vec{y})}{n(\vec{y})} = -\frac{\imag k}{4\PI}\euler^{-\imag k \hat{\vec{x}}\cdot\vec{y}}\hat{\vec{x}}\cdot\vec{n}(\vec{y})
\end{equation}
the formula in \Cref{Eq2:KirchhoffIntegral} simplifies in the far field to (cf.~\cite[p. 32]{Ihlenburg1998fea})
\begin{equation}\label{Eq2:KirchhoffIntegralFarField}
	p_0(\hat{\vec{x}}) = -\frac{1}{4\PI}\int_{\Gamma_1}\left[ \imag k p(\vec{y})\hat{\vec{x}}\cdot\vec{n}(\vec{y}) + \pderiv{p(\vec{y})}{n(\vec{y})}\right]\euler^{-\imag k \hat{\vec{x}}\cdot\vec{y}}\idiff \Gamma(\vec{y}).
\end{equation}
From the far field pattern, the \textit{target strength}, $\TS$, can be computed. It is defined by
\begin{equation}\label{Eq2:TS}
	\TS = 20\log_{10}\left(\frac{|p_0(\hat{\vec{x}})|}{|P_{\mathrm{inc}}|}\right)
\end{equation}
where $P_{\mathrm{inc}}$ is the amplitude of the incident wave at the geometric center of the scatterer (i.e. the origin). Note that $\TS$ is independent of $P_{\mathrm{inc}}$, which is a result of the linear dependency of the amplitude of the incident wave in scattering problems (i.e. doubling the amplitude of the incident wave will double the amplitude of the scattered wave).

%% file: contents/coupledFluidStructure.tex
\section{Acoustic-structure interaction}
\label{Sec2:coupledFluidStruct}
In~\cite[pp. 13-14]{Ihlenburg1998fea} Ihlenburg briefly derives the governing equations for the ASI problem. Building upon this the formulas are generalized to include an interior fluid domain $\Omega^-$. The pressure in the exterior and interior fluid domain are now denoted by $p_1$ and $p_2$ (see \Cref{Fig2:artificialBoundary}).
\begin{alignat}{3}
	\nabla^2 p_1 + k_1^2 p_1 &= 0 	&&\text{in}\quad \Omega^+\\
	\pderiv{p(\vec{x},\omega)}{r}-\imag k p(\vec{x},\omega) &= o\left(r^{-1}\right) &&\text{with}\quad r=|\vec{x}|\label{Eq2:sommerfeldCond2}\\
	\rho_{\mathrm{f},1} \omega^2 u_i n_i - \pderiv{p_1}{n} &= \pderiv{p_{\mathrm{inc}}}{n}\qquad &&\text{on}\quad \Gamma_1\label{Eq2:coupling2}\\
	\sigma_{ij}n_i n_j + p_1 &= -p_{\mathrm{inc}}\qquad &&\text{on}\quad \Gamma_1\label{Eq2:coupling1}\\
	\sigma_{ij,j} + \omega^2 \rho_{\mathrm{s}} u_i &= 0 \qquad &&\text{in}\quad \Omega_{\mathrm{s}}\label{Eq2:strongFormLinEl}\\
	\rho_{\mathrm{f},2} \omega^2 u_i n_i - \pderiv{p_2}{n} &= 0 \qquad &&\text{on}\quad \Gamma_2\label{Eq2:coupling2_2}\\
	\sigma_{ij}n_i n_j + p_2 &= 0\qquad &&\text{on}\quad \Gamma_2\label{Eq2:coupling1_2}\\
	\nabla^2 p_2 + k_2^2 p_2 &= 0 	&&\text{in}\quad \Omega^-.
\end{alignat} 
The first two equations represent the Helmholtz equation and Sommerfeld conditions, respectively, for the exterior domain. The wave numbers in the exterior and interior fluid domain are denoted by $k_1$ and $k_2$. The elasticity equation in \Cref{Eq2:strongFormLinEl} comes from momentum conservation (Newton's second law), while \Cref{Eq2:coupling2,Eq2:coupling1,Eq2:coupling2_2,Eq2:coupling1_2} represent the coupling equations and come from the continuity requirement of the displacement and pressures at the boundaries $\Gamma_m$. The final formula is simply the Helmholtz equation for the internal fluid domain. The function $p_{\mathrm{inc}}$ represents the incident plane wave in \Cref{Eq2:p_inc} (in the exterior domain). The mass densities of the solid and the fluid are denoted by $\rho_{\mathrm{s}}$ and $\rho_{\mathrm{f}}$, respectively, and $\sigma_{ij}(\vec{u})$ represents the stress components as a function of the displacement $\vec{u}=u_i\vec{e}_i$ in the solid. 

For the domain of the scatterer, $\Omega_{\mathrm{s}}$, it can be shown that the following weak formulation is obtained from the strong form in \Cref{Eq2:strongFormLinEl} (see for example~\cite{Ihlenburg1998fea})
\begin{equation}\label{Eq2:intermediateStepFSI}
	\int_{\Omega_{\mathrm{s}}} \left[v_{i,j}\sigma_{ij} - \rho_{\mathrm{s}}\omega^2 u_i\bar{v}_i\right]\idiff\Omega = \int_{\Gamma_1} v_i(\sigma_{ij} n_j)\idiff\Gamma + \int_{\Gamma_2} v_i(\sigma_{ij} n_j)\idiff\Gamma.
\end{equation}
where the normal vectors point out of $\Omega_{\mathrm{s}}$. The integrands on the right-hand side may be rewritten using \Cref{Eq2:coupling1,Eq2:coupling1_2} in the following way. Consider a point $\vec{P}$ on $\Gamma_1$ or $\Gamma_2$, with normal vector $\vec{n}=n_i\vec{e}_i$. Let $T_i$ be the components (in Cartesian coordinates) of the exterior traction vector $\vec{T}$. That is to say, $T_i = \sigma_{ij} n_j$. One can then create a local orthogonal coordinate system at this point with unit vectors $\vec{e}_\perp$, $\vec{e}_{\|_1}$ and $\vec{e}_{\|_2}$, where the latter two vectors represent basis vectors for the tangential plane of the surface at $\vec{P}$ (and $\vec{e}_\perp$ represents the normal unit vector on this plane at $\vec{P}$). 

As the scalar product is invariant to orthogonal transformations, the following holds
\begin{equation*}
	T_i v_i = T_x v_x + T_y v_y + T_z v_z = T_{\perp} v_{\perp} + T_{\|_1} v_{\|_1} + T_{\|_2} v_{\|_2}.
\end{equation*}
Since the acoustic pressure from the fluid only exerts forces normal to the surfaces $\Gamma_1$ and $\Gamma_2$, the static equilibrium conditions for the traction at $\vec{P}$ are given by
\begin{equation*}
	T_{\|_1}=0,\qquad T_{\|_2}=0,\quad\text{and}\quad T_{\perp} = -p_{\mathrm{tot},m},
\end{equation*}
where the total pressure is given by
\begin{equation*}
	p_{\mathrm{tot},m}= \begin{cases} p_{\mathrm{inc}} + p_1 & m = 1\\
	p_2 & m = 2.\end{cases}
\end{equation*}
The scalar product may therefore be written as
\begin{equation*}
	T_i v_i = -p_{\mathrm{tot},m} v_{\perp} = -p_{\mathrm{tot},m} v_i n_i.
\end{equation*}
\Cref{Eq2:intermediateStepFSI} can thus be rewritten as
\begin{equation}\label{Eq2:FSIeq1}
	\int_{\Omega_{\mathrm{s}}} \left[v_{i,j}\sigma_{ij} - \rho_{\mathrm{s}}\omega^2 u_i v_i\right]\idiff\Omega = -\int_{\Gamma_1} (p_{\mathrm{inc}} + p_1) v_i n_i\idiff\Gamma-\int_{\Gamma_2} p_2 v_i n_i\idiff\Gamma.
\end{equation}
Moreover, from \Cref{Eq2:weakformulationHelmholtz} one obtains
\begin{equation*}
	\int_{\Omega^+} \left[\nabla q_1\cdot\nabla p_1 - k_1^2q_1p_1\right]\idiff\Omega = -\int_{\Gamma_1} q_1 \pderiv{p_1}{n}\idiff\Gamma
\end{equation*}
and
\begin{equation*}
	\int_{\Omega^-} \left[\nabla q_2\cdot\nabla p_2 - k_2^2q_2p_2\right]\idiff\Omega = -\int_{\Gamma_2} q_2 \pderiv{p_2}{n}\idiff\Gamma
\end{equation*}
where the sign of the right-hand side must be changed in order to get a normal vector that points out of $\Omega_{\mathrm{s}}$. Using now \Cref{Eq2:coupling2,Eq2:coupling2_2}
\begin{equation}\label{Eq2:FSIeq2}
	\frac{1}{\rho_{\mathrm{f},1} \omega^2}\int_{\Omega^+} \left[\nabla q_1\cdot\nabla p_1 -  k_1^2 q_1p_1\right]\idiff\Omega = -\int_{\Gamma_1} q_1\left(u_i n_i -\frac{1}{\rho_{\mathrm{f},1} \omega^2}\pderiv{p_{\mathrm{inc}}}{n}\right)\idiff\Gamma
\end{equation}
and
\begin{equation}\label{Eq2:FSIeq3}
	\frac{1}{\rho_{\mathrm{f},2} \omega^2}\int_{\Omega^-} \left[\nabla q_2\cdot\nabla p_2 -  k_2^2 q_2p_2\right]\idiff\Omega = -\int_{\Gamma_2} q_2 u_i n_i\idiff\Gamma.
\end{equation}
Adding \Cref{Eq2:FSIeq1,Eq2:FSIeq2,Eq2:FSIeq3}
\begin{align*} %\label{Eq2:FSIbilinearForm}
\begin{split}
	&\frac{1}{\rho_{\mathrm{f},1} \omega^2}\int_{\Omega^+} \left[\nabla q_1\cdot\nabla p_1 -  k_1^2 q_1p_1\right]\idiff\Omega + \int_{\Gamma_1}\left[q_1 u_i n_i + p_1 v_i n_i\right]\idiff\Gamma + \int_{\Omega_{\mathrm{s}}} \left[v_{i,j}\sigma_{ij} - \rho_{\mathrm{s}}\omega^2 u_i v_i\right]\idiff\Omega \\
	+&\frac{1}{\rho_{\mathrm{f},2} \omega^2}\int_{\Omega^-} \left[\nabla q_2\cdot\nabla p_2 -  k_2^2 q_2p_2\right]\idiff\Omega +\int_{\Gamma_2}\left[q_2 u_i n_i + p_2 v_i n_i\right]\idiff\Gamma = \int_{\Gamma_1} \left[\frac{1}{\rho_{\mathrm{f},1} \omega^2}q_1\pderiv{p_{\mathrm{inc}}}{n} - p_{\mathrm{inc}} v_i n_i\right]\idiff\Gamma
\end{split}
\end{align*}
where $\vec{n}=\{n_1,n_2,n_3\}$ points outwards from the solid. Defining the Sobolev spaces $\bm{\calH}_w = \bm{\calS}\times H_w^{1+}(\Omega^+) \times H^1(\Omega^-)$ and $\bm{\calH}_{w^*} = \bm{\calS}\times H_{w^*}^1(\Omega^+) \times H^1(\Omega^-)$ where $\bm{\calS} = \{\vec{u}: u_i\in H^1(\Omega_{\mathrm{s}})\}$, the weak formulation for the ASI problem then becomes (with the notation $U=\{\vec{u},p_1,p_2\}$ and $V = \{\vec{v},q_1,q_2\}$): 
\begin{equation}
	\text{Find}\quad U\in\bm{\calH}_w\quad \text{such that} \quad B_{\mathrm{ASI}}(V,U) = L_{\mathrm{ASI}}(V),\quad \forall V\in\bm{\calH}_{w^*}
\end{equation}
where
\begin{align*}
	B_{\mathrm{ASI}}(V,U) &= \frac{1}{\rho_{\mathrm{f},1} \omega^2}\int_{\Omega^+} \left[\nabla q_1\cdot\nabla p_1 -  k_1^2 q_1p_1\right]\idiff\Omega + \int_{\Gamma_1}\left[q_1 u_i n_i + p_1 v_i n_i\right]\idiff\Gamma + \int_{\Omega_{\mathrm{s}}} \left[v_{i,j}\sigma_{ij} - \rho_{\mathrm{s}}\omega^2 u_i v_i\right]\idiff\Omega \\
	&{\hskip1em\relax}+\frac{1}{\rho_{\mathrm{f},2} \omega^2}\int_{\Omega^-} \left[\nabla q_2\cdot\nabla p_2 -  k_2^2 q_2p_2\right]\idiff\Omega +\int_{\Gamma_2}\left[q_2 u_i n_i + p_2 v_i n_i\right]\idiff\Gamma
\end{align*}
and
\begin{equation*}
	L_{\mathrm{ASI}}(V) = \int_{\Gamma_1} \left[\frac{1}{\rho_{\mathrm{f},1} \omega^2}q_1\pderiv{p_{\mathrm{inc}}}{n} - p_{\mathrm{inc}} v_i n_i\right]\idiff\Gamma.
\end{equation*}
Let $\bm{\calS}_h=\{\vec{u}: u_i\in \calV(\Omega_{\mathrm{s}})\}\subset\bm{\calS}$ where $\calV(\Omega_{\mathrm{s}})$ is the space spanned by the NURBS basis functions used to parameterize $\calV(\Omega_{\mathrm{s}})$, and correspondingly for $\calF^-_h=\{p_2: p_2\in \calV(\Omega^-)\}\subset H^1(\Omega^-)$. Moreover, define the spaces $\bm{\calH}_{h,w} = \bm{\calS}_h \times \calF^+_{h,w} \times \calF^-_h$ and $\bm{\calH}_{h,w^*} = \bm{\calS}_h\times \calF^+_{h,w^*} \times \calF^-_h$. The Galerkin formulation for the ASI problem then becomes: 
\begin{equation}
	\text{Find}\quad U_h\in\bm{\calH}_{h,w}\quad \text{such that} \quad B_{\mathrm{ASI}}(V_h,U_h) = L_{\mathrm{ASI}}(V_h),\quad \forall V_h\in\bm{\calH}_{h,w^*}.
\end{equation}
%The $H^1$-norm for the three domains $\Omega = \Omega^+\cup \Omega_{\mathrm{s}}\cup\Omega^-$ is defined by
%\begin{equation}\label{Eq2:H1Norm}
%	\Hnorm{U}{\Omega} = \sqrt{\Hnorm{p_1}{\Omega_{\mathrm{a}}}^2 + \Hnorm{\vec{u}}{\Omega_{\mathrm{s}}}^2 + \Hnorm{p_2}{\Omega^-}^2}.
%\end{equation}
%where
%\begin{equation*}
%	\Hnorm{p}{\Omega} = \sqrt{\int_{\Omega} |p|^2 + \sum_{i=1}^3\left|\pderiv{p}{x_i}\right|^2\idiff\Omega}\quad\text{and}\quad
%	\Hnorm{\vec{u}}{\Omega} = \sqrt{\int_\Omega \sum_{i=1}^3\left|u_i\right|^2 + \sum_{i=1}^3 \sum_{j=1}^3\left|\pderiv{u_i}{x_j}\right|^2\idiff\Omega}.
%\end{equation*}
As the bilinear forms treated in this work are not $V$-elliptic \cite[p. 46]{Ihlenburg1998fea}, they do not induce a well defined energy-norm. For this reason, the energy norm for the fluid domains $\Omega_{\mathrm{a}}$ are defined by
\begin{equation}\label{Eq2:energyNormFluids}
	\energyNorm{p_1}{\Omega_{\mathrm{a}}} = \sqrt{\int_{\Omega_{\mathrm{a}}} \left|\nabla p_1\right|^2 + k_1^2|p_1|^2 \idiff\Omega}\quad\text{and}\quad\energyNorm{p_2}{\Omega^-} = \sqrt{\int_{\Omega^-} \left|\nabla p_2\right|^2 + k_2^2|p_2|^2 \idiff\Omega}
\end{equation}
and for the solid domain (using Einstein summation convention)
\begin{equation}
	\energyNorm{\vec{u}}{\Omega_{\mathrm{s}}} = \sqrt{\int_{\Omega_{\mathrm{s}}} u_{(i,j)}c_{ijkl}\bar{u}_{(k,l)} + \rho_{\mathrm{s}}\omega^2|\vec{u}|^2\idiff\Omega}
\end{equation}
where
\begin{equation*}
	u_{(i,j)} = \frac{1}{2}\left(\pderiv{u_i}{x_j} + \pderiv{u_j}{x_i}\right)
\end{equation*}
and elastic coefficients expressed in terms of Young's modulus, $E$, and the Poisson's ratio, $\nu$, as~\cite[p. 110]{Cottrell2009iat}
\begin{equation*}
	c_{ijkl} = \frac{\nu E}{(1+\nu)(1-2\nu)}\delta_{ij}\delta_{kl} +\frac{E}{2(1+\nu)}(\delta_{ik}\delta_{jl} + \delta_{il}\delta_{jk}).
\end{equation*}
The energy norm for the coupled problem with $\Omega = \Omega_{\mathrm{a}}\cup \Omega_{\mathrm{s}}\cup\Omega^-$ is then defined by
\begin{align}\label{Eq2:energyNorm}
	\energyNorm{U}{\Omega} = \sqrt{\frac{1}{\rho_{\mathrm{f},1}\omega^2}\energyNorm{p_1}{\Omega_{\mathrm{a}}}^2 + \energyNorm{\vec{u}}{\Omega_{\mathrm{s}}}^2 + \frac{1}{\rho_{\mathrm{f},2}\omega^2}\energyNorm{p_2}{\Omega^-}^2}.
\end{align}
As the unconjugated formulations do not converge in the far field, the norm in the exterior domain is taken over the $\Omega_{\mathrm{a}}$ instead of $\Omega^+$. 

%% file: contents/resultsAndDiscussion.tex
\section{Numerical examples} 
\label{Sec2:resultsDisc}
Rigid scattering on a sphere and elastic scattering on a spherical shell are investigated in the following. These problems possess analytic solutions~\cite{Venas2019e3s} and are for this reason often used to verify numerical methods in acoustic scattering, e.g.~\cite{Gerdes1996so3,Ihlenburg1998fea,Simpson2014aib,Gerdes1998tcv,Gerdes1999otp,Coox2017aii}.  The mock shell is analyzed to investigate the infinite element formulations, and we end this section by analyzing a simplified submarine benchmark.

In this work, the test setting is chosen so that the present approach can be compared to other methods. In particular, the scattering on a rigid sphere example found in~\cite{Simpson2014aib} and the scattering on a spherical shell used in~\cite{Ihlenburg1998fea} are addressed. The latter problem will be investigated in depth and we shall build upon this problem to include both rigid scattering and scattering with full ASI on both sides of the shell.

The direction of the incident wave is along the $x$-axis while the symmetry of the parametrization of the domain is around the $z$-axis (to avoid exploitation of the symmetry of the problems).

We define the \textit{SAV index} by
\begin{equation}
	I_{\mathrm{SAV}} = \frac{L_{\Gamma_{\mathrm{a}}}}{2}\frac{|\Gamma_1|}{|\Omega_{\mathrm{a}}|}
\end{equation}
where $L_{\Gamma_{\mathrm{a}}}$ is the characteristic length of the artificial boundary, $|\Gamma_1|$ is the surface area of the scatterer and $|\Omega_{\mathrm{a}}|$ is the volume of the discretized fluid between $\Gamma_1$ and $\Gamma_{\mathrm{a}}$. The SAV index is based on a scaled surface-area-to-volume ratio (SA/V) such that the domain of computation is fitted in a unit sphere. It can be thought of as an efficiency index for the IEM compared to BEM, as problems with low $I_{\mathrm{SAV}}$ will be more suited for BEM, while high values of $I_{\mathrm{SAV}}$ will be more suited for IEM. If we for the sphere example place the artificial boundary, $\Gamma_{\mathrm{a}}$, at $r_{\mathrm{a}}=sR_0$, where $R_0$ is the outer radius of the scatterer, then the SAV index is given by
\begin{equation}
	I_{\mathrm{SAV}} = \frac{3s}{s^3-1}.
\end{equation}
The IEM is optimal for the sphere problem in the sense that the SAV index can be arbitrarily large. In fact, the infinite elements can be attached directly onto the scatterer (such that ${I_{\mathrm{SAV}}=\infty}$) as done in~\cite{Shirron2002aie}. This, however, is not the case for more complex geometries. A typical SAV index for submarines like the one depicted in \Cref{Fig2:model3_in_waterInf} is approximately~5, so by choosing $s>1$, the SAV index can be adjusted for a fairer comparison with methods like BEM. In the numerical experiments on spherical shells we use $s=\frac{32+\PI}{32-\PI}\approx 1.2$ (such that the aspect ratio of the elements in the tensor product meshes are minimal), resulting in $I_{\mathrm{SAV}}\approx 4.5$.
\begin{figure}
	\centering
	\begin{subfigure}{0.3\textwidth}
		\centering
		\includegraphics[width=0.8\textwidth]{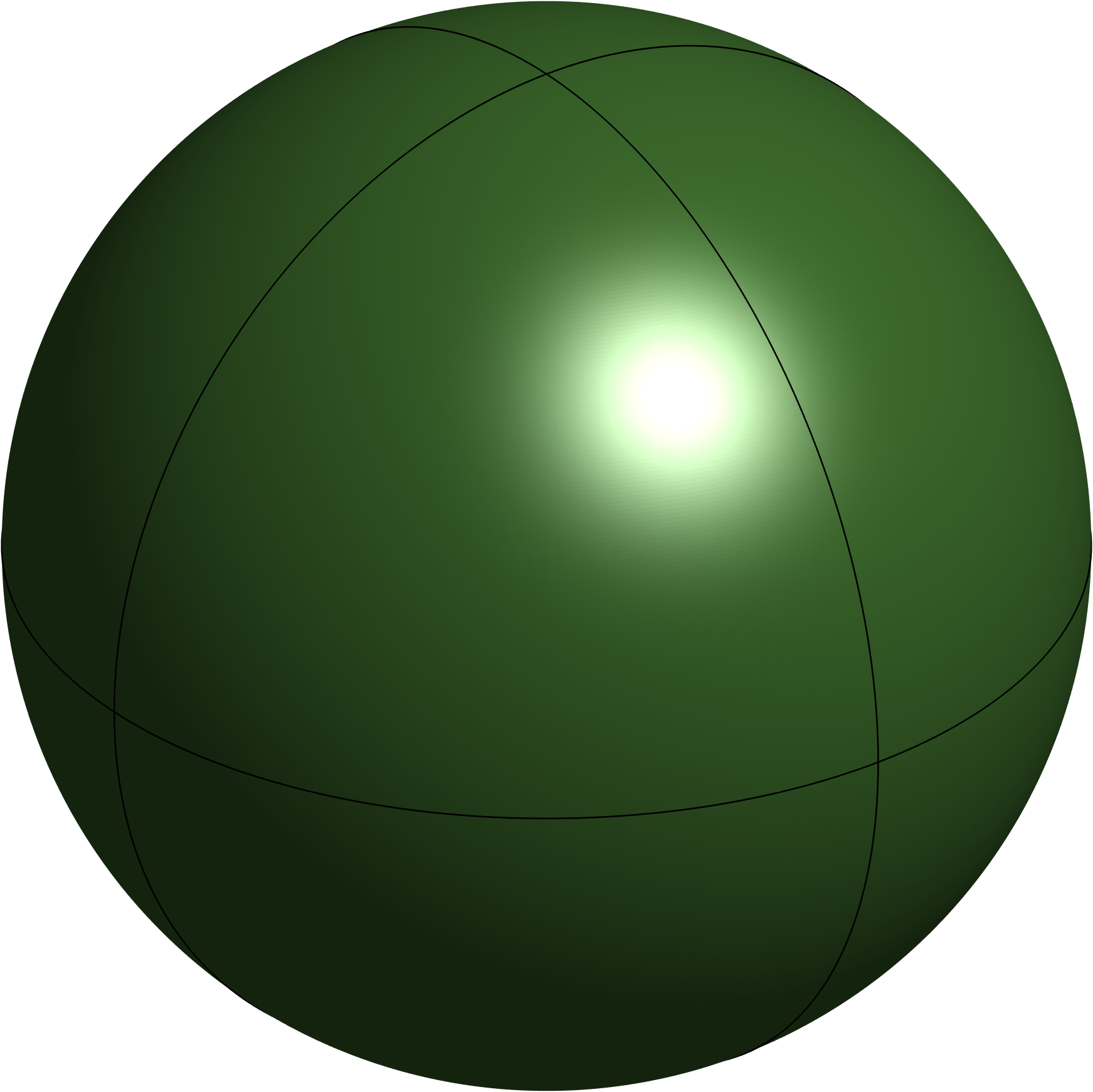}
		\caption{Mesh ${\cal M}_{1,\check{p},\check{k}}^{\textsc{iga}}$}
		\label{Fig2:SphericalShellMeshes1}
    \end{subfigure}
    ~
	\begin{subfigure}{0.3\textwidth}
		\centering
		\includegraphics[width=0.8\textwidth]{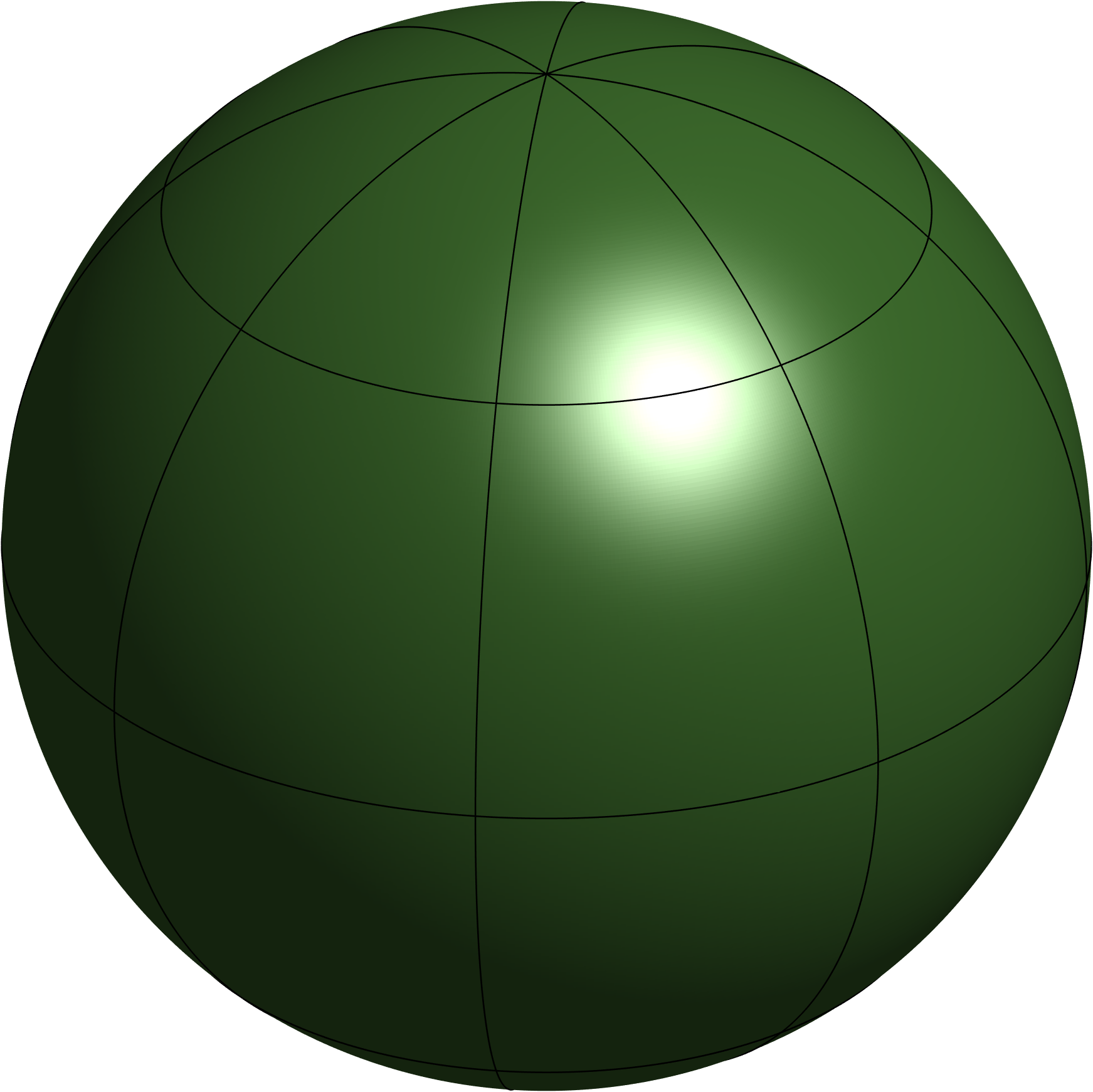}
		\caption{Mesh ${\cal M}_{2,\check{p},\check{k}}^{\textsc{iga}}$}
		\label{Fig2:SphericalShellMeshes2}
    \end{subfigure}
    ~
	\begin{subfigure}{0.3\textwidth}
		\centering
		\includegraphics[width=0.8\textwidth]{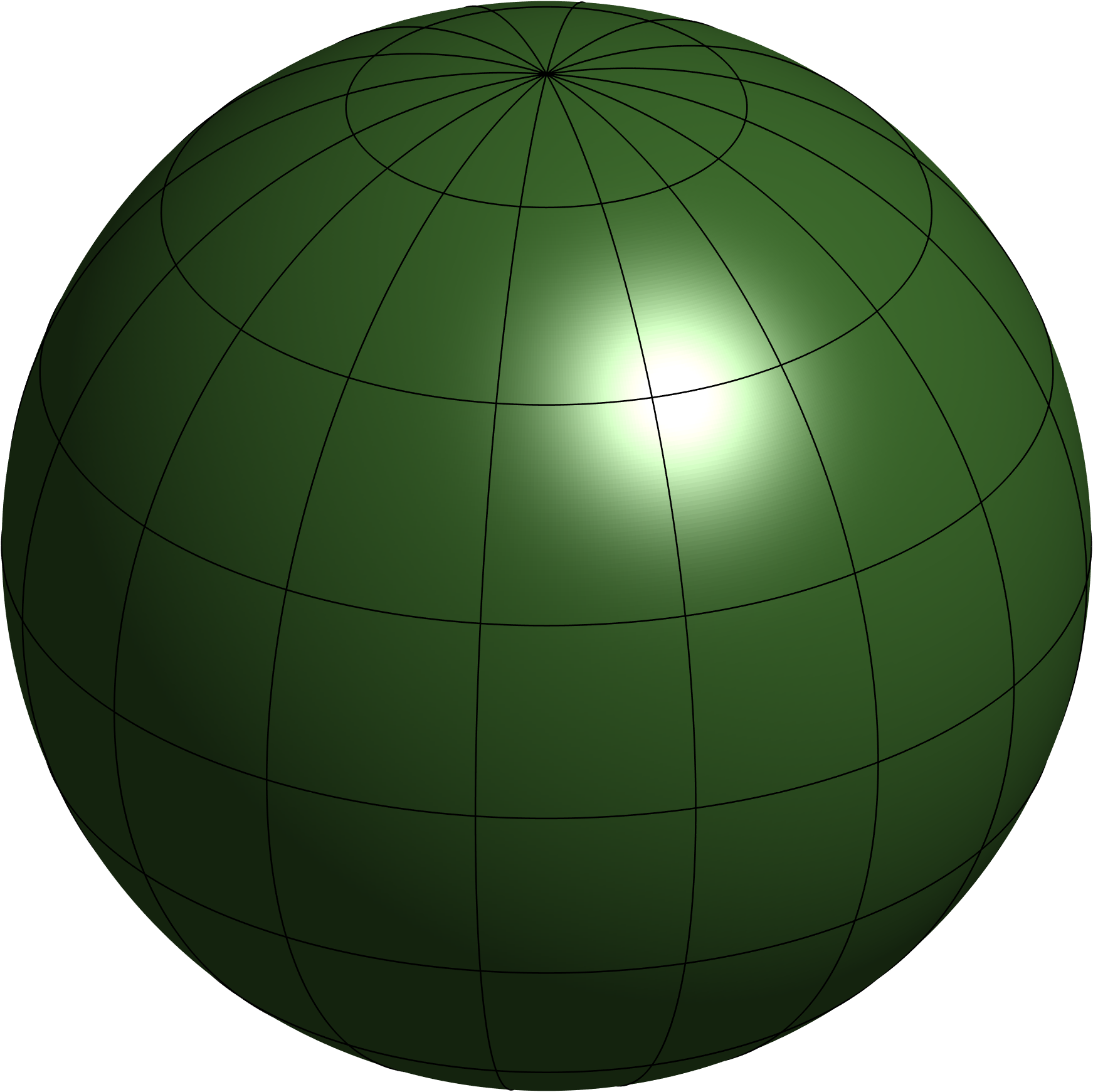}
		\caption{Mesh ${\cal M}_{3,\check{p},\check{k}}^{\textsc{iga}}$}
		\label{Fig2:SphericalShellMeshes3}
    \end{subfigure}
	\caption{\textbf{Numerical examples}: Illustration of the first three meshes, using two successive refinements from the coarse mesh ${\cal M}_{1,\check{p},\check{k}}^{\textsc{iga}}$.}
	\label{Fig2:SphericalShellMeshes}
\end{figure}

The meshes will be generated from a standard discretization of a sphere using NURBS as seen in \Cref{Fig2:SphericalShellMeshes}. We shall denote by ${\cal M}_{m,\check{p},\check{k}}^{\textsc{iga}}$, mesh number $m$ with polynomial order $\check{p}$ and continuity $\check{k}$ across element boundaries\footnote{Except for some possible $C^0$ lines in the initial CAD geometry.}. For the corresponding FEM meshes we denote by ${\cal M}_{m,\check{p},\mathrm{s}}^{\textsc{fem}}$ and ${\cal M}_{m,\check{p},\mathrm{i}}^{\textsc{fem}}$ the subparametric and isoparametric FEM meshes, respectively. The construction of NURBS meshes are illustrated in \Cref{Fig2:SphericalShellMeshes}. The initial mesh is depicted as mesh ${\cal M}_{1,\check{p},\check{k}}^{\textsc{iga}}$ in \Cref{Fig2:SphericalShellMeshes1} and is refined only in the angular directions for the first 3 refinements (that is, mesh ${\cal M}_{4,\check{p},\check{k}}^{\textsc{iga}}$ only have one element thickness in the radial direction). Mesh ${\cal M}_{m,\check{p},\check{k}}^{\textsc{iga}}$, $m=5,6,7$, have 2, 4 and 8 elements in its thickness, respectively. This is done to obtain low aspect ratios for the elements. All the meshes will then be nested and the refinements are done uniformly. We shall use the same polynomial order in all parameter directions; $\check{p}_\upxi=\check{p}_\upeta=\check{p}_\upzeta$. Finally, unless otherwise stated, we shall use the BGU formulation and $N=4$ basis functions in the radial direction of the infinite elements.

\subsection{Simpson benchmark}
\begin{figure}
	\centering
	\includegraphics[width=\textwidth]{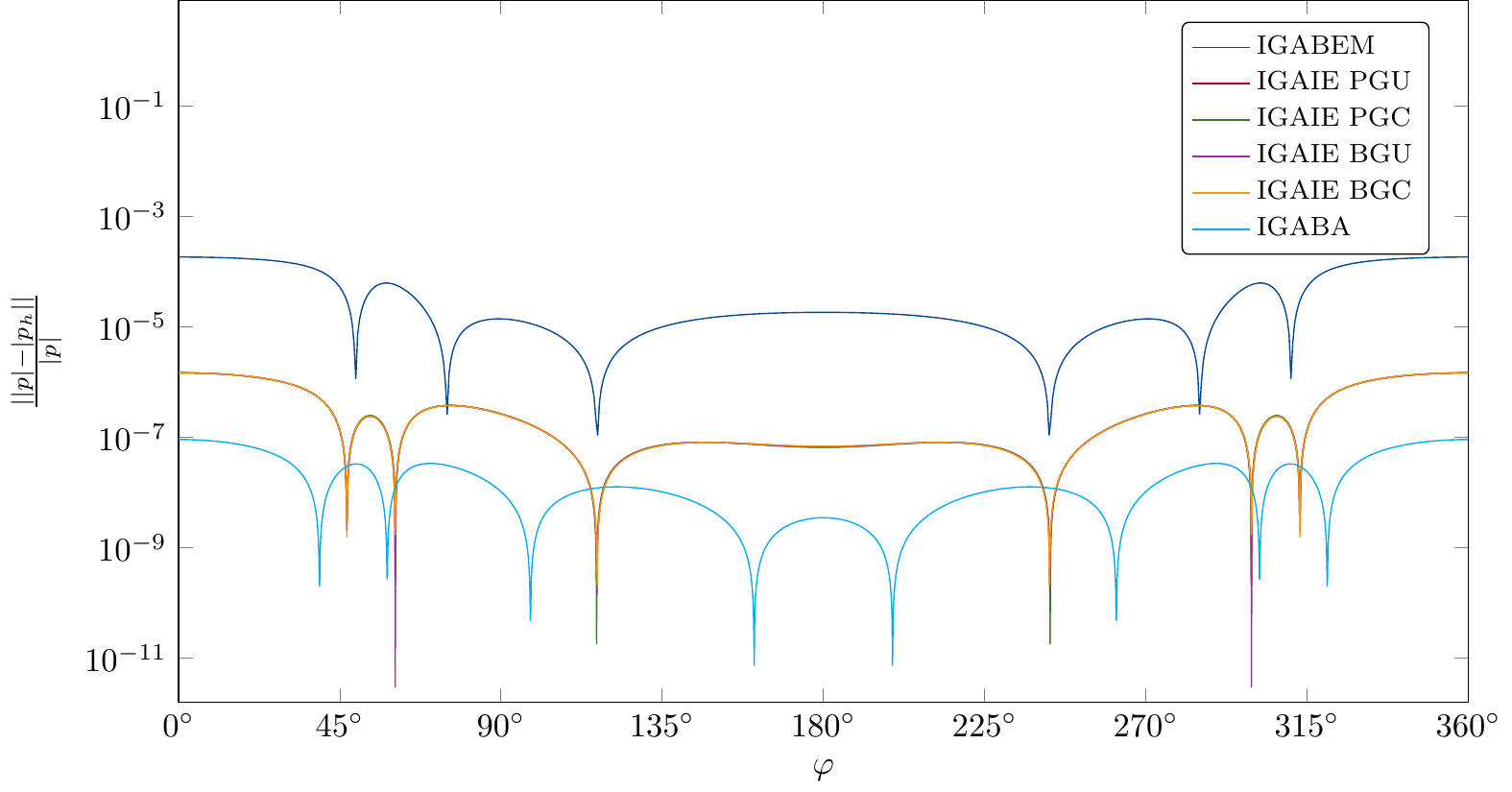}
	\caption{\textbf{Simpson benchmark}: The relative error in the modulus of the pressure is plotted on a circle (azimuth direction, $\varphi$) in the $xy$-plane at $r=\SI{5}{m}$. All simulations were computed on mesh ${\cal M}_{3,3,2}^{\textsc{iga}}$. The IGAIE formulations here produce roughly the same result.}
	\label{Fig2:simpsonPlot}
\end{figure}

The configuration presented by Simpson et al.~\cite{Simpson2014aib} is considered: a rigid sphere of radius $R_0=\SI{0.5}{m}$ is impinged by an incident plane wave and the total pressure is measured at a distance $r=\SI{5}{m}$ from the origin. This is a low frequency problem with $k=\SI{2}{m^{-1}}$. It is emphasized that the trace of the NURBS discretization of the domain $\Omega_{\mathrm{a}}$ at the surface $\Gamma_1$ reduces to the exact same NURBS discretization used in~\cite{Simpson2014aib} to discretize the boundary $\Gamma_1$. From \Cref{Fig2:simpsonPlot} we observe that the IGA infinite element method (IGAIE) exploits the available degrees of freedom at $\Gamma_1$ more effectively than the IGA boundary element method (IGABEM) in~\cite{Simpson2014aib}\footnote{Due to low resolution of the plots in~\cite[Fig. 17]{Simpson2014aib}, the results was reproduced and sampled at 3601 points (rather than 30 points) using our own IGABEM implementation.}. By projecting the analytic solution onto this set of NURBS basis functions at $\Gamma_1$ (the best approximation in the $L_2$-norm by least squares projection, IGA best approximation, IGABA), it is revealed that even more accuracy can potentially be made. This is an inherent problem for Galerkin FEM when solving the Helmholtz equation and is related to the pollution effect \cite{Babuska1995agf}. All IEM formulations (PGU, PGC, BGU and BGC) gave approximately the same result in this case.

\subsection{Ihlenburg benchmark}
\begin{figure}
	\centering
	\includegraphics[width=\textwidth]{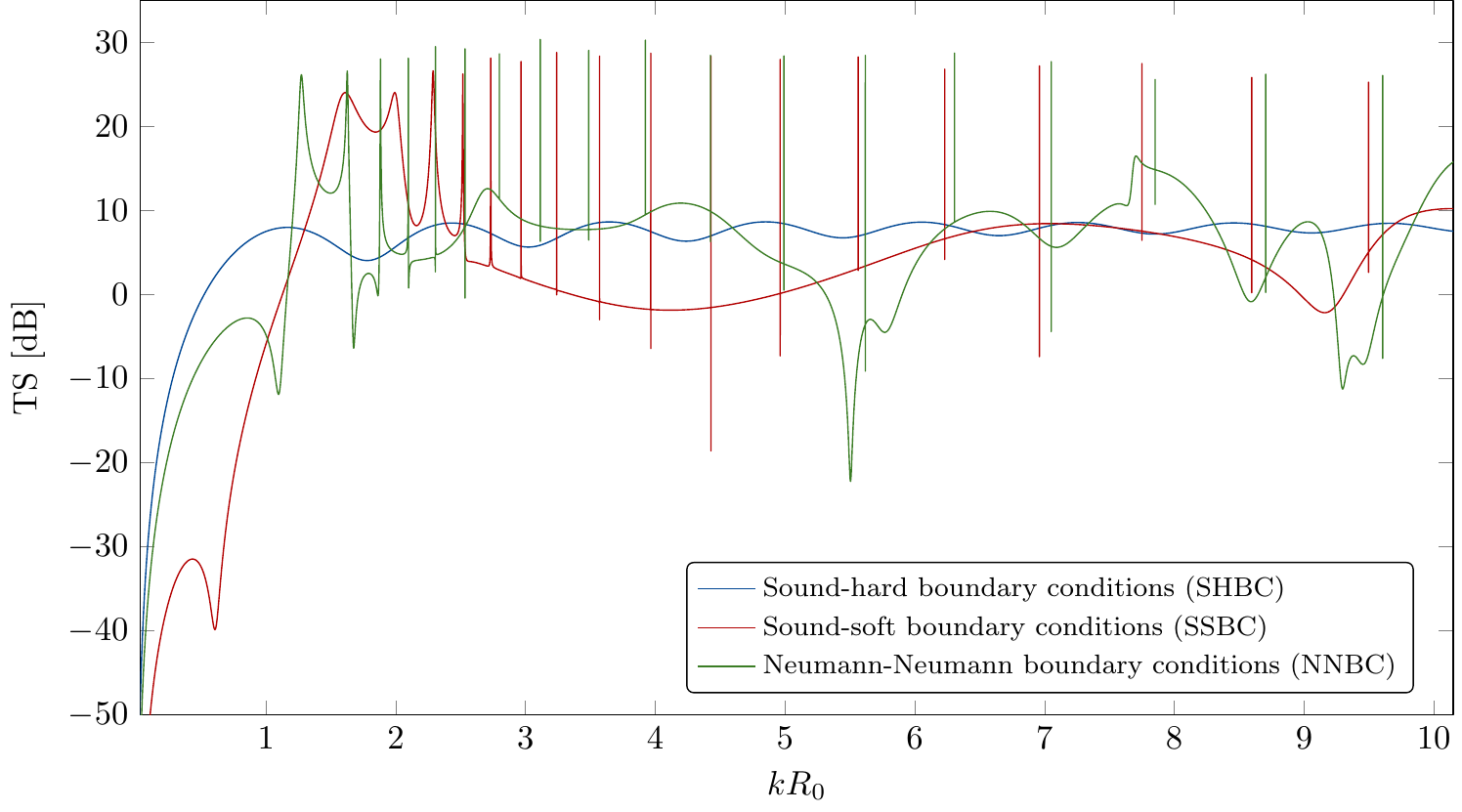}
	\caption{\textbf{Ihlenburg benchmark}: Analytic solutions to the scattering problem on a spherical shell with parameters given in \Cref{Tab2:IhlenburgParameters}. The far field pattern of backscattered pressure is plotted against the wave number $k$. A single Neumann condition at the outer radius, $R_0$, corresponds to the rigid scattering case with $\vec{u}=\zerovec$ and $p_2=0$. ASI at $R_0$ and Neumann at $R_1$ models $p_2=0$. Note that Ihlenburg~\cite[p. 192]{Ihlenburg1998fea} plots the far field pattern in \Cref{Eq2:farfield} instead of the target strength, $\TS$, in \Cref{Eq2:TS}.}
	\label{Fig2:ihlenburgTSexact}
\end{figure}

\begin{figure}
	\centering
	\begin{subfigure}[t]{0.48\textwidth}
		\centering
		\includegraphics[width=\textwidth]{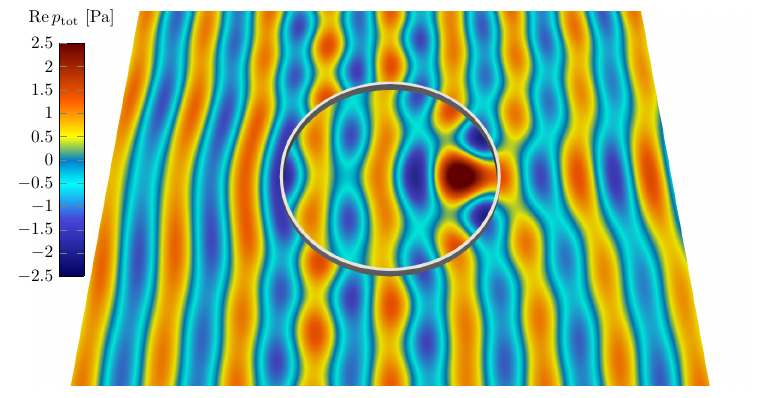}
		\caption{Plot of the real part of the total pressure.}
	\end{subfigure}
	~
	\begin{subfigure}[t]{0.48\textwidth}
		\centering
		\includegraphics[width=\textwidth]{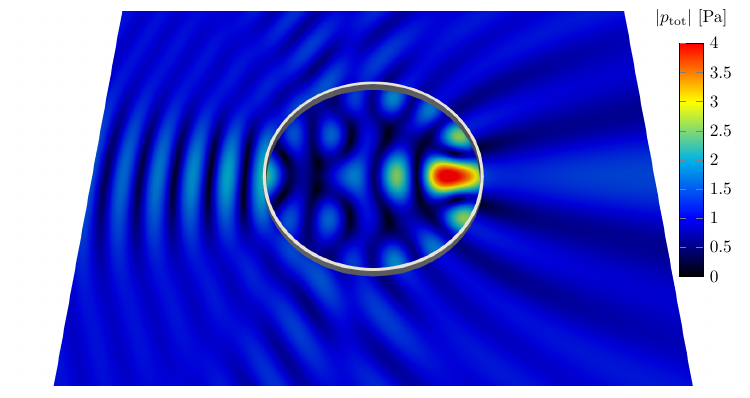}
		\caption{Plot of the modulus of the total pressure.}
	\end{subfigure}
	
	\caption{\textbf{Ihlenburg benchmark with NNBC}: The analytic solution with ASI at both $R_0$ and $R_1$ with $kR_0=10.15$ is plotted in the $xy$-plane. The solid domain is cut open for visualization purposes.}
	\label{Fig2:ihlenburg3Dexact}
\end{figure}

Three benchmark solutions based on the model problem after Ihlenburg~\cite[p. 191]{Ihlenburg1998fea} with parameters given in \Cref{Tab2:IhlenburgParameters}, are investigated. 
\begin{table}
	\centering
	\caption{\textbf{Ihlenburg benchmark}: Parameters for the Ihlenburg benchmark problems.}
	\label{Tab2:IhlenburgParameters}
	\begin{tabular}{l l}
		\toprule
		Parameter & Description\\
		\midrule
		$P_{\mathrm{inc}}=\SI{1}{Pa}$ & Amplitude of incident wave\\
		$E = \SI{2.07e11}{Pa}$ & Young's modulus\\
		$\nu = 0.3$ & Poisson's ratio\\
		$\rho_{\mathrm{s}} = \SI{7669}{kg.m^{-3}}$ & Density of solid\\
		$\rho_{\mathrm{f}} = \SI{1000}{kg.m^{-3}}$ & Density of water\\
		$c_{\mathrm{f}} = \SI{1524}{m.s^{-1}}$ & Speed of sound in water\\
		$R_0=\SI{5.075}{m}$ & Outer radius\\
		$R_1=\SI{4.925}{m}$ & Inner radius\\
		\bottomrule
	\end{tabular}
\end{table}
The parameters for the fluid domains are the speed of sound in water $c_{\mathrm{f}}$ and the fluid density $\rho_{\mathrm{f}}$, and the parameters for the solid domain are the Young's modulus, $E$, the Poisson's ratio $\nu$ and the solid density $\rho_{\mathrm{s}}$. The first benchmark is a simple rigid scattering case (with sound-hard boundary conditions, SHBC) on a sphere with radius $R_0$. The second benchmark problem on a spherical shell has ASI conditions at the outer radius, $R_0$, and homogeneous Neumann condition at the inner radius, $R_1$ (sound-soft boundary conditions, SSBC). This case can be thought of as an approximation of a scattering problem on a spherical shell with an internal fluid with very low density. The third and final benchmark is a further extension with ASI conditions on both sides of the spherical shell (Neumann-Neumann conditions on both surfaces of the shell, NNBC). All of these benchmarks have analytic solutions~\cite{Venas2019e3s} (see \Cref{Fig2:ihlenburgTSexact,Fig2:ihlenburg3Dexact}), which enables computation of the error in the energy norm. As we use the same parameters in both fluids, we denote the common wave number in these fluids by $k=k_1=k_2$.
For each experiment, we use the same NURBS order everywhere. Denote by $\check{p}_\upxi = \check{p}_{\upxi,\mathrm{f}} = \check{p}_{\upxi,\mathrm{s}}$ the common NURBS order in the fluid and the solid in the $\xi$-direction. Similarly $\check{p}_\upeta = \check{p}_{\upeta,\mathrm{f}} = \check{p}_{\upeta,\mathrm{s}}$ and $\check{p}_\upzeta = \check{p}_{\upzeta,\mathrm{f}} = \check{p}_{\upzeta,\mathrm{s}}$. Moreover, we denote by $\check{p} = \check{p}_\upxi = \check{p}_\upeta = \check{p}_\upzeta$ the common polynomial orders in all domains.

In order to compare $C^0$ FEM and IGA on the scattering problem, we shall transform the NURBS mesh to a $C^0$ FEM mesh. We use the technique described in \Cref{Sec2:NURBStransformation} to get an isoparametric B-spline approximation of the geometry (isoparametric FEM). This parametrization will have $C^0$ continuity at element boundaries and correspondingly $G^0$ continuity of the geometry representation (i.e. with kinks). The geometric approximation error is of one order higher than the finite element approximation of the solution \cite{Strang1973aao}, so one could expect the $C^0$-IGA meshes (with $\hat{k}=0$) to produce the same accuracy as the isoparametric FEM meshes of higher order ($\hat{p}\geq 2$). It should be noted that the FEM analysis would then use the Bernstein basis instead of the classical Lagrange basis. However, both of these set of functions spans the same spaces, such that the results should be identical in the absence of round-off errors. 

\begin{figure}
	\centering
	\includegraphics[width=\textwidth]{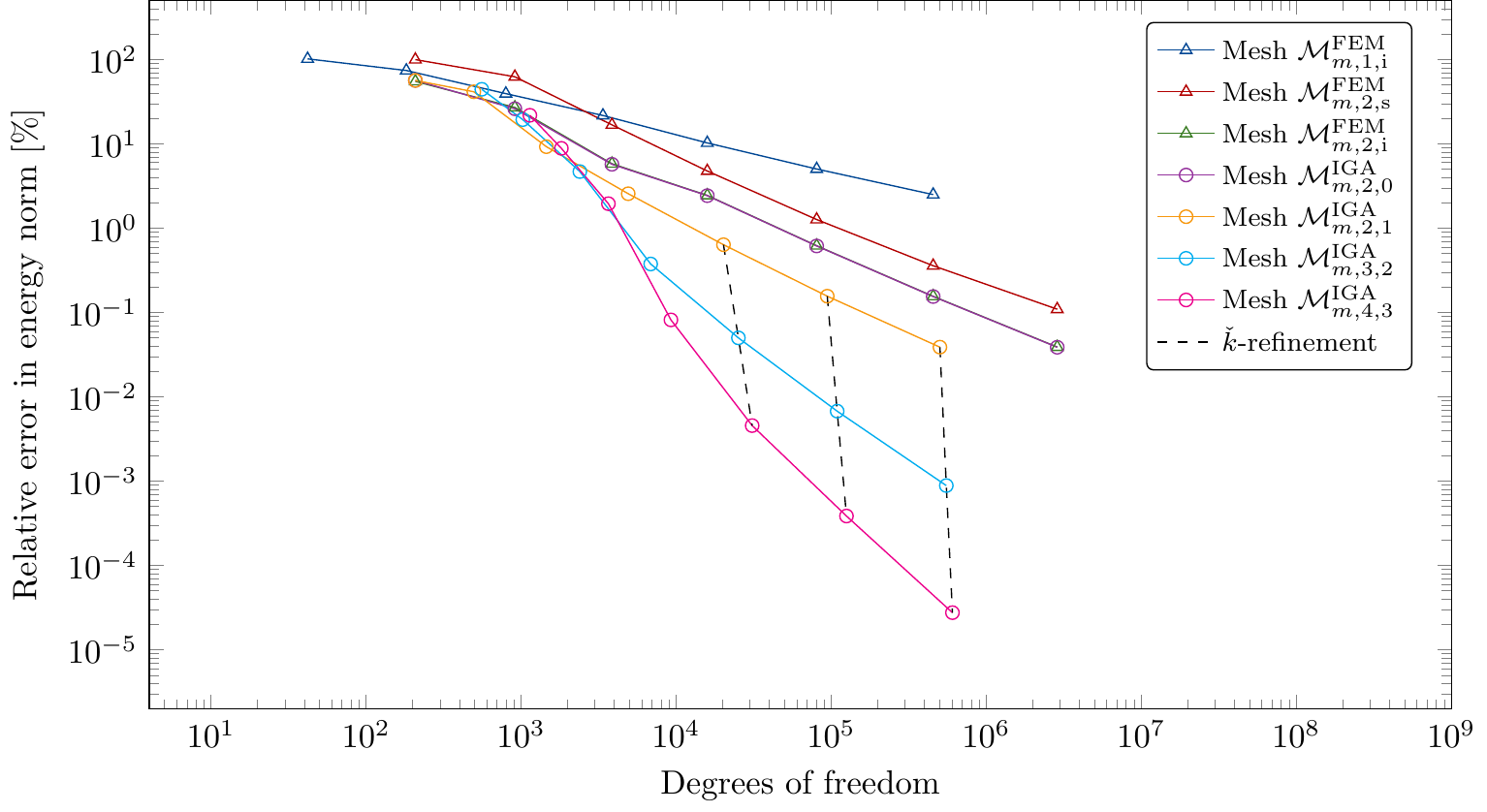}
	\caption{\textbf{Ihlenburg benchmark with SHBC}: Convergence analysis on the rigid scattering case with $k=\SI{1}{m^{-1}}$ and mesh ${\cal M}_m$, $m=1,\dots,7$, using $N=6$. The relative energy error (\Cref{Eq2:energyNormFluids}) is plotted against the degrees of freedom.}
	\label{Fig2:EnergyErrorPlotsDofs}
	\par\bigskip
	\includegraphics[width=\textwidth]{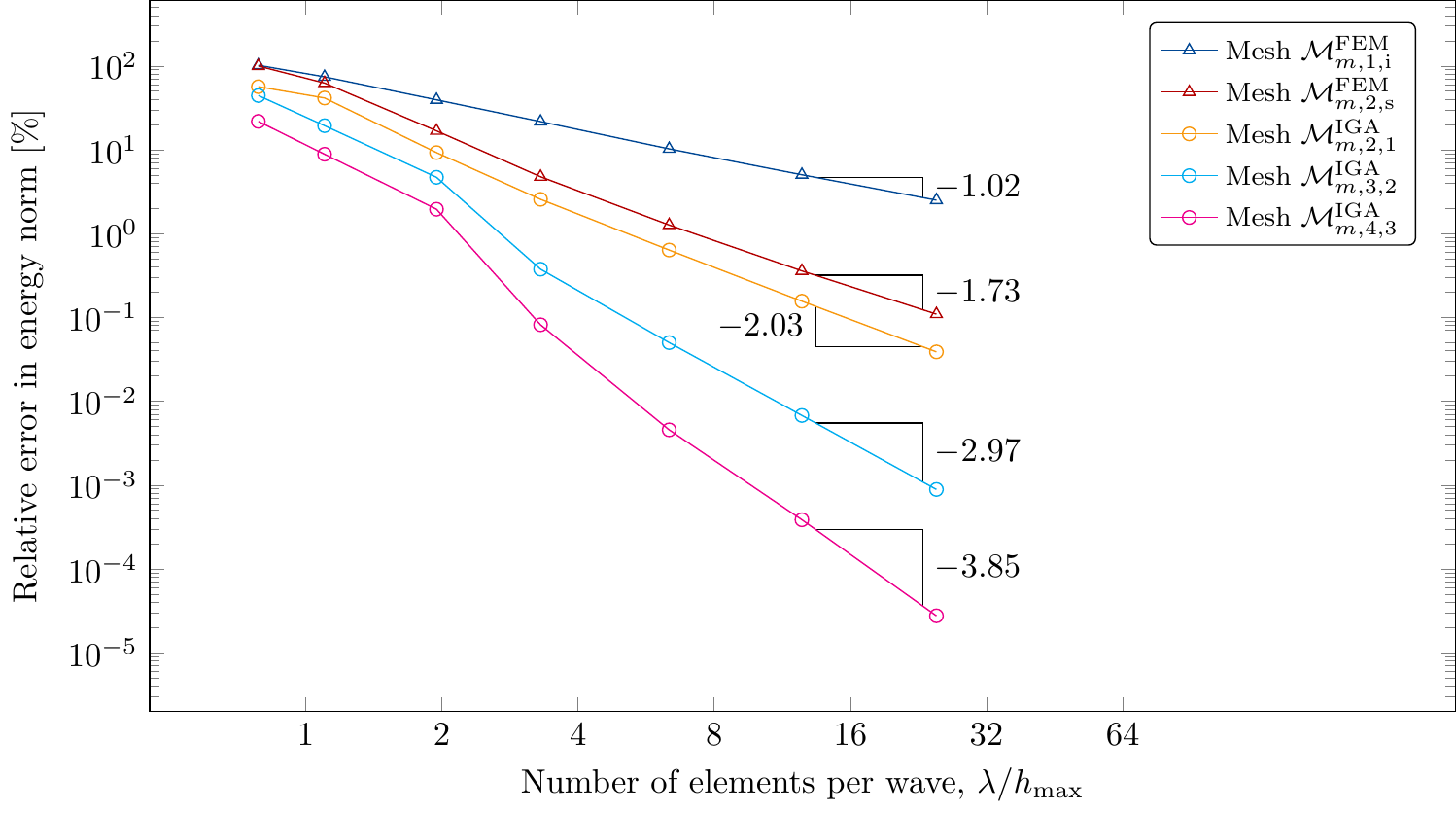}
	\caption{\textbf{Ihlenburg benchmark with SHBC}: Convergence analysis on the rigid scattering case with $k=\SI{1}{m^{-1}}$ and mesh ${\cal M}_m$, $m=1,\dots,7$, using $N=6$. The relative energy error (\Cref{Eq2:energyNormFluids}) is plotted against the number of elements per wave.}
	\label{Fig2:EnergyErrorPlotsh}
\end{figure}

In \Cref{Fig2:EnergyErrorPlotsDofs} we illustrate $h$-refinement through the error in the energy norm for the first benchmark example (rigid scattering).
Predicted convergence rates are not obtained until the aspect ratio of the elements are reduced sufficiently (that is, from mesh ${\cal M}_4$ and onward). By comparing the results of mesh ${\cal M}_{m,2,\mathrm{i}}^{\textsc{fem}}$ and mesh ${\cal M}_{m,2,0}^{\textsc{iga}}$ it can be concluded that the geometry error of mesh ${\cal M}_{m,2,\mathrm{i}}^{\textsc{fem}}$ has almost no impact on the accuracy. However, when using maximum continuity, we get significantly better results. Expected convergence rates are visualized in \Cref{Fig2:EnergyErrorPlotsh} where we now plot the energy norm against $\lambda/h_{\mathrm{max}}$ (corresponding to the number of elements per wave) with $\lambda$ being the wavelength $\lambda=2\pi/k$. A key observation is that the number of elements per wave (needed to obtain a given accuracy) is greatly reduced with higher order IGA methods compared to the classical linear FEM (where 10 elements per wavelength is typically desired for engineering precision, \cite[p. 182]{Ihlenburg1998fea}). The result for the subparametric meshes ${\cal M}_{m,2,\mathrm{s}}^{\textsc{fem}}$ indicates that the convergence rate is reduced due to the reduced accuracy in the geometric representation. This is to be expected as shown in \cite[p. 202]{Strang1973aao}.

\begin{figure}
	\centering
	\begin{subfigure}{0.3\textwidth}
		\centering
		\includegraphics[width=0.6\textwidth]{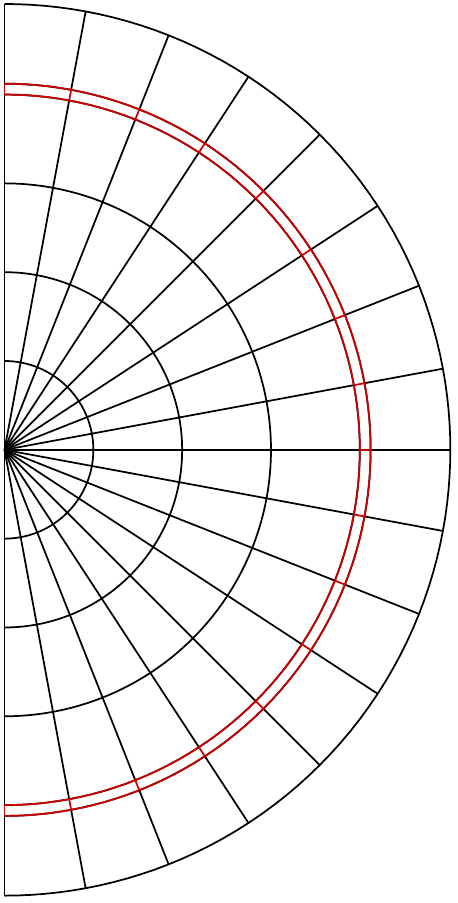}
		\caption{Mesh ${\cal M}_{4,\check{p},\check{k}}^{\textsc{iga}}$}
		\label{Fig2:SphericalShellMeshes1NNBC}
    \end{subfigure}
    ~
	\begin{subfigure}{0.3\textwidth}
		\centering
		\includegraphics[width=0.6\textwidth]{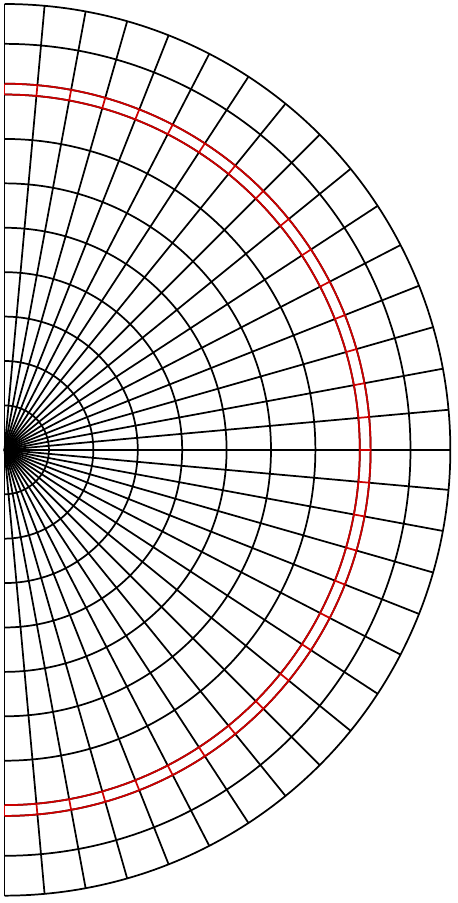}
		\caption{Mesh ${\cal M}_{5,\check{p},\check{k}}^{\textsc{iga}}$}
		\label{Fig2:SphericalShellMeshes2NNBC}
    \end{subfigure}
    ~
	\begin{subfigure}{0.3\textwidth}
		\centering
		\includegraphics[width=0.6\textwidth]{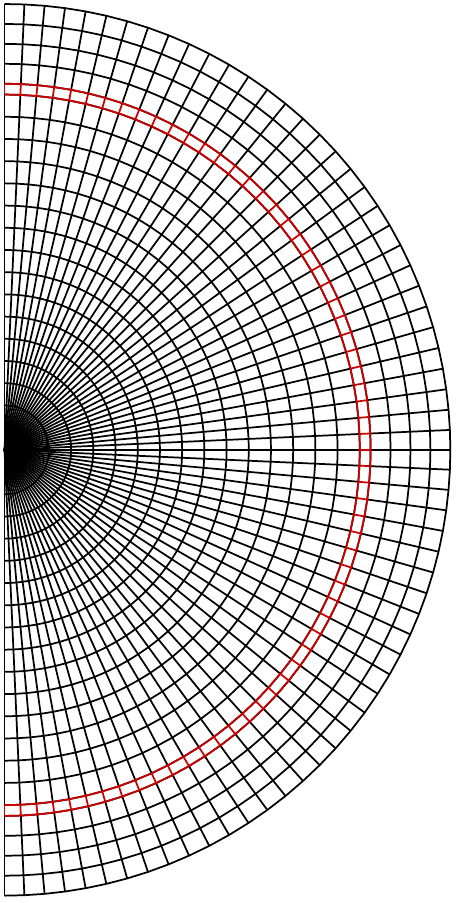}
		\caption{Mesh ${\cal M}_{6,1,\mathrm{i}}^{\textsc{fem}}$}
		\label{Fig2:SphericalShellMeshes3NNBC}
    \end{subfigure}
	\caption{\textbf{Ihlenburg benchmark with NNBC}: Illustration of some meshes for the full ASI problem in the $xz$-plane ($x>0$), where the mesh lines for the solid domain is colored red. The full mesh is obtained by rotation around the $z$-axis. Mesh ${\cal M}_{5,2,\mathrm{i}}^{\textsc{fem}}$ is visually indistinguishable from ${\cal M}_{5,\check{p},\check{k}}^{\textsc{iga}}$.}
	\label{Fig2:SphericalShellMeshesNNBC}
\end{figure}
\begin{figure}
	\centering
	\includegraphics[width=\textwidth]{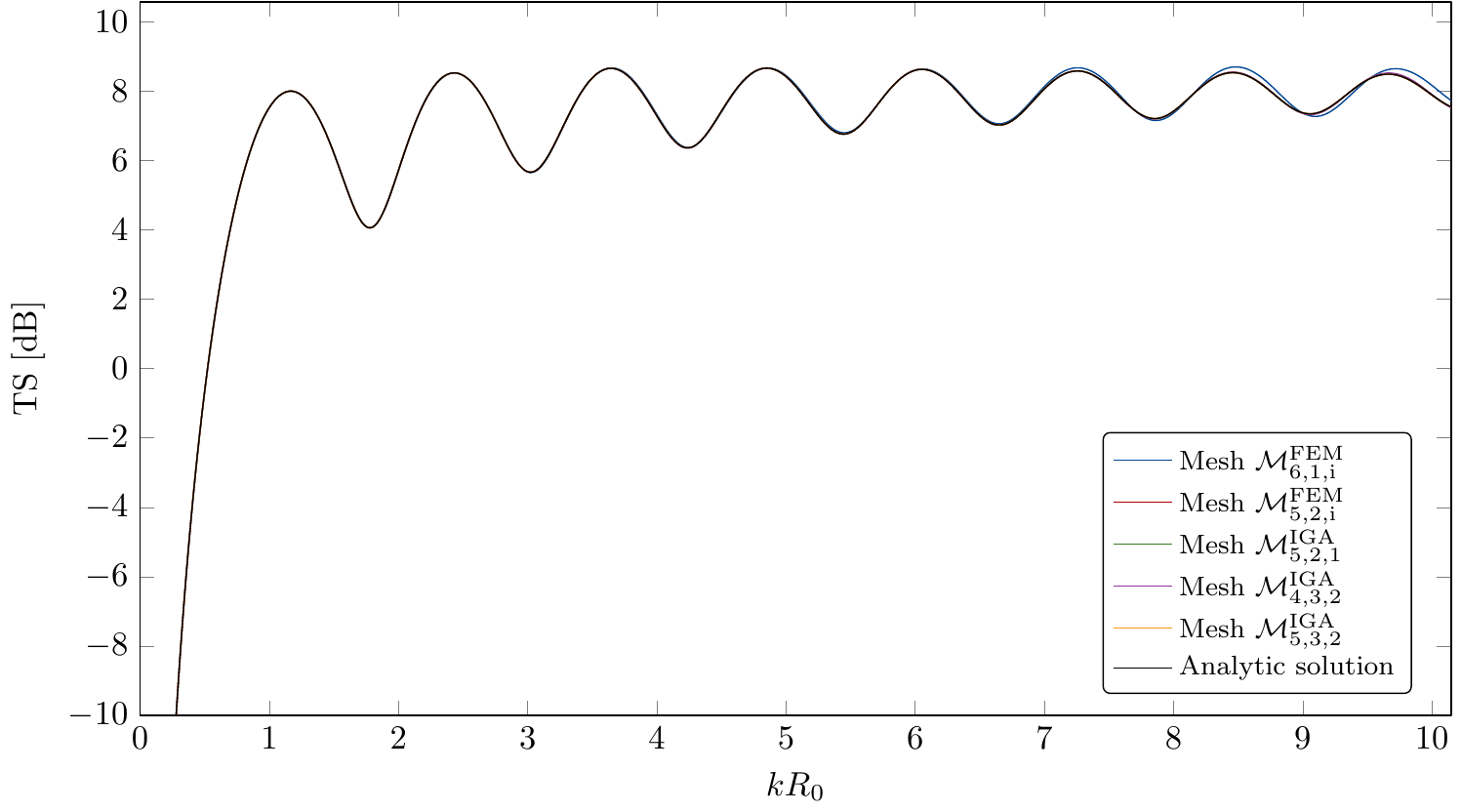}
	\caption{\textbf{Ihlenburg benchmark with SHBC}: The target strength (TS) in \Cref{Eq2:TS} is plotted against $kR_0$.}
	\label{Fig2:TSPlotSHBC}
	\par\bigskip
	\par\bigskip
	\includegraphics[width=\textwidth]{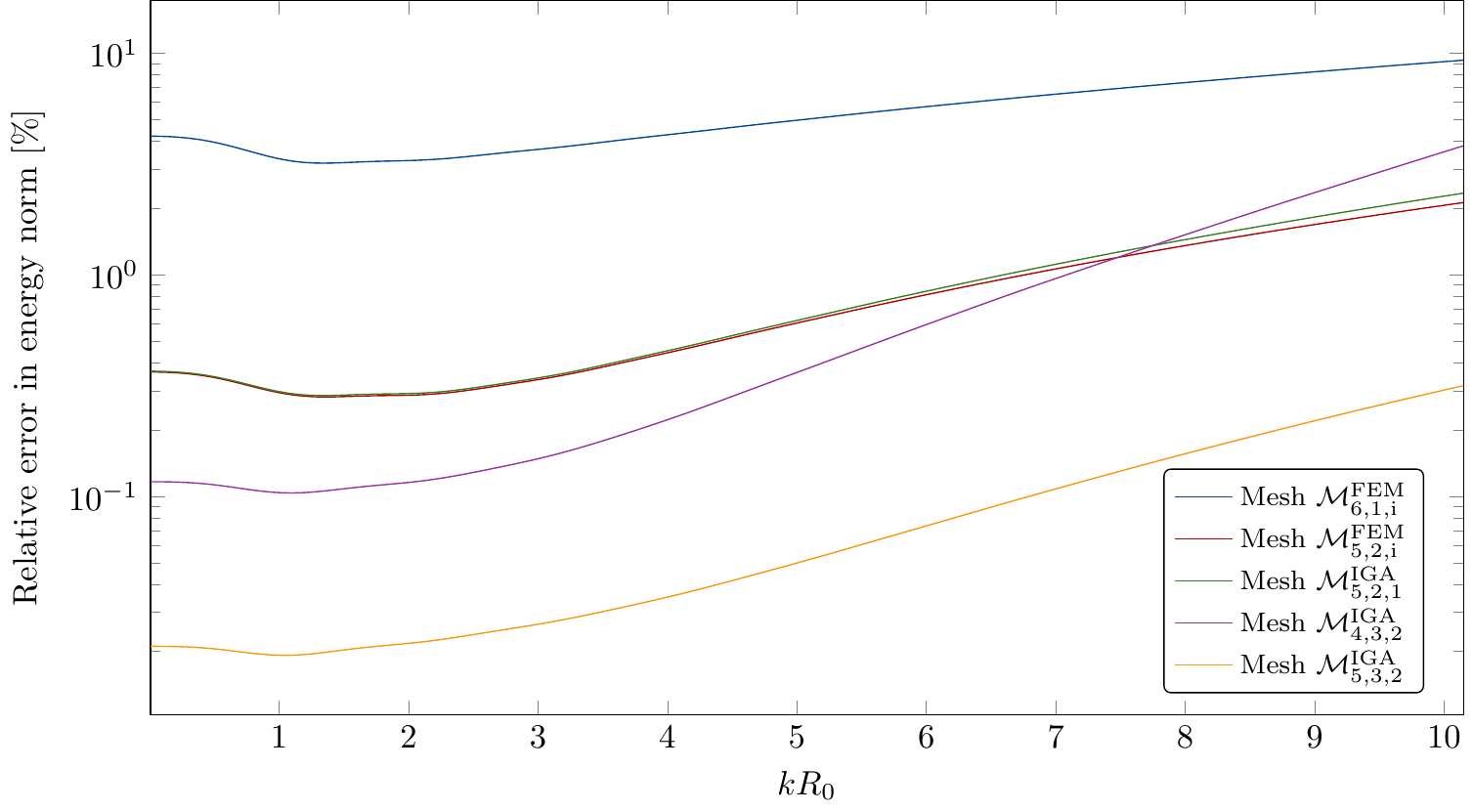}
	\caption{\textbf{Ihlenburg benchmark with SHBC}: The relative energy norm (\Cref{Eq2:energyNorm}) is plotted against $kR_0$.}
	\label{Fig2:errorPlotSHBC}
\end{figure}
\begin{figure}
	\centering
	\includegraphics[width=\textwidth]{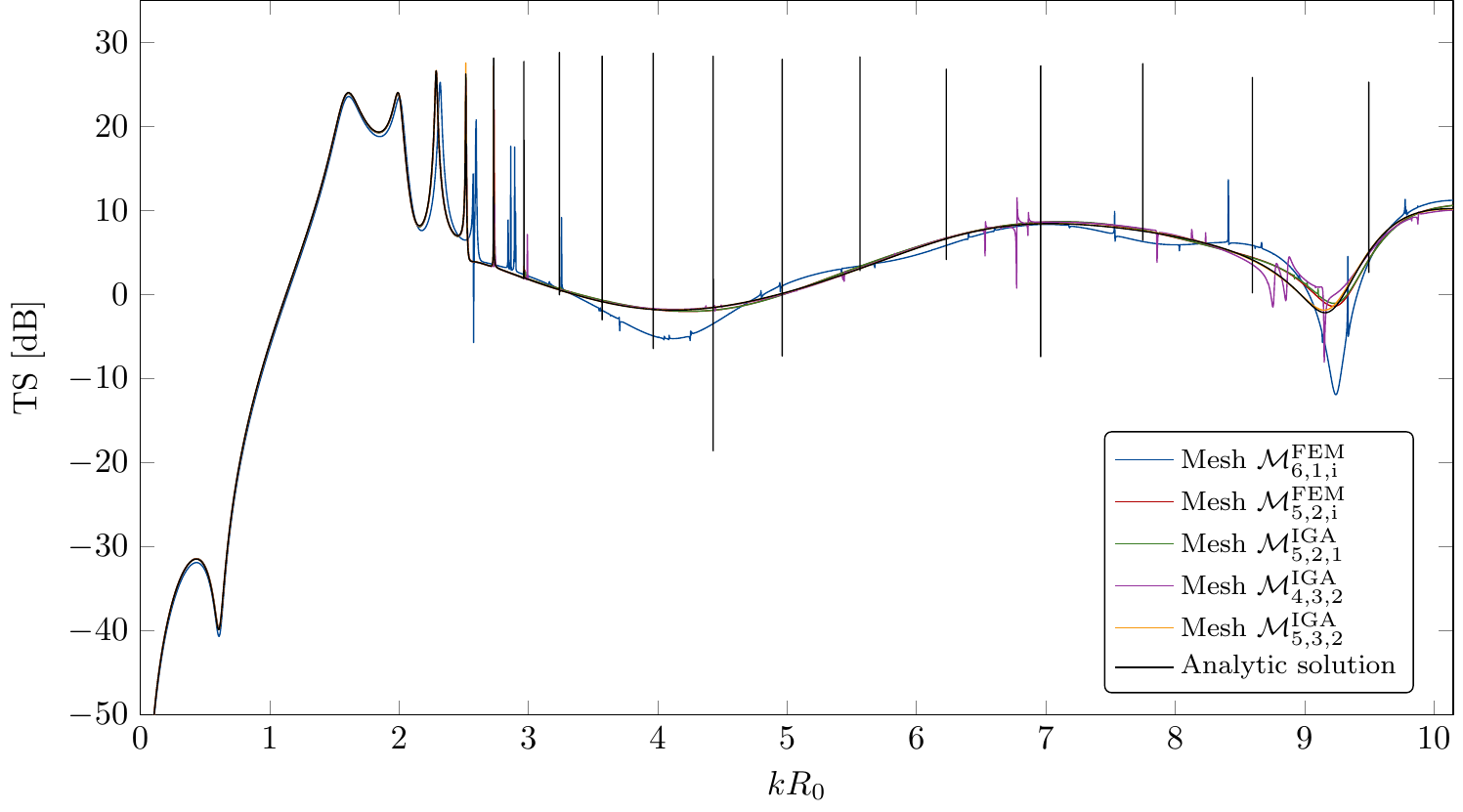}
	\caption{\textbf{Ihlenburg benchmark with SSBC}: ASI problem with the internal pressure modeled to be $p_2=0$. The target strength (TS) in \Cref{Eq2:TS} is plotted against $kR_0$.}
	\label{Fig2:TSPlotSSBC}
	\par\bigskip
	\par\bigskip
	\includegraphics[width=\textwidth]{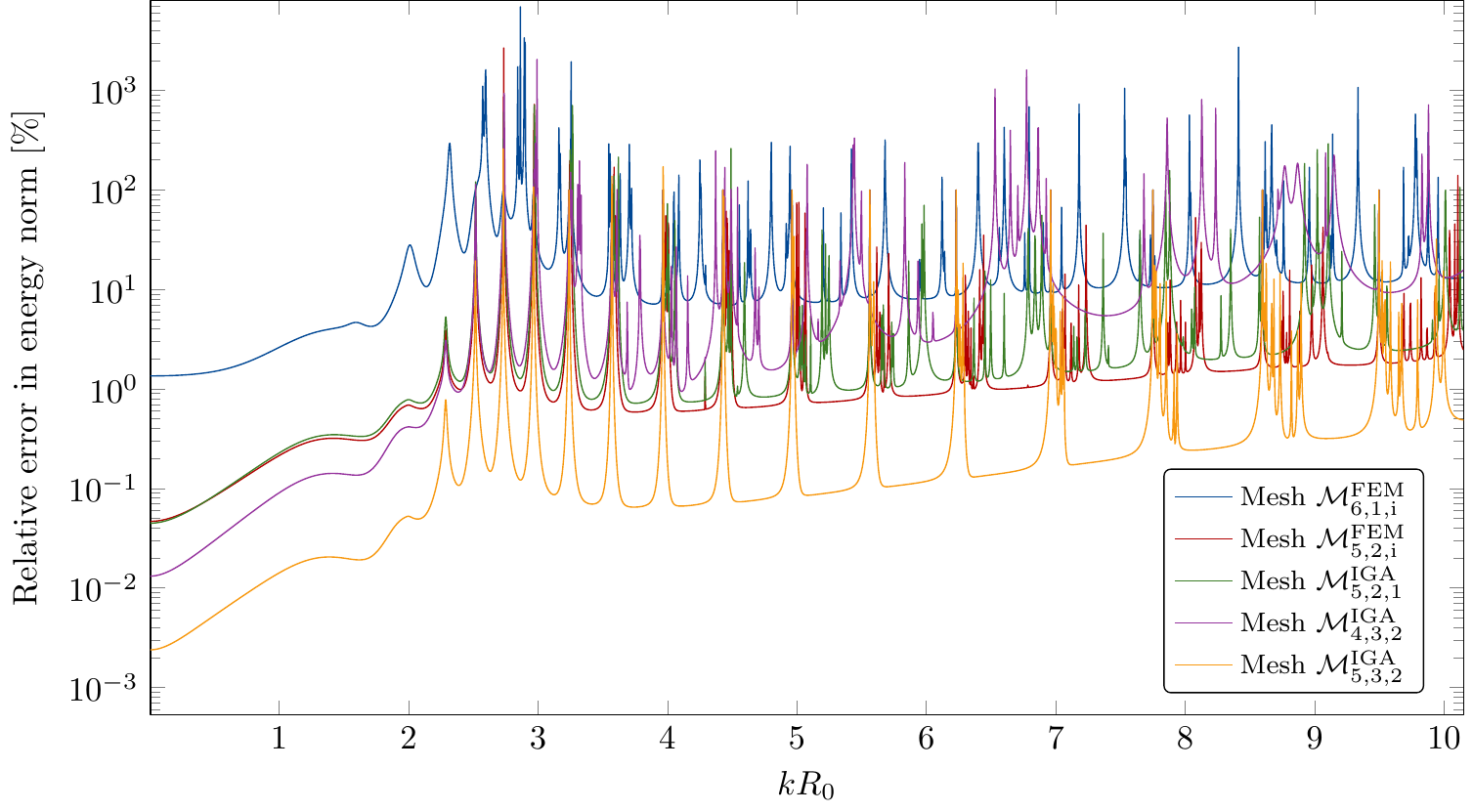}
	\caption{\textbf{Ihlenburg benchmark with SSBC}: ASI problem with the internal pressure modeled to be $p_2=0$. The relative energy norm (\Cref{Eq2:energyNorm}) is plotted against $kR_0$.}
	\label{Fig2:errorPlotSSBC}
\end{figure}
\begin{figure}
	\centering
	\includegraphics[width=\textwidth]{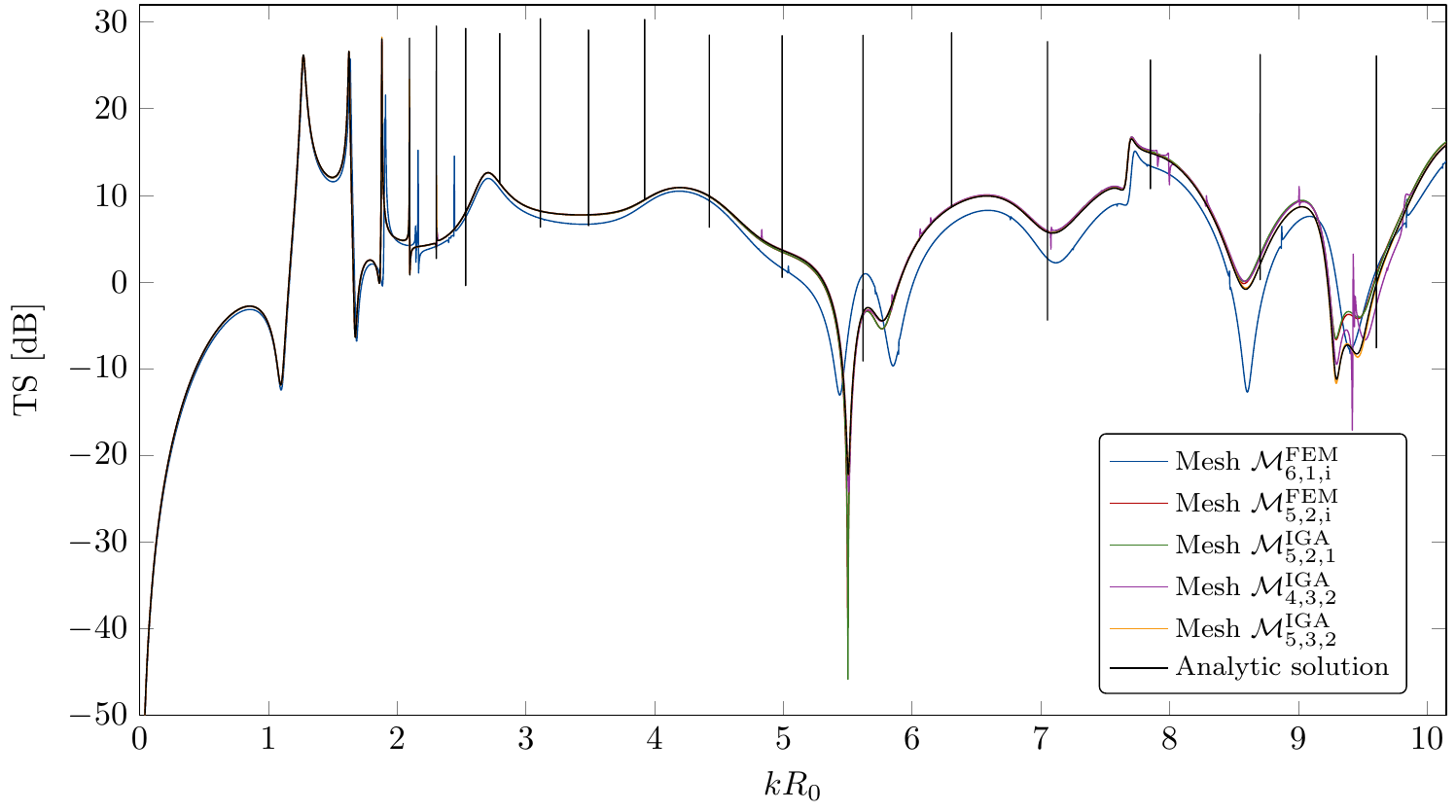}
	\caption{\textbf{Ihlenburg benchmark with NNBC}: The target strength (TS) in \Cref{Eq2:TS} is plotted against $kR_0$.}
	\label{Fig2:TSPlotNNBC}
	\par\bigskip
	\par\bigskip
	\includegraphics[width=\textwidth]{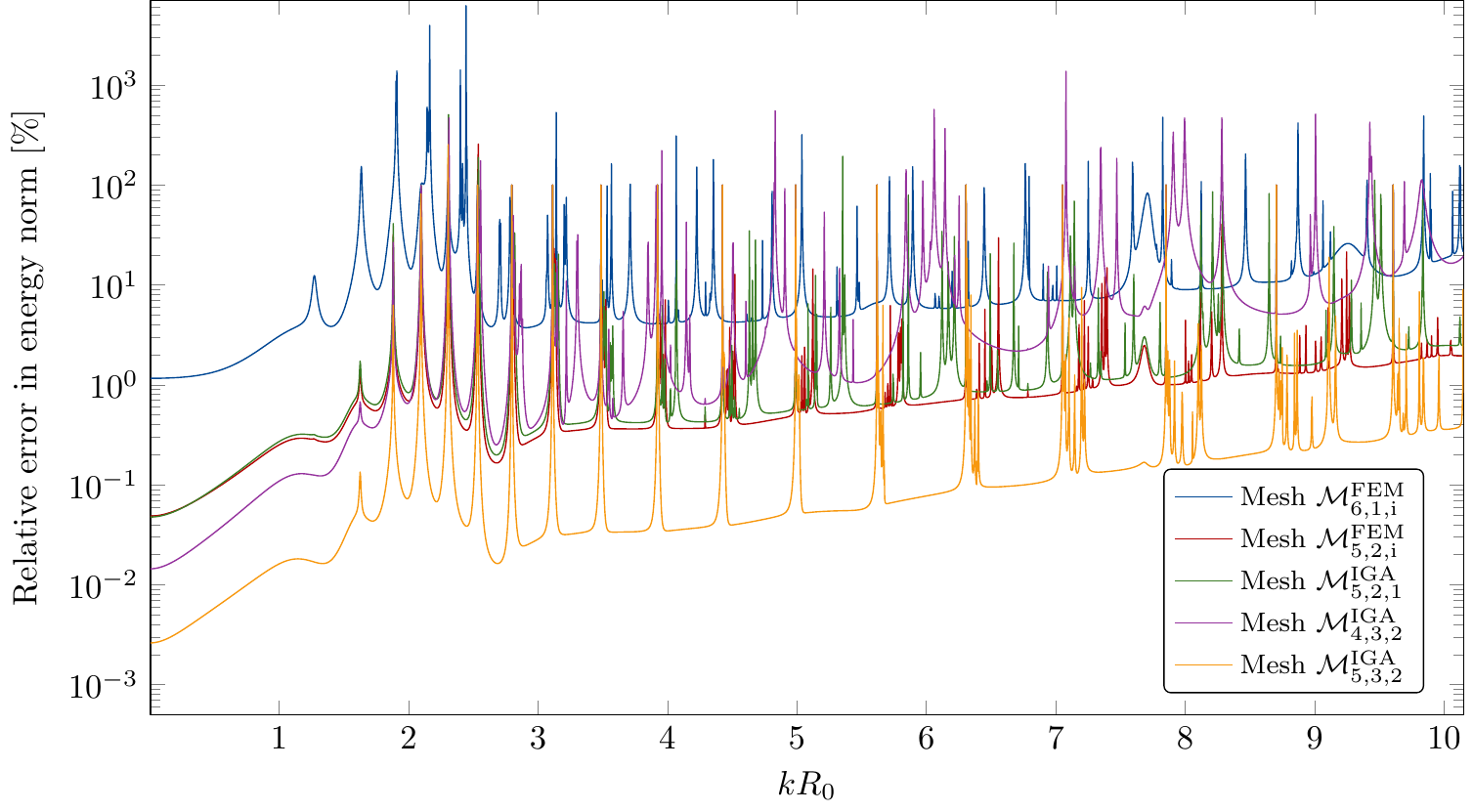}
	\caption{\textbf{Ihlenburg benchmark with NNBC}: The relative energy norm (\Cref{Eq2:energyNorm}) is plotted against $kR_0$.}
	\label{Fig2:errorPlotNNBC}
\end{figure}

Approaching the ASI problems, we illustrate some meshes in \Cref{Fig2:SphericalShellMeshesNNBC} for the full ASI problem. The corresponding meshes for the SSBC problem (with $p_2=0$) are obtained by removing the mesh inside the solid domain. In \Cref{Fig2:TSPlotSHBC,Fig2:errorPlotSHBC,Fig2:TSPlotSSBC,Fig2:errorPlotSSBC,Fig2:TSPlotNNBC,Fig2:errorPlotNNBC} the target strength, $\TS$, and the error in the energy norm is plotted against the scaled wave number, $kR_0$, in all of the three Ihlenburg benchmarks. As each frequency sweep is computed with a different number of degrees of freedom, one should draw the conclusions based on comparing both the accuracy of the results and the related computational costs.

Some data from simulations at $k=\SI{1}{m^{-1}}$ are reported in \Cref{Tab2:dataRigidScattering} (simulation run with 12 processors of the type Intel(R) Xeon(R) CPU E5-4650 2.70GHz). It should be noted that all simulations were done using the same code, such that the computational time for the FEM simulations can be optimized. However, this is actually the case for the IGA code as well since the implementation does not utilize optimized quadrature rules. The integration is done with $(\check{p}+1)^3$ quadrature points per element when building the system. For higher order splines spaces this is significantly more quadrature point than what is needed for exact integration (on meshes with affine geometry mapping\footnote{Using the same quadrature scheme on truly isoparametric elements will according to \cite[p. 256]{Ciarlet1991bee} give a numerical integration error of the same order as the finite element discretization error. Thus, the argument for optimal quadrature scheme also holds for isoparametric elements as well.}). In \cite{Hughes2010eqf,Johannessen2017oqf}, it is shown that the optimal number of quadrature points is half the number of degrees of freedom of the splines space under consideration. That is, the number of quadrature points in the IGA 3D tensor product meshes can be reduced by a factor up to $2^3(\check{p}+1)^3$ for meshes with maximal continuity. Thus, the efficiency of the IGA simulation may be improved significantly.

A particular interesting observation is that IGA obtains roughly the same accuracy as FEM when the same number of elements is used, even though this corresponds to far less degrees of freedom for the IGA simulation. Moreover, even better result can be obtained with less degrees of freedom if the polynomial degree is increased in the IGA simulations. This, however, only occurs when the mesh resolves the number of waves per element. When the mesh is sufficiently resolved, one order of magnitude improvement in the accuracy is obtained by increasing the polynomial degree. Since another magnitude of accuracy is obtained by the use of higher order elements in FEM/IGA, the IGA offers several orders of magnitude better accuracy than classical linear FEM. 

The peaks in the frequency sweeps represent eigenmodes. The quality of the numerical approximation of the corresponding frequencies is reduced for higher frequencies, resulting in fictitious modes. This typically does not pose that much of a problem as the bandwidth of these eigenmodes becomes very small, with a corresponding reduction in the energy they represent. Note that mesh ${\cal M}_{4,3,2}^{\textsc{iga}}$ performs particularly poorly on the partial ASI problem due to a fictitious mode at $k=\SI{1}{m^{-1}}$ for this mesh. The improvement offered by IGA concerning the accuracy in the eigenmodes is investigated in~\cite{Venas2015iao}. 

It should be noted that the meshes used throughout this work are not optimal. This is in particular the case for the full ASI problem where the density of elements becomes large at the origin. These meshes were used as they naturally arise from tensor product NURBS meshes of spherical shells and spheres. One could thus obtain increased performance for the FEM solutions using standard meshing of the domain. However, locally refined meshes can also be obtained with the IGA method, for example using LR B-splines \cite{Johannessen2014iau}.

\begin{table}
	\centering
	\caption{\textbf{Ihlenburg benchmark}: Data for some simulations on the rigid scattering problem with $k=\SI{1}{m^{-1}}$. The errors are given in the energy norm (\Cref{Eq2:energyNorm}). For each simulation, the mesh number, the polynomial order, $\check{p}$, the number of mesh elements $n_{\mathrm{el}}$ (not including the infinite elements) and the number of degrees of freedom $n_{\mathrm{dof}}$, is reported. The elapsed times for building the system $t_{\mathrm{sys}}$ and for solving the system $t_{\mathrm{sol}}$ (using LU-factorization) are also included (times in seconds). Finally, the relative error in the energy norm is given in percentage.}
	\label{Tab2:dataRigidScattering}
	\begin{subtable}[t]{\linewidth}
		\caption{Sound-hard boundary conditions (SHBC).}
		\label{Tab2:dataRigidScatteringSHBC}
		\centering
		\bgroup
		\def\arraystretch{1.1}%
		\begin{tabular}{l S[table-format = 6.0] S[table-format = 5.0] S[table-format = 2.1,round-mode=places,round-precision=1] S[table-format = 2.1,round-mode=places,round-precision=1] S[table-format = 1.2]}
			\hline
			 		   & {$n_{\mathrm{el}}$} & {$n_{\mathrm{dof}}$} &{$t_{\mathrm{sys}}$ [s]}	& {$t_{\mathrm{sol}}$ [s]} 		& {Relative energy error [\%]} \\
			\hline
			\input{\contents/IGAFEM_SHBC.tex}
			\hline
		\end{tabular}
		\egroup
	\end{subtable}
	\par\bigskip
	\begin{subtable}[t]{\linewidth}
		\caption{Sound-soft boundary conditions (SSBC).}
		\label{Tab2:dataRigidScatteringSSBC}
		\centering
		\bgroup
		\def\arraystretch{1.1}%
		\begin{tabular}{l S[table-format = 6.0] S[table-format = 6.0] S[table-format = 2.1,round-mode=places,round-precision=1] S[table-format = 3.1,round-mode=places,round-precision=1] S[table-format = 2.2]}
			\hline
			 		   & {$n_{\mathrm{el}}$} & {$n_{\mathrm{dof}}$} &{$t_{\mathrm{sys}}$ [s]}	& {$t_{\mathrm{sol}}$ [s]} 		& {Relative energy error [\%]} \\
			\hline
			\input{\contents/IGAFEM_SSBC.tex}
		    \hline
		\end{tabular}
		\egroup
	\end{subtable}	
	\par\bigskip
	\begin{subtable}[t]{\linewidth}
		\caption{Neumann-Neumann boundary conditions (NNBC).}
		\label{Tab2:dataRigidScatteringNNBC}
		\centering
		\bgroup
		\def\arraystretch{1.1}%
		\begin{tabular}{l S[table-format = 6.0] S[table-format = 6.0] S[table-format = 2.1,round-mode=places,round-precision=1] S[table-format = 3.1,round-mode=places,round-precision=1] S[table-format = 1.2]}
			\hline
			 		  & {$n_{\mathrm{el}}$} & {$n_{\mathrm{dof}}$} &{$t_{\mathrm{sys}}$ [s]}	& {$t_{\mathrm{sol}}$ [s]} 		& {Relative energy error [\%]} \\
			\hline
			\input{\contents/IGAFEM_NNBC.tex}
		    \hline
		\end{tabular}
		\egroup
	\end{subtable}	
\end{table}

\begin{figure}
	\centering
	\includegraphics[width=0.7\textwidth]{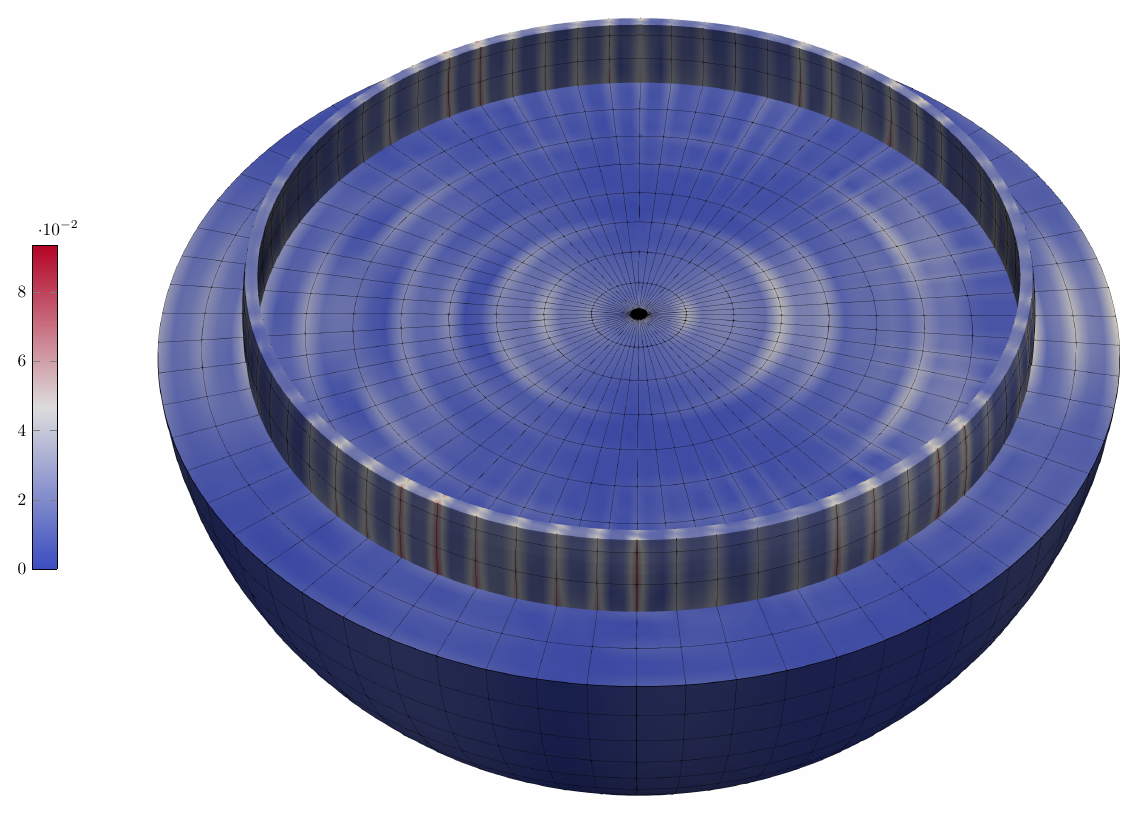}
	\caption{\textbf{Ihlenburg benchmark with NNBC}: Simulation of the full ASI problem on mesh ${\cal M}_{5,2,1}^{\textsc{iga}}$. Pointwise evaluation of the square root of the integrand of the volume integrals in the energy norm $\energyNorm{U-U_h}{\Omega}$ in \Cref{Eq2:energyNorm} with $k=\SI{2}{m^{-1}}$ (error in the infinite elements in $\Omega_{\mathrm{a}}^+$ is not shown) is here visualized, where $U$ is the set of analytic solutions in both fluid domains and the solid domain, and $U_h$ is the corresponding numerical solution. The values are scaled by the square root of the maximum of the corresponding integrand values of $\energyNorm{U}{\Omega}$. Both fluid domains are cut open at the $xy$-plane (at $z=0$), and the solid domain is cut open at $z=\SI{1.1}{m}$.}
	\label{Fig2:energyErrorPlot}
\end{figure}

In \Cref{Fig2:energyErrorPlot} we visualize the distribution of the error of the full ASI problem. The error is observed to be largest at element boundaries where the continuity is reduced. Since second order basis functions are used and the error in the velocity/stress dominates the error in the pressure/displacement, the results are in agreement with what was observed in \cite{Kumar2017spr}, i.e., that the error in the derivative of the primary solution is largest at the element boundaries.

\subsection{Radial pulsation from a mock shell}
\label{Subsec2:mockShell}
By construction of the fundamental solution of the Helmholtz equation ($\Phi_k(\vec{x},\vec{y})$ in \Cref{Eq2:FreeSpaceGrensFunction}), the function $p(\vec{x}) = \Phi_k(\vec{x},\vec{y})$ is a solution to \Cref{Eq2:HelmholtzEqn,Eq2:HelmholtzEqnNeumannCond,Eq2:sommerfeldCond} whenever $\vec{y}\in\R^3\setminus\overline{\Omega^+}$ and for the Neumann boundary condition $g(\vec{x})=\partial_n\Phi_k(\vec{x},\vec{y})$ on $\Gamma_1$.  Hence, we have an exact reference solution for the exterior Helmholtz problem for arbitrary geometries $\Gamma_1$ which encloses the point $\vec{y}$. It is emphasized that this solution is non-physical for non-spherical geometries $\Gamma_1$. General solutions may be constructed by separation of variables (cf.~\cite[p. 26]{Ihlenburg1998fea})
\begin{equation}
	p(\vec{x}) = \sum_{n=0}^\infty\sum_{m=-n}^n C_{nm} h_n^{(1)}(kR) P_n^{|m|}(\cos\vartheta)\euler^{\imag m\varphi} 
\end{equation}
with
\begin{equation*}
	R = |\vec{x}-\vec{y}|,\quad \vartheta=\arccos\left(\frac{x_3-y_3}{R}\right),\quad\varphi = \operatorname{atan2}(x_2-y_2,x_1-y_1)
\end{equation*}
where $h_n^{(1)}$ is the $n^{\mathrm{th}}$ spherical Hankel function of first kind and $P_n^m$ are the associated Legendre functions. In fact, the solution $p(\vec{x}) = \Phi_k(\vec{x},\vec{y})$ is a special case of this general form with 
\begin{equation}
	C_{nm} = \begin{cases}
		\frac{\imag k}{4\PI} & n = 0,\,\,m=0\\
		0 & \text{otherwise}.
		\end{cases}
\end{equation}
The complexity of this problem setup does not scale with the complexity of the model as it is independent of $\Gamma_1$. However, it preserves two important properties of acoustic scattering, namely the radial decay and the oscillatory nature. Thus, this problem setup represents a general way of constructing manufactured solutions, that can be utilized to verify the correctness of the implemented code for solving the Helmholtz equation. A special case of this general setup is the pulsating sphere example in~\cite{Simpson2014aib}.

From the first limit of \Cref{Eq2:Phi_k_limits}, the far field is given by $p_0(\hat{\vec{x}})=\frac{1}{4\PI}\euler^{-\imag k \hat{\vec{x}}\cdot\vec{y}}$. Thus, the target strength is a constant, $\TS=-20\log_{10}(4\PI)\approx -21.984$ (where we define $P_{\mathrm{inc}}=\SI{1}{Pa}$ in \Cref{Eq2:TS} for this problem). 

Consider the case $\vec{y}=\frac{R_0}{4}(1,1,1)$ and the boundary $\Gamma_1$ given by a \textit{mock shell} composed of a cylinder with hemispherical endcaps (with axis of symmetry along the $x$-axis such that the center of the spherical endcaps are located at $x = 0$ and $x = -L$). The cylinder has radius $R_0=\SI{1}{m}$ and length $L=\frac{\PI}{2}R_0$. The analytic solution is given by
\begin{equation}
	p(\vec{x}) = \frac{\euler^{\imag kR}}{4\PI R},\quad R=|\vec{x}-\vec{y}|
\end{equation}
and the Neumann condition is then
\begin{equation}
	g(\vec{x}) = \frac{\euler^{\imag kR}}{4\PI R^3}(\imag kR -1) (\vec{x}-\vec{y})\cdot\vec{n}(\vec{x}).
\end{equation}

\begin{figure}
	\centering    
	\begin{subfigure}[b]{0.49\textwidth}
		\centering
		\includegraphics[width=\textwidth]{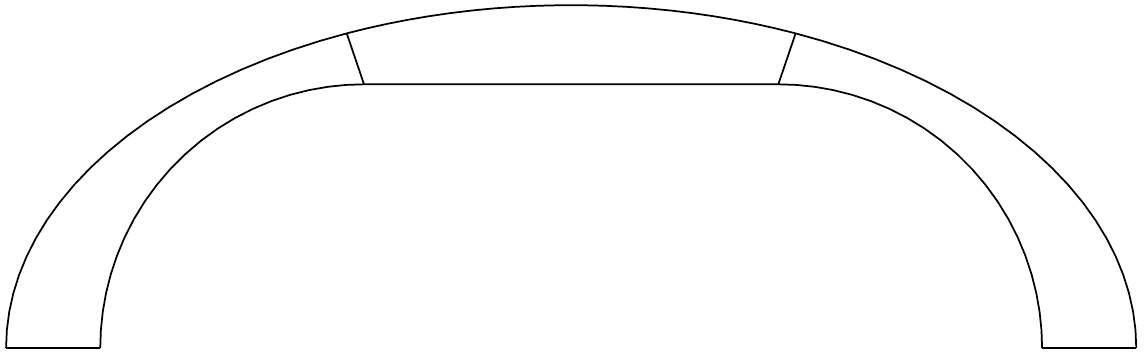}
		\caption{Mesh 1.}
	\end{subfigure}
	~    
	\begin{subfigure}[b]{0.49\textwidth}
		\includegraphics[width=\textwidth]{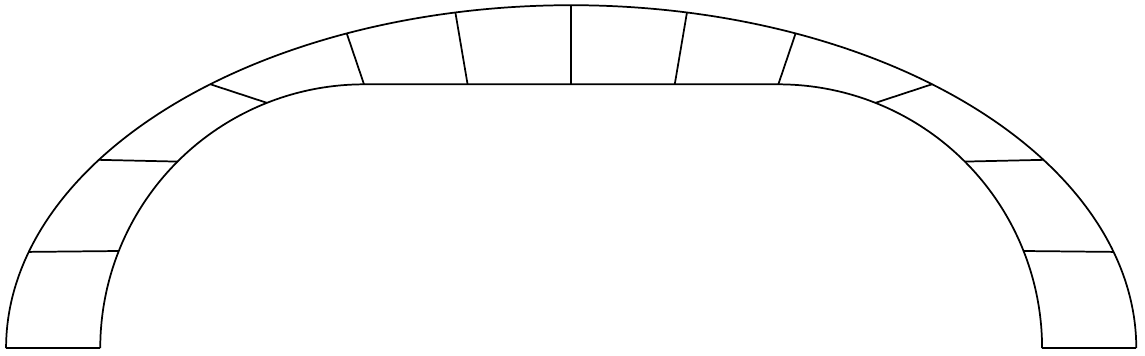}
		\caption{Mesh 3.}
	\end{subfigure}
	\par\bigskip
	\begin{subfigure}[b]{0.49\textwidth}
		\centering
		\includegraphics[width=\textwidth]{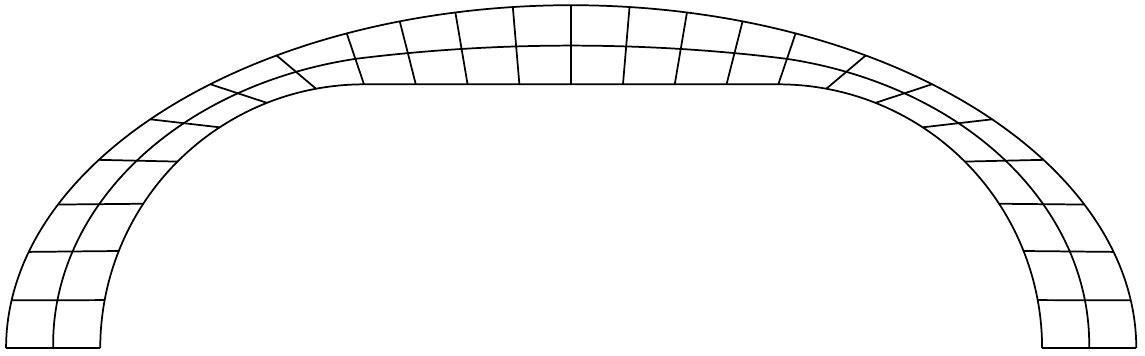}
		\caption{Mesh 4.}
	\end{subfigure}
	~    
	\begin{subfigure}[b]{0.49\textwidth}
		\includegraphics[width=\textwidth]{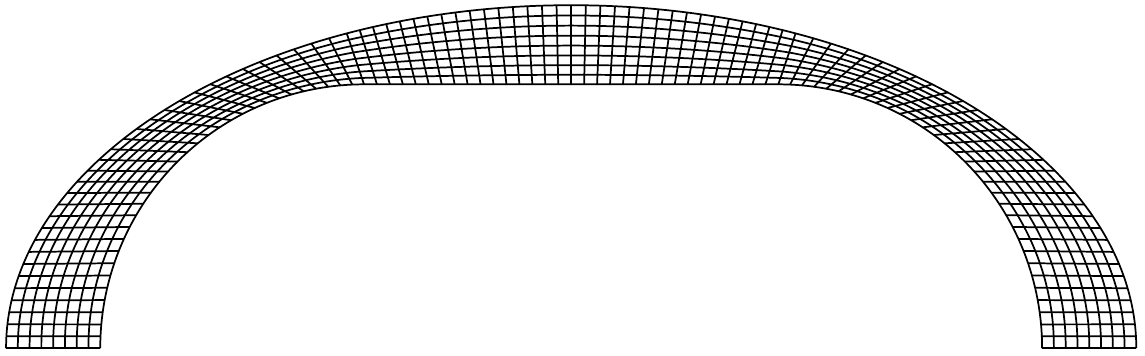}
		\caption{Mesh 6.}
	\end{subfigure}
	\caption{\textbf{Radial pulsation from a mock shell}: Meshes for the fluid domain between the scatterer and the artificial boundary. The meshes are constructed from the initial mesh 1, which is rotated around the axis of symmetry using four elements. Mesh 2 and 3 are refined only in the angular direction, while the more refined meshes also refine in the radial direction to obtain smallest aspect ratio. The meshes are nested.}
	\label{Fig2:MS_meshes}
\end{figure}

\begin{figure}
	\centering
	\includegraphics[width=\textwidth]{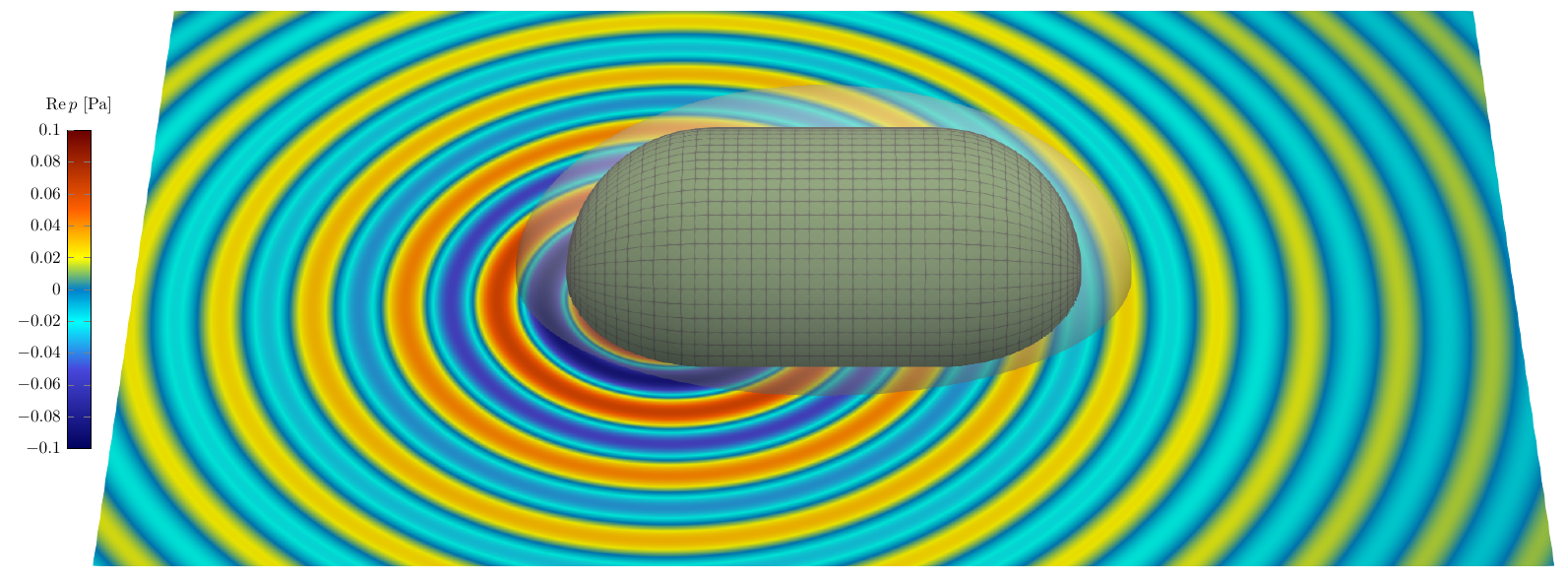}
	\caption{\textbf{Radial pulsation from a mock shell}: Visualization of numerical solution in the $xy$-plane using BGU with $N=6$ on mesh 5.}
	\label{Fig2:MS_visualization}
\end{figure}

\begin{figure}
	\centering    
	\begin{subfigure}{0.49\textwidth}
		\centering
		\includegraphics[width=\textwidth]{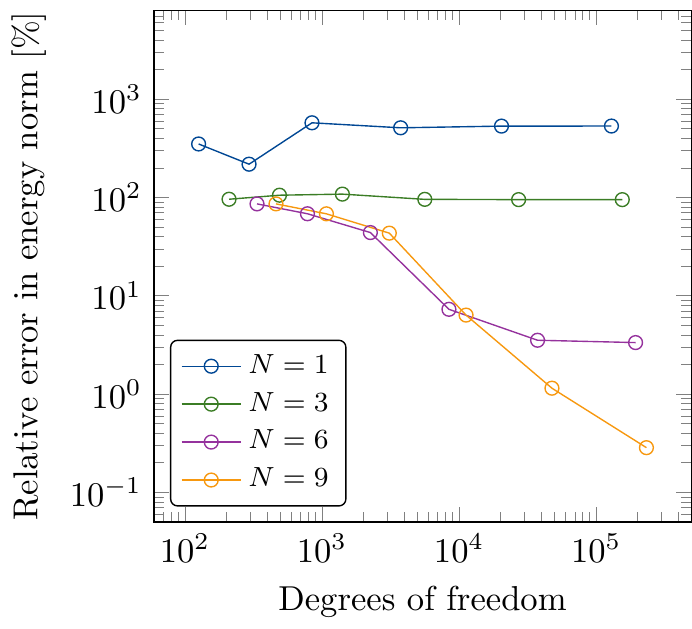}
	\caption{BGC}
	\end{subfigure}
	~    
	\begin{subfigure}{0.49\textwidth}
		\centering
		\includegraphics[width=\textwidth]{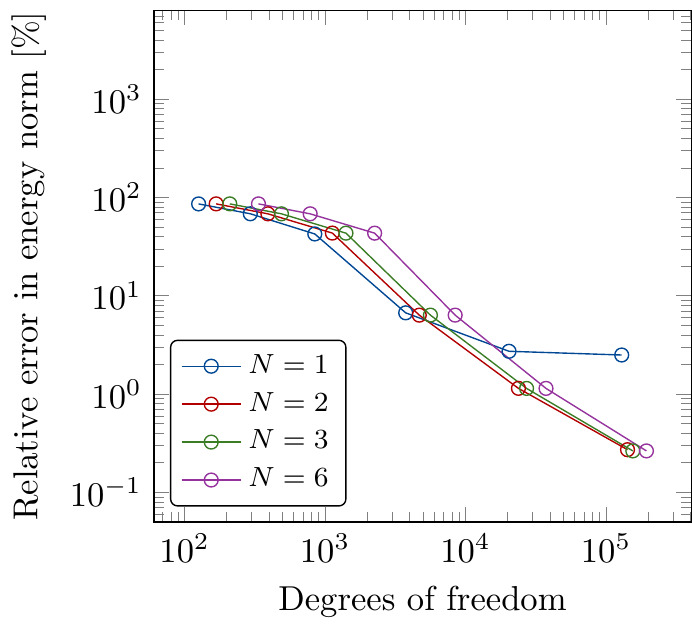}
		\caption{BGU}
	\end{subfigure}
	\par\bigskip
	\par\bigskip
	\begin{subfigure}{0.49\textwidth}
		\centering
		\includegraphics[width=\textwidth]{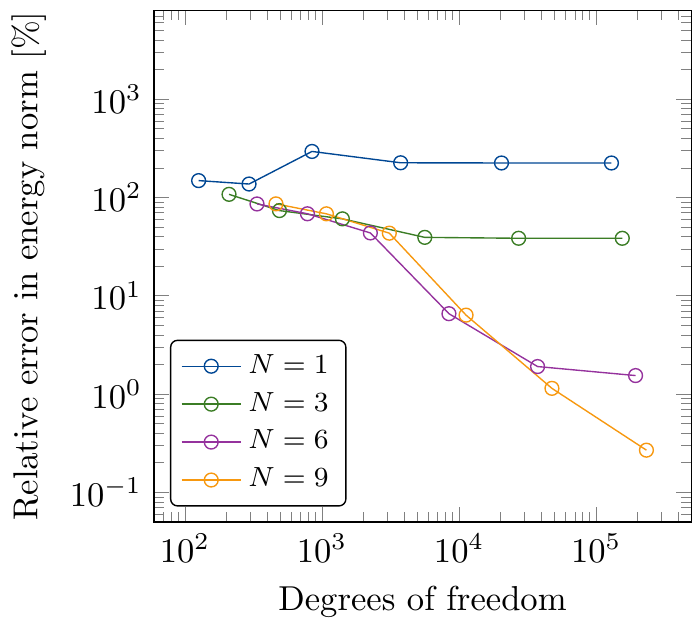}
		\caption{PGC}
	\end{subfigure}
	~    
	\begin{subfigure}{0.49\textwidth}
		\centering
		\includegraphics[width=\textwidth]{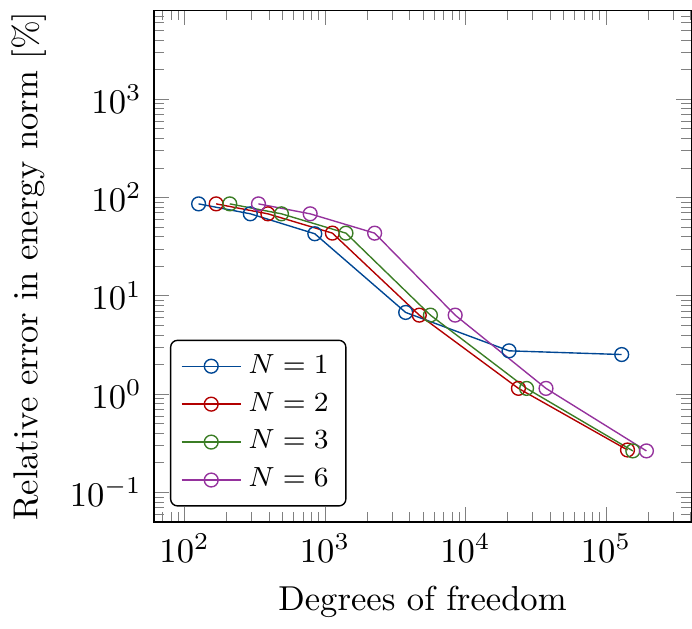}
		\caption{PGU}
	\end{subfigure}
	\caption{\textbf{Radial pulsation from a mock shell}: Convergence plots for the four infinite element formulations. The relative error in the energy norm (\Cref{Eq2:energyNormFluids}) is plotted against the number of degrees of freedom.}
	\label{Fig2:MS_convergence}
\end{figure}

This example is used to illustrate the differences of the infinite element formulations using the prolate ellipsoidal elements after Burnett~\cite{Burnett1994atd}. The mesh construction is illustrated in \Cref{Fig2:MS_meshes}, and an illustration of the solution is presented in \Cref{Fig2:MS_visualization}.
Convergence plots are shown in \Cref{Fig2:MS_convergence}. Gerdes did a similar comparison in~\cite{Gerdes1998tcv} where scattering on a sphere was investigated. Our results verify these findings, namely lower errors for the unconjugated formulations (cf. \Cref{Fig2:MS_convergence}). Good results can be obtained using only a single radial shape function in case of unconjugated formulations. For the conjugated versions, on the other hand, $N > 6$ functions are needed to obtain similar accuracy and more degrees of freedom are required to get an asymptotic behavior.

\begin{figure}
	\begin{subfigure}{0.49\textwidth}
		\centering
		\includegraphics[width=0.85\textwidth]{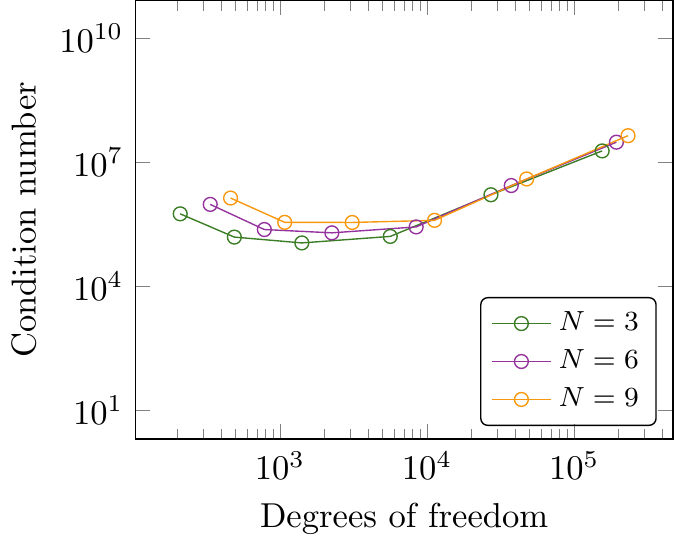}
		\caption{BGC with the shifted Chebyshev basis}
	\end{subfigure}
	~    
	\begin{subfigure}{0.49\textwidth}
		\centering
		\includegraphics[width=0.85\textwidth]{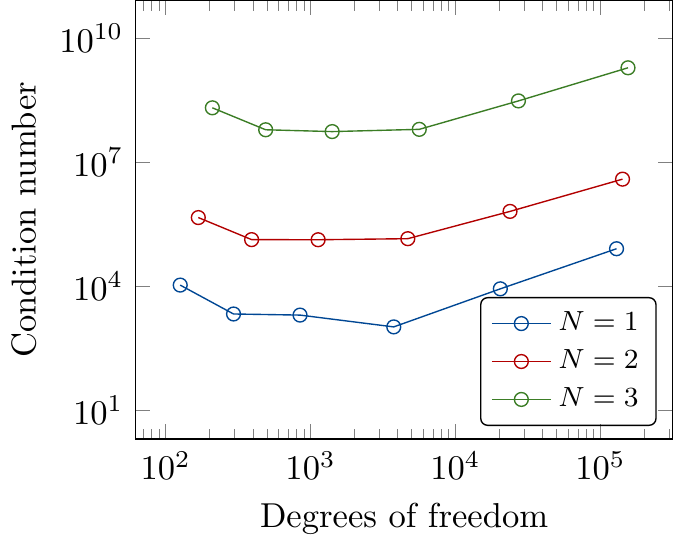}
		\caption{BGU with the shifted Chebyshev basis}
	\end{subfigure}
	\par\bigskip
	\par\bigskip
	\begin{subfigure}{0.49\textwidth}
		\centering
		\includegraphics[width=0.85\textwidth]{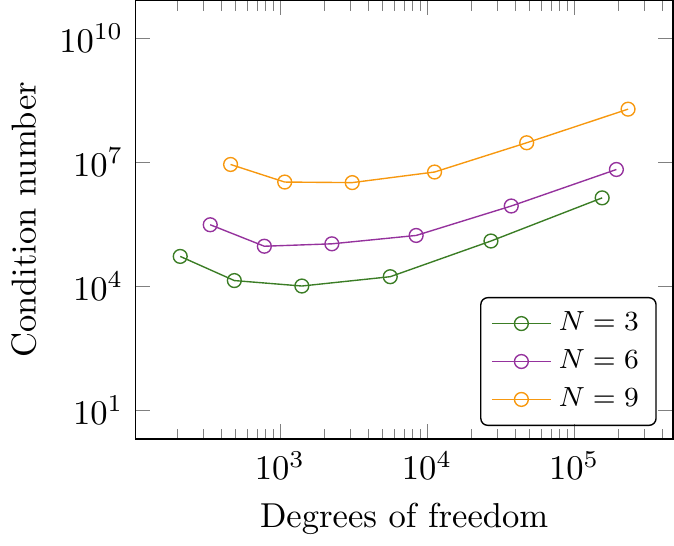}
		\caption{BGC with the Bernstein basis}
	\end{subfigure}
	~    
	\begin{subfigure}{0.49\textwidth}
		\centering
		\includegraphics[width=0.85\textwidth]{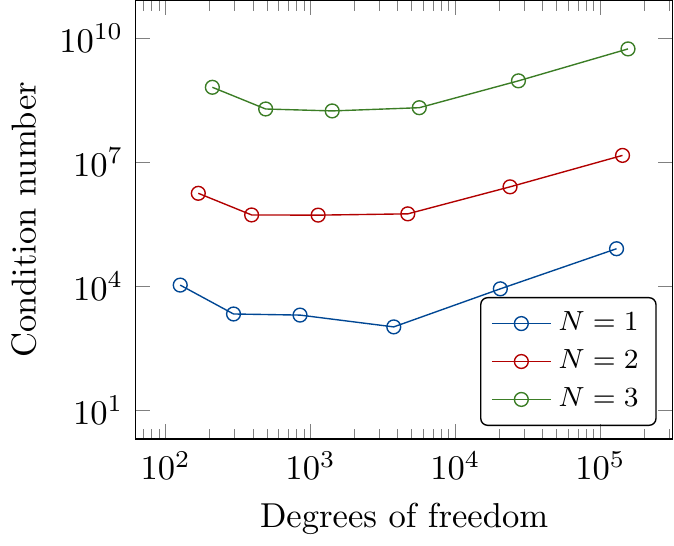}
		\caption{BGU with the Bernstein basis}
	\end{subfigure}
	\par\bigskip
	\par\bigskip
	\begin{subfigure}{0.49\textwidth}
		\centering
		\includegraphics[width=0.85\textwidth]{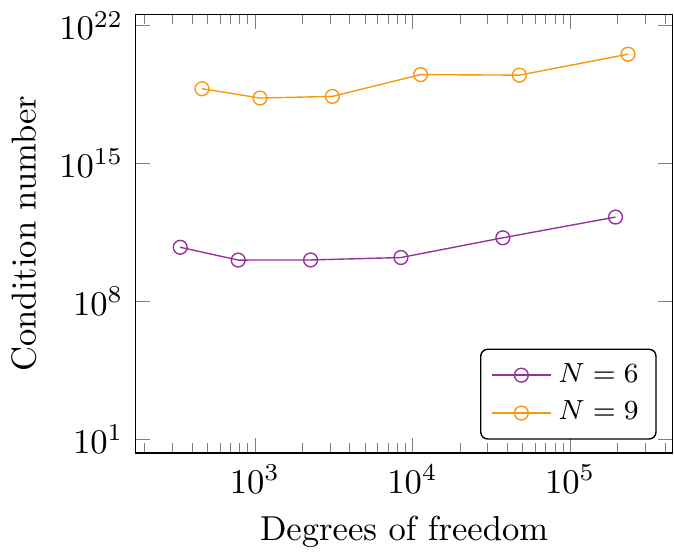}
		\caption{BGC with the Lagrange basis}
	\end{subfigure}
	~    
	\begin{subfigure}{0.49\textwidth}
		\centering
		\includegraphics[width=0.85\textwidth]{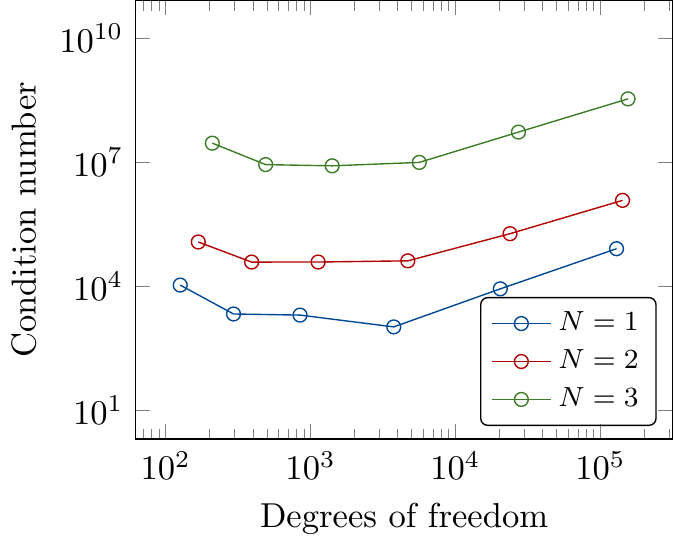}
		\caption{BGU with the Lagrange basis}
	\end{subfigure}
	\caption{\textbf{Radial pulsation from a mock shell}: Convergence plots for the BGC and BGU formulations using three different sets of radial shape functions (Chebyshev, Bernstein and Lagrange). The condition number (1-norm condition number estimate provided by \texttt{condest} in MATLAB) is plotted against the number of degrees of freedom.}
	\label{Fig2:MS_condNumbersB}
\end{figure}

\begin{figure}
	\begin{subfigure}{0.49\textwidth}
		\centering
		\includegraphics[width=0.85\textwidth]{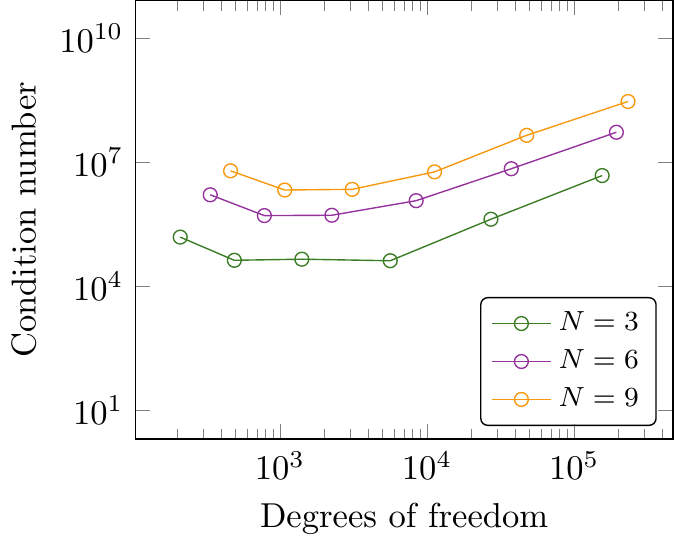}
		\caption{PGC with the shifted Chebyshev basis}
	\end{subfigure}
	~    
	\begin{subfigure}{0.49\textwidth}
		\centering
		\includegraphics[width=0.85\textwidth]{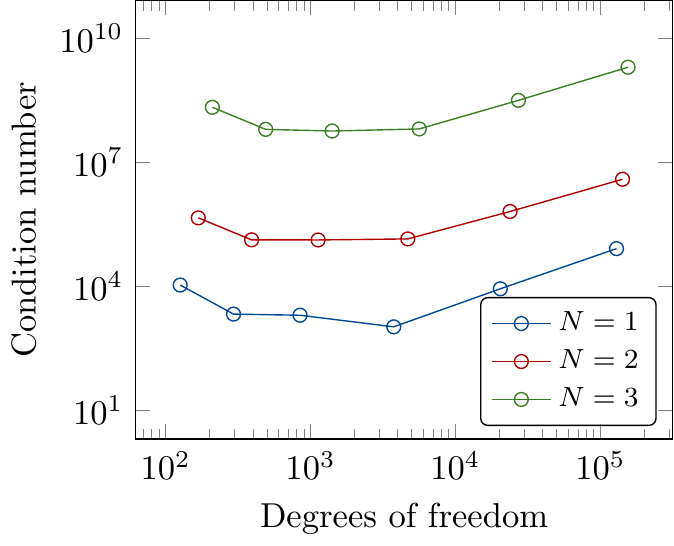}
		\caption{PGU with the shifted Chebyshev basis}
	\end{subfigure}
	\par\bigskip
	\par\bigskip
	\begin{subfigure}{0.49\textwidth}
		\centering
		\includegraphics[width=0.85\textwidth]{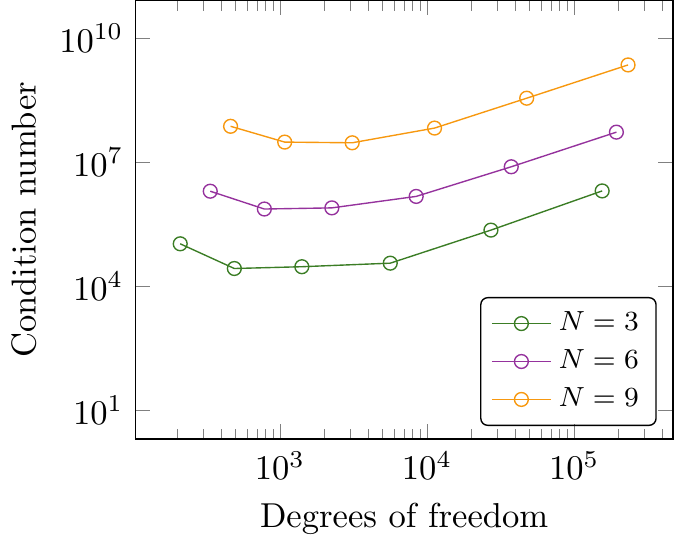}
		\caption{PGC with the Bernstein basis}
	\end{subfigure}
	~    
	\begin{subfigure}{0.49\textwidth}
		\centering
		\includegraphics[width=0.85\textwidth]{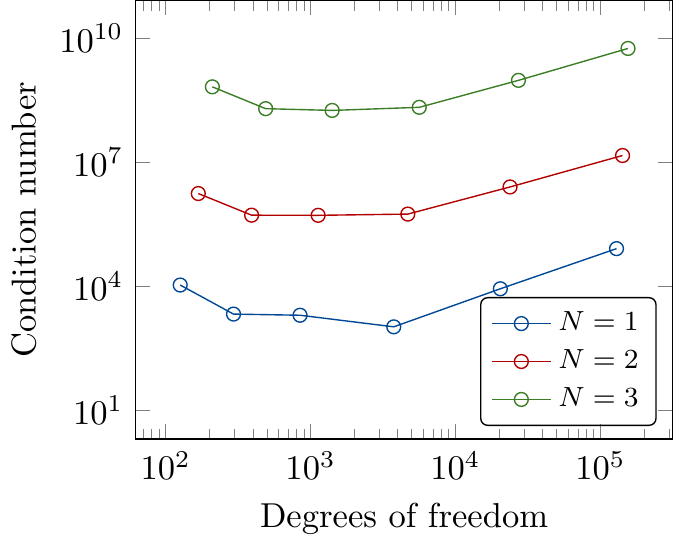}
		\caption{PGU with the Bernstein basis}
	\end{subfigure}
	\par\bigskip
	\par\bigskip
	\begin{subfigure}{0.49\textwidth}
		\centering
		\includegraphics[width=0.85\textwidth]{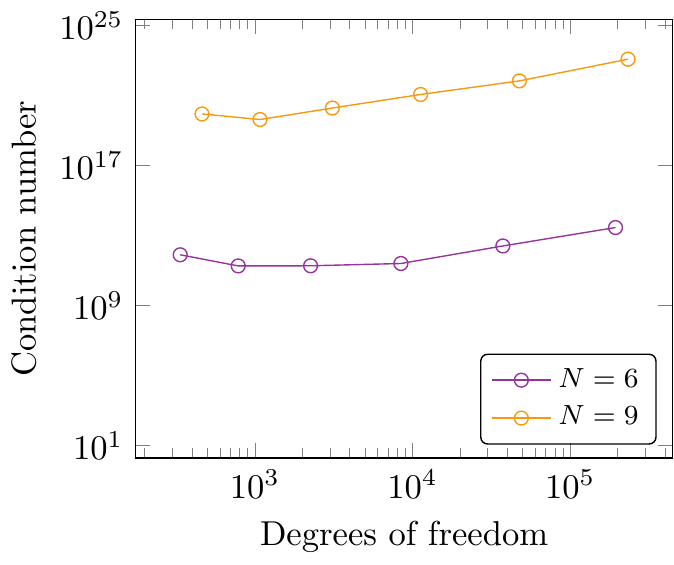}
		\caption{PGC with the Lagrange basis}
	\end{subfigure}
	~    
	\begin{subfigure}{0.49\textwidth}
		\centering
		\includegraphics[width=0.85\textwidth]{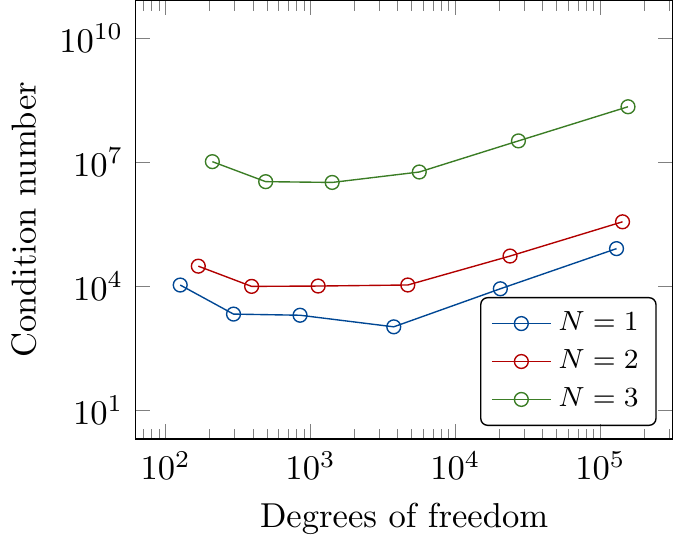}
		\caption{PGU with the Lagrange basis}
	\end{subfigure}
	\caption{\textbf{Radial pulsation from a mock shell}: Convergence plots for the PGC and PGU formulations using three different sets of radial shape functions (Chebyshev, Bernstein and Lagrange). The condition number (1-norm condition number estimate provided by \texttt{condest} in MATLAB) is plotted against the number of degrees of freedom.}
	\label{Fig2:MS_condNumbersP}
\end{figure}

In \Cref{Fig2:MS_condNumbersB} and \Cref{Fig2:MS_condNumbersP} the condition number is investigated for the different formulations and basis functions in the radial shape functions. The condition number for the unconjugated versions increases more rapidly as a function of $N$ compared to the corresponding formulations in the conjugated case. The condition number of the Lagrange basis increases particularly fast with $N$, making it useless\footnote{In the case of $r_n = nr_{\mathrm{a}}$.} for the conjugated formulations. However, the Lagrange basis yields the best result for the unconjugated formulations for small $N$. The Chebyshev basis seems to give the best condition numbers for the conjugated formulations for large $N$ (which is required for acceptable results).
%			N=1		N=2		N=3		N=6		N=9
%BGU		same	LCB		LCB		CBL		BCL
%PGU		 "		 "		 "		 "		 "
%BGC		 "		BCL		BCL		BCL		CBL
%PGC		 "		 "		 "		 "		 "
The unconjugated formulations perform quite similar, both in terms of the condition numbers and the error. The BGU formulation has the additional advantage of producing symmetric matrices, and reduces the memory requirement. It is clear that the choice of basis functions in the infinite elements plays a crucial role for the condition number, and more research is required to find the optimal set of basis functions. Based on the findings in this work, it is recommended to use the BGU formulation alongside the Lagrange basis (in the radial direction) in the infinite elements. However, if larger $N$ is needed for accuracy, the Chebyshev basis is recommended.

\subsection{Stripped BeTSSi submarine}
\begin{figure}
	\begin{subfigure}{\textwidth}
		\centering
		\includegraphics[width=0.95\textwidth]{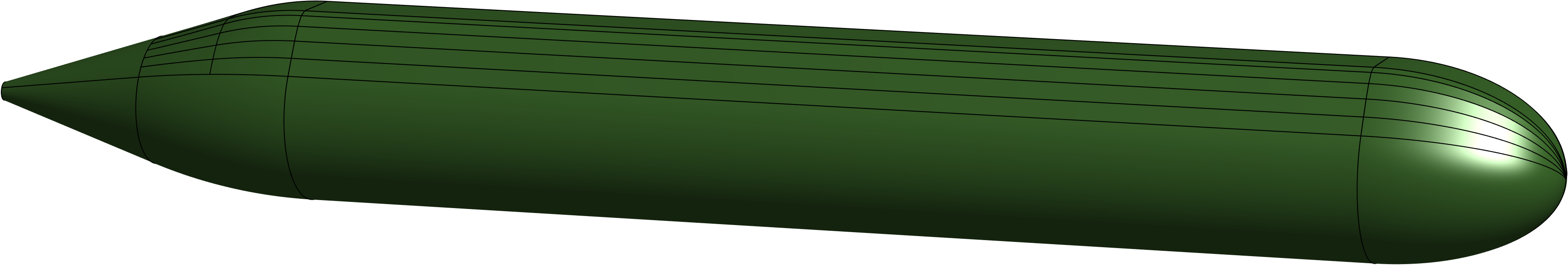}
		\caption{CAD model.}
		\label{Fig2:BC_stripped}
	\end{subfigure}
	\par\bigskip
	\par\bigskip
	\begin{subfigure}{\textwidth}
		\centering
		\includegraphics[width=0.95\textwidth]{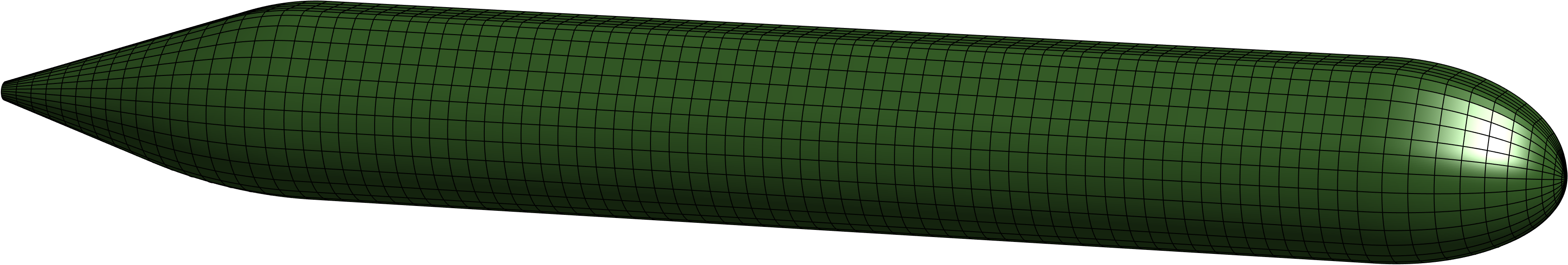}
		\caption{Surface mesh for mesh ${\cal M}_{1,\hat{p},\hat{k}}^{\mathrm{IGA}}$}
	\end{subfigure}	
	\par\bigskip
	\par\bigskip
	\begin{subfigure}{\textwidth}
		\centering
		\includegraphics[width=0.95\textwidth]{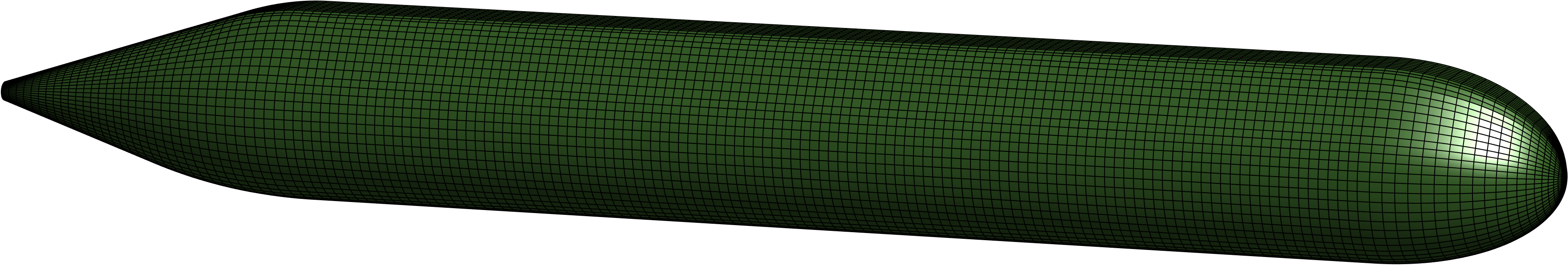}
		\caption{Surface mesh for mesh ${\cal M}_{2,\hat{p},\hat{k}}^{\mathrm{IGA}}$.}
	\end{subfigure}	
	\par\bigskip
	\par\bigskip
	\begin{subfigure}{\textwidth}
		\centering
		\includegraphics[width=0.95\textwidth]{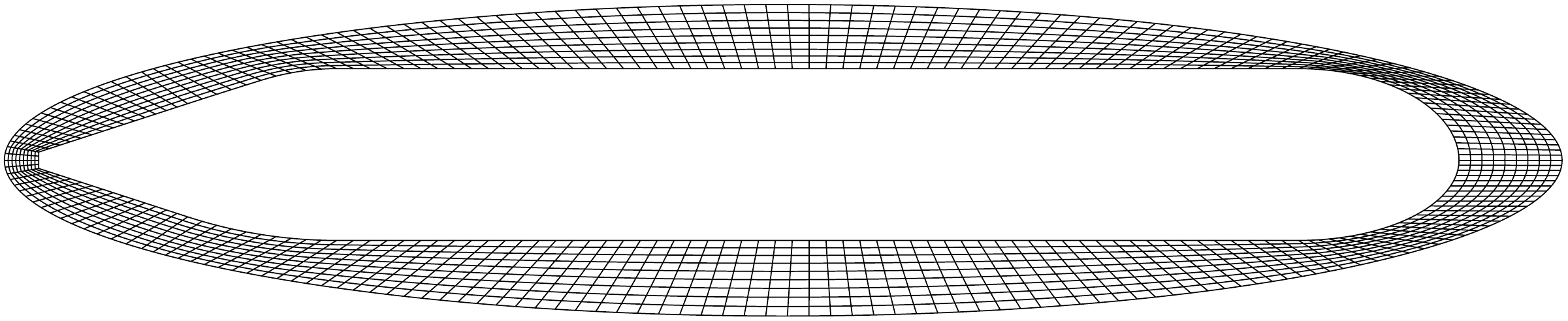}
		\caption{Crossection in the $xz$-plane for mesh ${\cal M}_{1,\hat{p},\hat{k}}^{\mathrm{IGA}}$.}
	\end{subfigure}	
	\par\bigskip
	\par\bigskip
	\begin{subfigure}{\textwidth}
		\centering
		\includegraphics[width=0.95\textwidth]{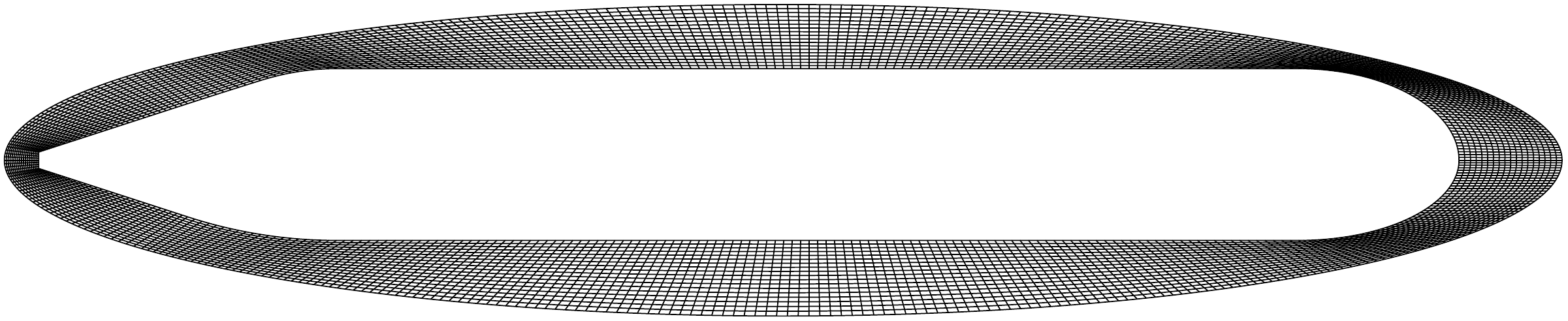}
		\caption{Crossection in the $xz$-plane for mesh ${\cal M}_{2,\hat{p},\hat{k}}^{\mathrm{IGA}}$.}
	\end{subfigure}	
	\caption{\textbf{Stripped BeTSSi submarine}: CAD model and meshes used for computations.}
	\label{Fig2:BC_strippedAndMesh}
\end{figure}

Finally, we consider the \textit{stripped BeTSSi submarine}\footnote{Based upon the BeTSSi submarine which originates from the BeTSSi workshops~\cite{Gilroy2013bib}.} described in \Cref{Sec2:BeTSSi_description}, and let a plane wave, with the direction of incidence given by
\begin{equation}
	\vec{d}_{\mathrm{s}} = -\begin{bmatrix}
		\cos\beta_{\mathrm{s}}\cos\alpha_{\mathrm{s}}\\
		\cos\beta_{\mathrm{s}}\sin\alpha_{\mathrm{s}}\\
		\sin\beta_{\mathrm{s}}
	\end{bmatrix}, \quad\text{where}\quad \alpha_{\mathrm{s}} = \ang{240},\,\beta_{\mathrm{s}} = \ang{0},
\end{equation}
be scattered by this submarine. The CAD model is given in \Cref{Fig2:BC_strippedAndMesh} alongside computational meshes. Again, we shall denote by ${\cal M}_{m,\check{p},\check{k}}^{\textsc{iga}}$, mesh number $m$ with polynomial order $\check{p}$ and continuity $\check{k}$ across element boundaries of the NURBS parametrization. 

\begin{figure}
	\centering    
	\begin{subfigure}[b]{\textwidth}
		\centering
		\includegraphics[width=\textwidth]{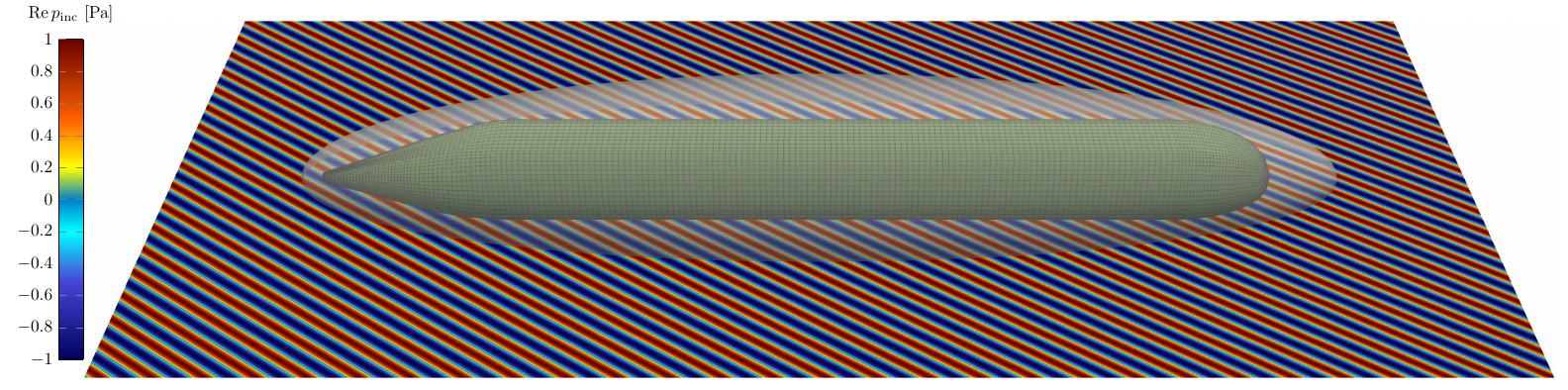}
		\caption{Real part of the incident wave $p_{\mathrm{inc}}(\vec{x})=P_{\mathrm{inc}}\euler^{\imag k\vec{d}_{\mathrm{s}}\cdot \vec{x}}$.}
	\end{subfigure}
	\par\bigskip
	\begin{subfigure}[b]{\textwidth}
		\centering
		\includegraphics[width=\textwidth]{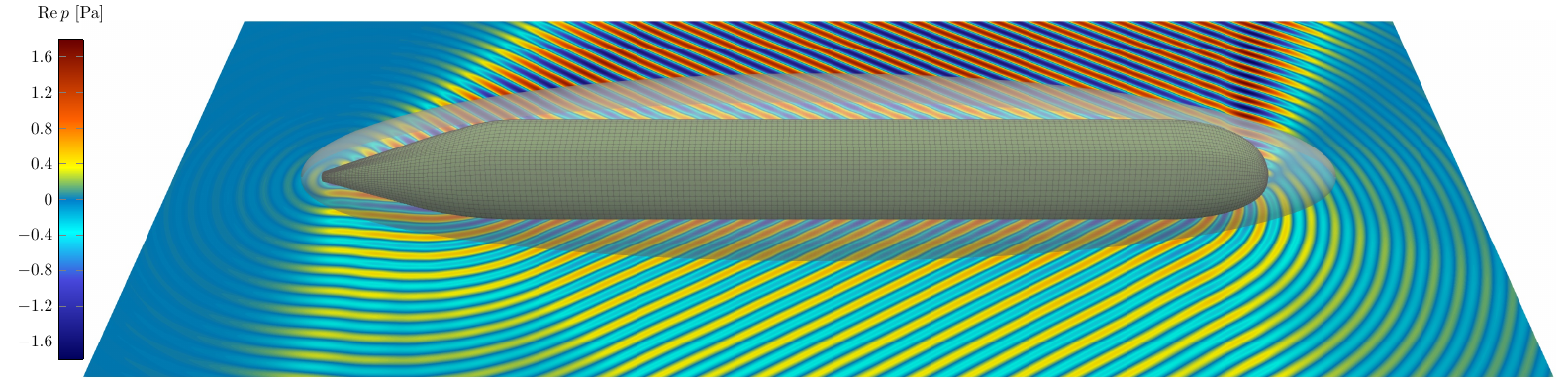}
		\caption{Real part of the scattered pressure $p(\vec{x})$.}
	\end{subfigure}
	\par\bigskip
	\begin{subfigure}[b]{\textwidth}
		\centering
		\includegraphics[width=\textwidth]{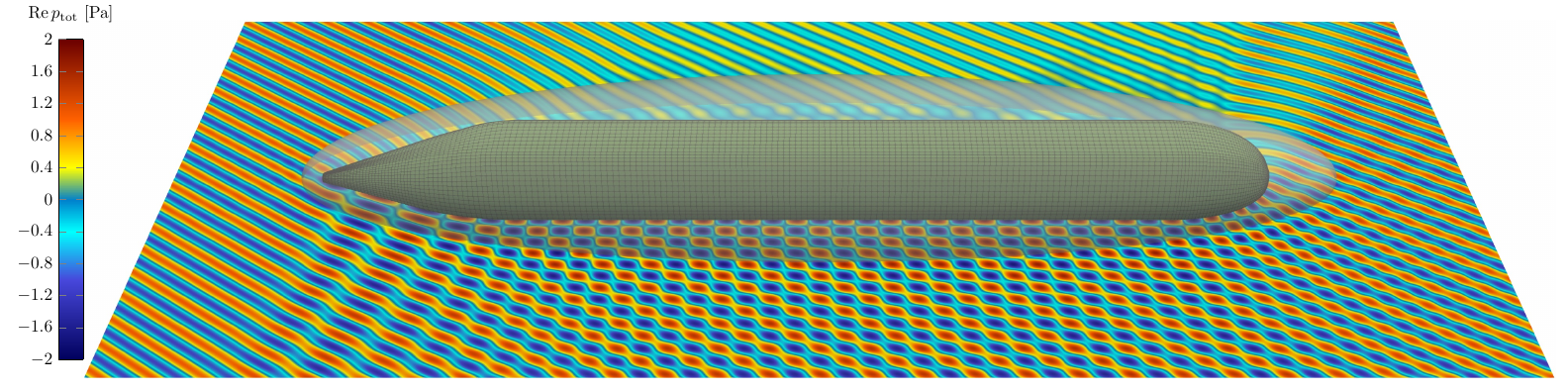}
		\caption{Real part of the total pressure $p_{\mathrm{tot}}(\vec{x})=p_{\mathrm{inc}}(\vec{x})+p(\vec{x})$.}
	\end{subfigure}
	\par\bigskip
	\begin{subfigure}[b]{\textwidth}
		\centering
		\includegraphics[width=\textwidth]{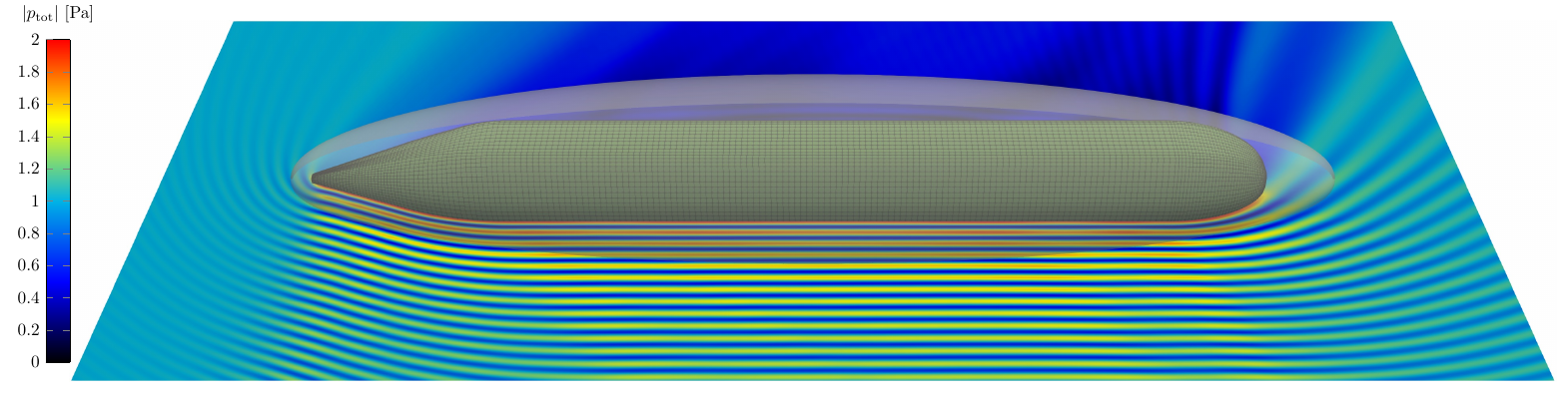}
		\caption{Modulus of the total pressure $p_{\mathrm{tot}}(\vec{x})=p_{\mathrm{inc}}(\vec{x})+p(\vec{x})$.}
	\end{subfigure}
	\caption{\textbf{Stripped BeTSSi submarine with SHBC}: The simulation at $f=\SI{1000}{Hz}$ is visualized in the $xy$-plane, and is computed on mesh ${\cal M}_{2,3,2}^{\mathrm{IGA}}$ and the BGU formulation with $N=4$. The numerical evaluations outside the (transparent) prolate ellipsoidal artificial boundary are evaluated with \Cref{Eq2:KirchhoffIntegral}.}
	\label{Fig2:BC_NearField}
\end{figure}

\begin{figure}
	\centering    
	\includegraphics[width=\textwidth]{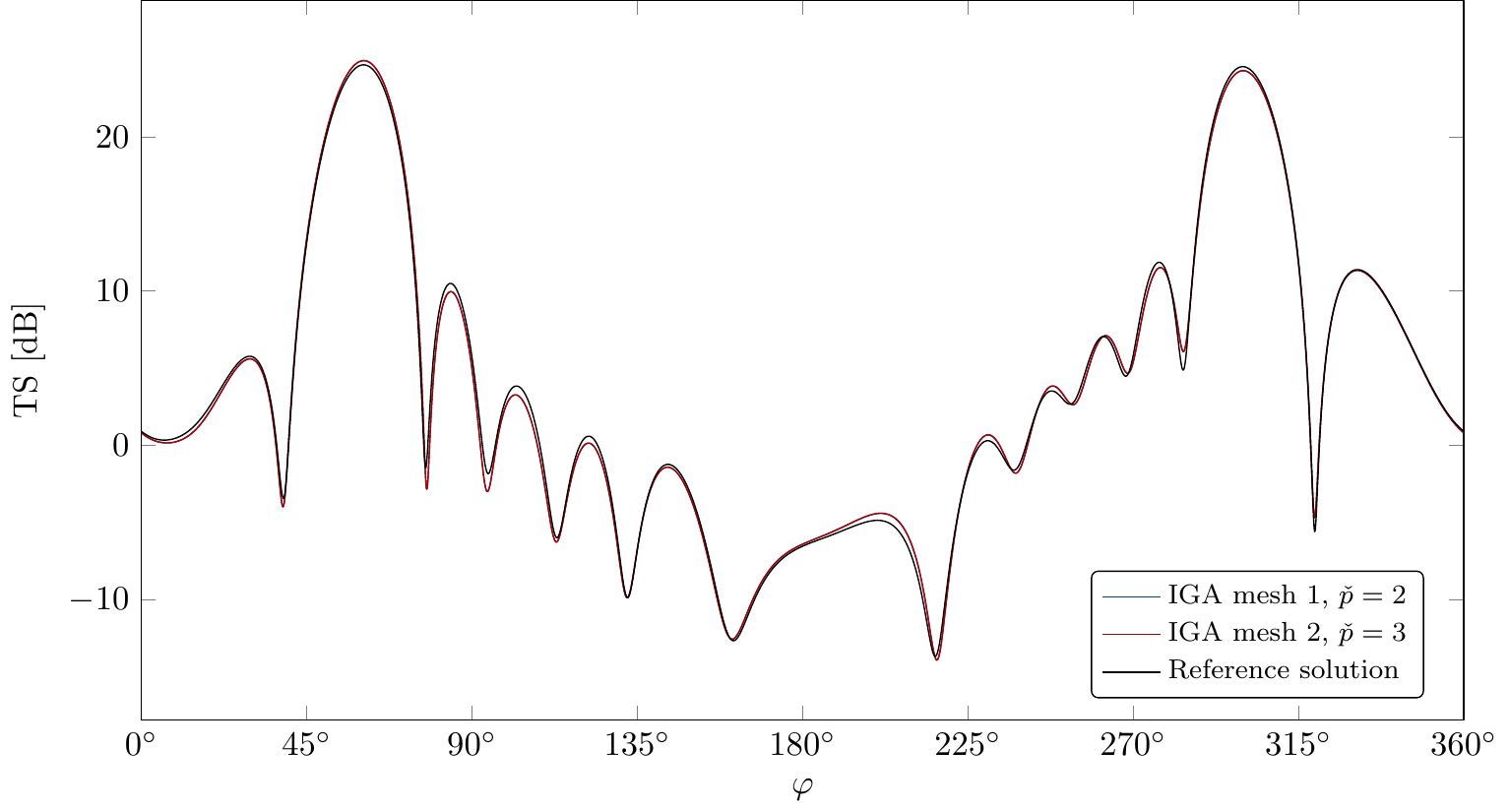}
	\caption{\textbf{Stripped BeTSSi submarine with SHBC}: Computation of target strength (\Cref{Eq2:TS}) at $f=\SI{100}{Hz}$ as a function of the azimuth angle in the spherical coordinate system. The two IGA results (both using $N=3$) are visually indistinguishable meaning that mesh ${\cal M}_{1,2,1}^{\mathrm{IGA}}$ is well resolved for this frequency.}
	\label{Fig2:FarField100}
\end{figure}

\begin{figure}
	\centering    
	\includegraphics[width=\textwidth]{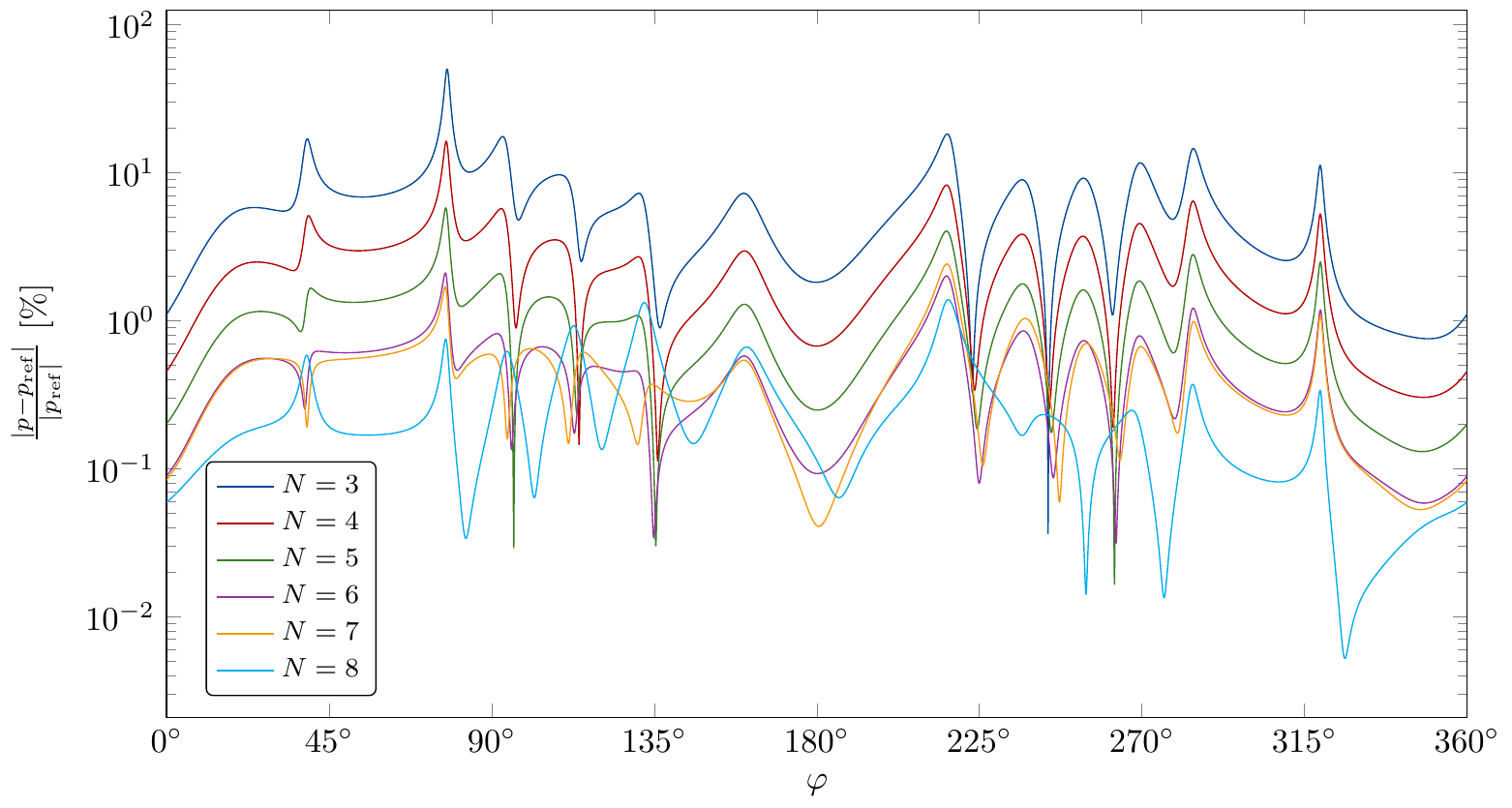}
	\caption{\textbf{Stripped BeTSSi submarine with SHBC}: Computation of the relative error in the far field (\Cref{Eq2:KirchhoffIntegral}) compared to a reference solution at $f=\SI{100}{Hz}$. The computations are done using IGA on mesh ${\cal M}_{2,3,2}^{\mathrm{IGA}}$ using the BGU formulation.}
	\label{Fig2:FarField100error}
\end{figure}

\begin{figure}
	\centering    
	\begin{subfigure}[b]{\textwidth}
		\centering
		\includegraphics[width=\textwidth]{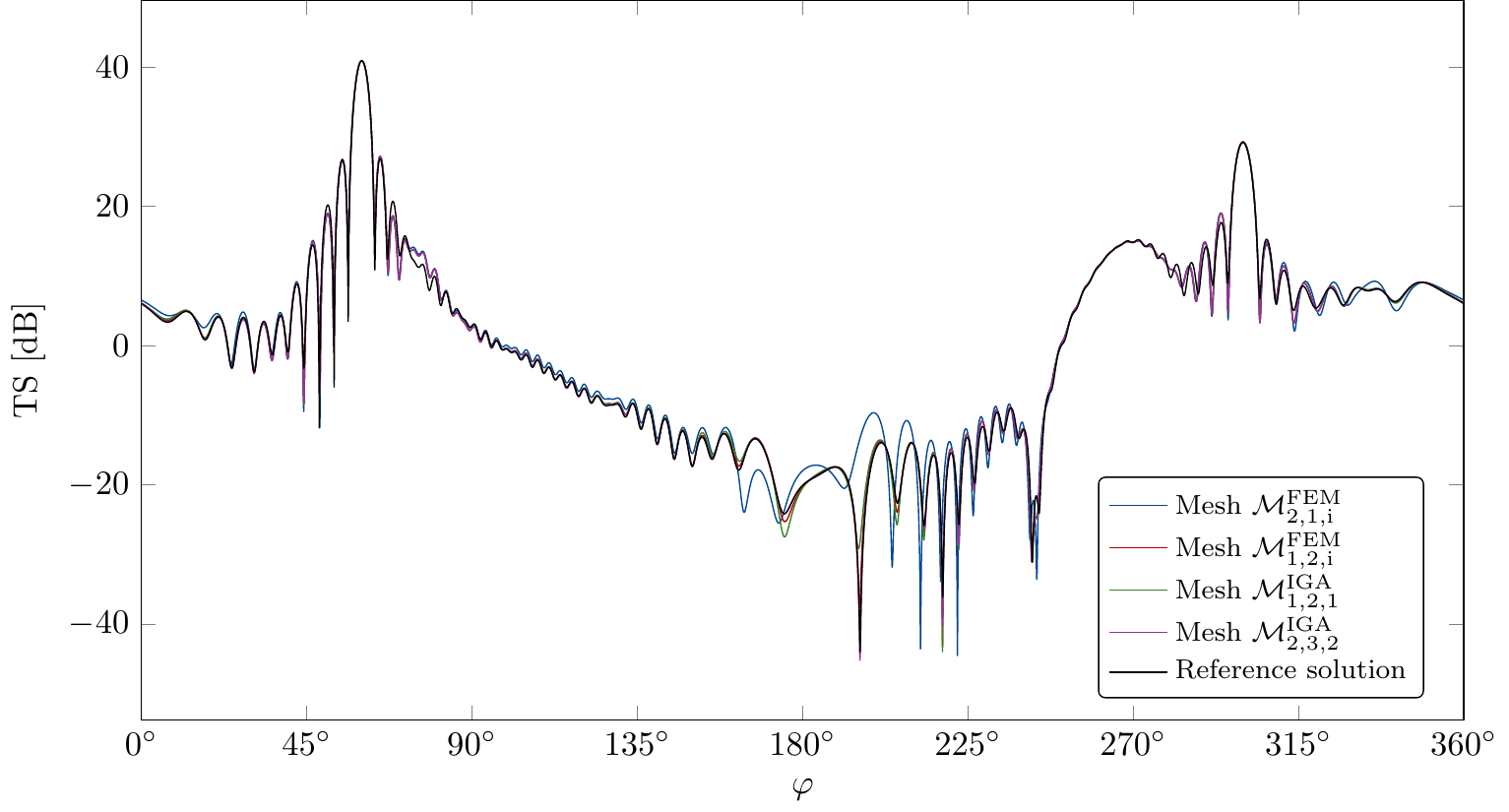}
		\caption{$f=\SI{500}{Hz}$}
	\end{subfigure}
	\par\bigskip
	\par\bigskip
	\begin{subfigure}[b]{\textwidth}
		\centering
		\includegraphics[width=\textwidth]{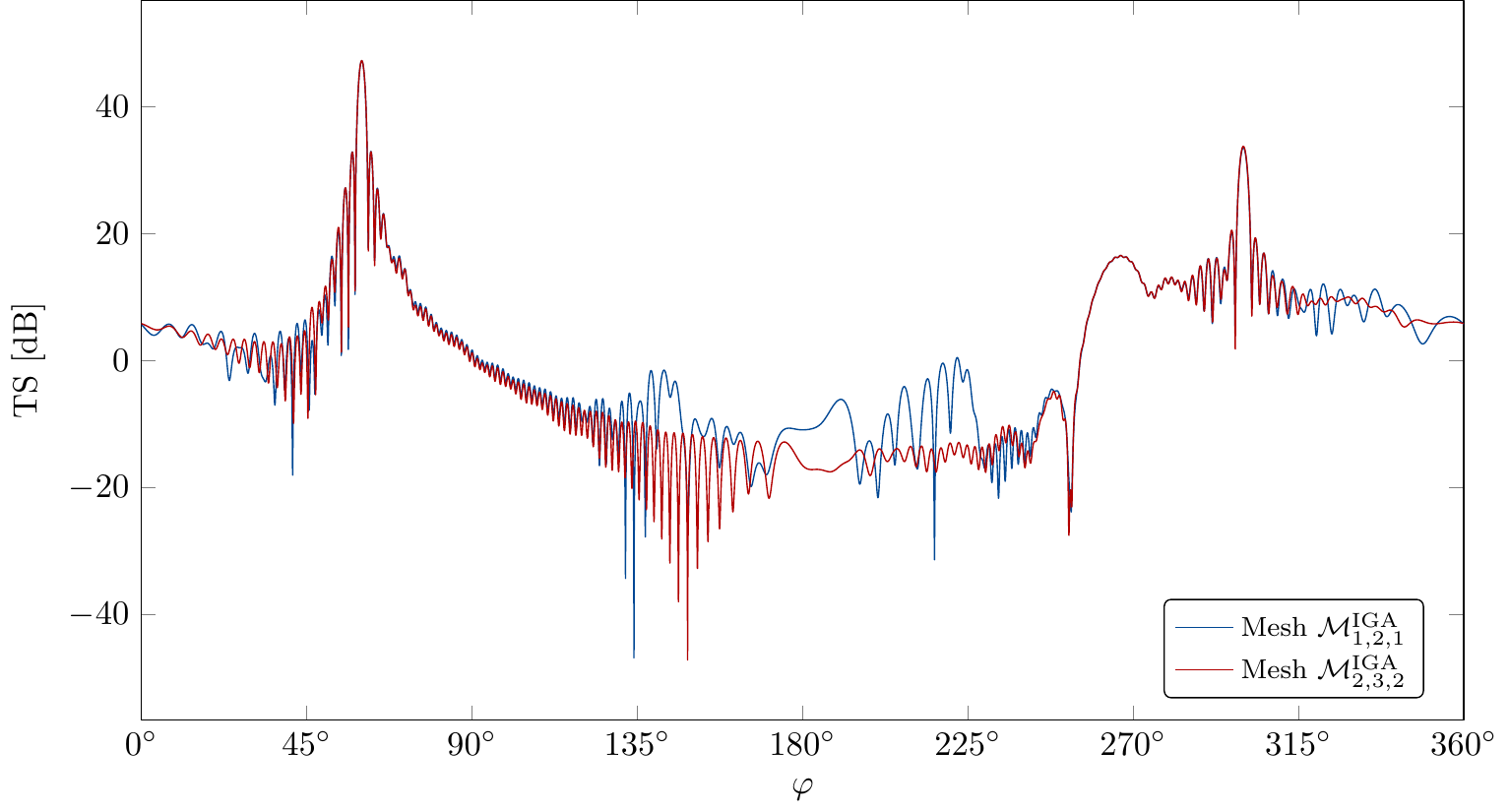}
		\caption{$f=\SI{1000}{Hz}$}
	\end{subfigure}
	\caption{\textbf{Stripped BeTSSi submarine with SHBC}: Computation of target strength (\Cref{Eq2:TS}) as a function of the azimuth angle in the spherical coordinate system. The numerical evaluations are evaluated with \Cref{Eq2:TS}.}
	\label{Fig2:FarField}
\end{figure}

The near field at $f=\SI{1000}{Hz}$ is visualized in \Cref{Fig2:BC_NearField}.
The low frequency problem at $f = \SI{100}{Hz}$ is considered in~\Cref{Fig2:FarField100}. In this case, mesh ${\cal M}_{1,2,1}^{\mathrm{IGA}}$ resolves this frequency, but the solution slightly deviates from the reference solution computed by IGABEM on a fine mesh. The reason for this is that $N$ is too low. Although $N=3$ was enough for engineering precision (below 1\%) in the mock shell example, it does not suffice for the more complicated geometry like the stripped BeTSSi submarine. Consider the relative error for the far field at the well resolved mesh ${\cal M}_{2,3,2}^{\mathrm{IGA}}$. In this case the error will originate from the low resolution (governed by $N$) in the radial direction for the infinite elements. As illustrated in \Cref{Fig2:FarField100error} an order of magnitude in accuracy is gained by increasing $N$. This effect was also observed by the verification test in \Cref{Subsec2:mockShell} applied to the stripped BeTSSi submarine.
In \Cref{Fig2:FarField} the target strength is plotted for $f=\SI{500}{Hz}$ and $f=\SI{1000}{Hz}$. A reference solution (using IGABEM) is added for the $f=\SI{500}{Hz}$ case, and illustrates again the pollution of low $N$. The IGA mesh 1 resolves the frequency $f=\SI{500}{Hz}$ quite well using only about 5 elements per wave length. This corresponds to about 5 dofs per wave length in each dimensional direction compared to the classical 10-12 dofs per wave length needed for FEM methods.

%% file: contents/IGAFEM_SHBC.tex
Mesh ${\cal M}_{6,1,\mathrm{i}}^{\textsc{fem}}$		& 32768	& 56462	& 7.70	& 26.50	& 5.04	\cr
Mesh ${\cal M}_{5,2,\mathrm{i}}^{\textsc{fem}}$		& 4096	& 56462	& 5.07	& 20.48	& 0.62	\cr
Mesh ${\cal M}_{5,2,1}^{\textsc{iga}}$				& 4096	& 13476	& 4.79	& 4.32	& 0.64	\cr
Mesh ${\cal M}_{4,3,2}^{\textsc{iga}}$				& 512	& 4572	& 3.06	& 1.64	& 0.38	\cr
Mesh ${\cal M}_{5,3,2}^{\textsc{iga}}$				& 4096	& 17654	& 14.75	& 8.30	& 0.05	\cr

%% file: contents/IGAFEM_SSBC.tex
Mesh ${\cal M}_{6,1,\mathrm{i}}^{\textsc{fem}}$		& 40960	& 104858	& 13.66	& 75.28	& 7.66	\cr
Mesh ${\cal M}_{5,2,\mathrm{i}}^{\textsc{fem}}$		& 6144	& 129056	& 15.85	& 106.07	& 1.35	\cr
Mesh ${\cal M}_{5,2,1}^{\textsc{iga}}$				& 6144	& 33690	& 13.35	& 34.68	& 0.99	\cr
Mesh ${\cal M}_{4,3,2}^{\textsc{iga}}$				& 1024	& 13716	& 11.83	& 9.04	& 41.30	\cr
Mesh ${\cal M}_{5,3,2}^{\textsc{iga}}$				& 6144	& 47918	& 52.03	& 69.94	& 0.09	\cr

%% file: contents/IGAFEM_NNBC.tex
Mesh ${\cal M}_{6,1,\mathrm{i}}^{\textsc{fem}}$		& 172032	& 233915	& 27.80	& 429.52	& 6.55	\cr
Mesh ${\cal M}_{5,2,\mathrm{i}}^{\textsc{fem}}$		& 22528	& 258113	& 31.05	& 462.66	& 0.53	\cr
Mesh ${\cal M}_{5,2,1}^{\textsc{iga}}$				& 22528	& 53905	& 22.46	& 64.91	& 0.71	\cr
Mesh ${\cal M}_{4,3,2}^{\textsc{iga}}$				& 3072	& 18289	& 17.10	& 17.28	& 1.47	\cr
Mesh ${\cal M}_{5,3,2}^{\textsc{iga}}$				& 22528	& 73139	& 93.22	& 145.74	& 0.05	\cr

%% file: contents/conclusion.tex
\clearpage
\section{Conclusions}
\label{Sec2:conclusions}
This article addresses acoustic scattering characterized by sound waves reflected by man-made elastic objects. The present approach is characterized by:
\begin{itemize}
	\item The fluid surrounding (inside and in the vicinity outside) the solid scatterer is discretized by using isogeometric analysis (IGA).
	\item The unbounded domain outside the artificial boundary circumscribing the scatterer is handled by use of the infinite element method (IEM).
	\item The elastic scatterer is discretized by using IGA.
	\item The coupled acoustic structure interaction (ASI) problem is solved as a monolithic problem.
\end{itemize}

The main finding of the present study is that the use of IGA significantly increases the accuracy compared to the use of $C^0$ finite element analysis (FEA) due to increased inter-element continuity of the spline basis functions.

Furthermore, the following observations are made
\begin{itemize}
	\item IGA and the four presented IEM formulations work well on acoustic scattering for low frequencies. Among the infinite element formulations, the unconjugated version seems to give the best results.
	\item IGA's ability to represent the geometry exactly was observed to be of less importance for accuracy when comparing to higher order ($\hat{p}\geq 2$) isoparametric FEA. However, a more significant improvement offered by IGA is due to higher continuity of the spline basis functions in the solution space.
	\item The IGA framework enables roughly the same accuracy per element (compared to higher order isoparametric FEA) even though the number of degrees of freedom is significantly reduced.
	\item IGA is more computationally efficient than FEA to obtain highly accurate solutions. That is, when the mesh is sufficiently resolved, a given accuracy is obtained computationally faster using IGA.
	\item As for the FEA, IGA also suffers from the pollution effect at high frequencies. This will always be a problem, and for the higher frequency spectrum, the methods must be extended correspondingly. The XIBEM~\cite{Peake2013eib,Peake2015eib} (extended isogeometric boundary element method) is such an extension for the boundary element method. This technique (and similar enrichment strategies) could be applied to IEM as well, and is suggested as future work.
	\item The IEM suffers from high condition numbers when the number of radial shape functions in the infinite elements ($N$) is large. This becomes a problem for more complex geometries as $N$ must be increased to achieve higher precision.
\end{itemize}
The main disadvantages of using IGA with IEM is the need for a surface-to-volume parametrization between the scatterer and the artificial boundary, $\Omega_{\mathrm{a}}$. In this paper, the scatterer has been simple enough to discretize $\Omega_{\mathrm{a}}$ using a single 3D NURBS patch. For more complex geometries, this becomes more involved, and is a topic of active research to this date in the IGA community~\cite{Engvall2016itb,Engvall2017iut,Xia2017iaw}. The surface-to-volume parametrization and the conditioning are the main open issues of IGA with IEM and should be explored in future research.

\section*{Acknowledgements}
This work was supported by the Department of Mathematical Sciences at the Norwegian University of Science and Technology and by the Norwegian Defence Research Establishment.

The stripped BeTSSi submarine simulations were performed on resources provided by UNINETT Sigma2 - the National Infrastructure
for High Performance Computing and Data Storage in Norway (reference number: NN9322K/4317).

The authors would like to thank the reviewers for detailed response and many constructive comments.

%% file: contents/infiniteElementsBilinearForm.tex
\clearpage
\section{Derivation of bilinear form in infinite elements}
\label{Sec2:AppendixDerivationOfBilinearForm}
In this appendix, the integrals in the bilinear forms for the infinite elements will be separated for the PGU case\footnote{The other three formulations has been derived in~\cite{Venas2015iao}. For the more general ellipsoidal coordinate system, refer to~\cite{Burnett1998aea}.}. For generality, the derivation is done in the prolate spheroidal coordinate system.

\subsection{The prolate spheroidal coordinate system}
\label{Sec2:prolateSphericalCoordinateSystem}
The prolate spheroidal coordinate system is an extension of the spherical coordinate system. It is defined by the relations
\begin{align}
	x &= \sqrt{r^2 - \Upsilon^2}\sin\vartheta\cos\varphi\\
	y &= \sqrt{r^2 - \Upsilon^2}\sin\vartheta\sin\varphi\\
	z &= r\cos\vartheta
\end{align}
with foci located at $z = \pm \Upsilon$ and $r\geq \Upsilon$. Note that the coordinate system reduces to the spherical coordinate system when $\Upsilon=0$. the following inverse formulas may be derived
\begin{align}\label{Eq2:XtoProl}
\begin{split}
	r &= \frac{1}{2}(d_1+d_2)\\
	\vartheta &= \arccos\left(\frac{z}{r}\right)\\
	\varphi &= \operatorname{atan2}(y,x)
\end{split}
\end{align}
where 
\begin{align*}
	d_1 &= d_1(x,y,z) = \sqrt{x^2+y^2+(z+\Upsilon)^2}\\
	d_2 &= d_2(x,y,z) = \sqrt{x^2+y^2+(z-\Upsilon)^2}
\end{align*}
and 
\begin{equation*}
	\operatorname{atan2}(y,x) = \begin{cases}
	\arctan(\frac{y}{x}) & \mbox{if } x > 0\\
	\arctan(\frac{y}{x}) + \PI & \mbox{if } x < 0 \mbox{ and } y \ge 0\\
	\arctan(\frac{y}{x}) - \PI & \mbox{if } x < 0 \mbox{ and } y < 0\\
	\frac{\PI}{2} & \mbox{if } x = 0 \mbox{ and } y > 0\\
	-\frac{\PI}{2} & \mbox{if } x = 0 \mbox{ and } y < 0\\
	\text{undefined} & \mbox{if } x = 0 \mbox{ and } y = 0.
	\end{cases}
\end{equation*}
The derivatives are found to be
\begin{equation}\label{Eq2:dXdProlateSphericalCoordinates}
\begin{alignedat}{4}
	\pderiv{x}{r} &= \frac{r\sin\vartheta\cos\varphi}{\sqrt{r^2-\Upsilon^2}},\qquad	&&\pderiv{y}{r} = \frac{r\sin\vartheta\sin\varphi}{\sqrt{r^2-\Upsilon^2}},\qquad	&&\pderiv{z}{r} = \cos\vartheta\\
	\pderiv{x}{\vartheta} &=\sqrt{r^2-\Upsilon^2}\cos\vartheta\cos\varphi,\qquad	&&\pderiv{y}{\vartheta} = \sqrt{r^2-\Upsilon^2}\cos\vartheta\sin\varphi,\qquad	&&\pderiv{z}{\vartheta} = -r\sin\vartheta\\
	\pderiv{x}{\varphi} &= -\sqrt{r^2-\Upsilon^2}\sin\vartheta\sin\varphi,\qquad	&&\pderiv{y}{\varphi} = \sqrt{r^2-\Upsilon^2}\sin\vartheta\cos\varphi,\qquad	&&\pderiv{z}{\varphi} =0
\end{alignedat}
\end{equation}
and
\begin{equation}\label{Eq2:dProlateSphericalCoordinatesdX}
\begin{alignedat}{4}
	\pderiv{r}{x} &= \frac{x(d_1+d_2)}{2d_1d_2},\qquad	&&\pderiv{r}{y} = \frac{y(d_1+d_2)}{2d_1d_2},\qquad	&&\pderiv{r}{z} = \frac{z(d_1+d_2)+\Upsilon(d_2-d_1)}{2d_1d_2}\\
	\pderiv{\vartheta}{x} &= \frac{xz}{d_1d_2\sqrt{r^2-z^2}},\qquad	&&\pderiv{\vartheta}{y} = \frac{yz}{d_1d_2\sqrt{r^2-z^2}},\qquad	&&\pderiv{\vartheta}{z} = \frac{1}{\sqrt{r^2-z^2}}\left(\frac{z^2}{d_1d_2}+\frac{\Upsilon z(d_2-d_1)}{d_1d_2(d_1+d_2)}-1\right)\\
	\pderiv{\varphi}{x} &= -\frac{y}{x^2+y^2},\qquad	&&\pderiv{\varphi}{y} = \frac{x}{x^2+y^2},\qquad	&&\pderiv{\varphi}{z} = 0.
\end{alignedat}
\end{equation}
The general nabla operator can be written as
\begin{equation}\label{Eq2:GeneralNablaOperator}
	\nabla = \frac{\vec{e}_{\mathrm{r}}}{h_{\mathrm{r}}} \pderiv{}{r} + \frac{\vec{e}_{\upvartheta}}{h_{\upvartheta}} \pderiv{}{\vartheta} + \frac{\vec{e}_{\upvarphi}}{h_{\upvarphi}} \pderiv{}{\varphi}
\end{equation}
where
\begin{equation*}
	\vec{e}_{\mathrm{r}} = \frac{1}{h_{\mathrm{r}}}\left[\pderiv{x}{r}, \pderiv{y}{r}, \pderiv{z}{r}\right]^\transpose,\qquad
	\vec{e}_{\upvartheta} = \frac{1}{h_{\upvartheta}}\left[\pderiv{x}{\vartheta}, \pderiv{y}{\vartheta}, \pderiv{z}{\vartheta}\right]^\transpose,\qquad
	\vec{e}_{\upvarphi} = \frac{1}{h_{\upvarphi}}\left[\pderiv{x}{\varphi}, \pderiv{y}{\varphi}, \pderiv{z}{\varphi}\right]^\transpose
\end{equation*}
and
\begin{align*}
	h_{\mathrm{r}} &= \sqrt{\frac{r^2-\Upsilon^2\cos^2\vartheta}{r^2-\Upsilon^2}}\\
	h_{\upvartheta} &= \sqrt{r^2-\Upsilon^2\cos^2\vartheta}\\
	h_{\upvarphi} &= \sqrt{r^2-\Upsilon^2}\sin\vartheta.
\end{align*}
The Jacobian determinant (for the mapping from Cartesian coordinates to prolate spheroidal coordinates) may now be written as
\begin{equation}
	J_1 = h_{\mathrm{r}} h_{\upvartheta} h_{\upvarphi} = \left(r^2-\Upsilon^2\cos^2\vartheta\right)\sin\vartheta.
\end{equation}
As any normal vector at a surface with constant radius $r=\gamma$ can be written as $\vec{n} = \vec{e}_{\upvartheta}\times\vec{e}_{\upphi}=\vec{e}_{\mathrm{r}}$
\begin{equation}
	\partial_n p = \vec{n}\cdot\nabla p = \vec{e}_{\mathrm{r}}\cdot\nabla p = \frac{1}{h_{\mathrm{r}}} \pderiv{p}{r}.
\end{equation}
The surface Jacobian determinant at a given (constant) $r=\gamma$ is
\begin{equation}
	J_S = h_{\upvartheta} h_{\upvarphi} = \sqrt{r^2-\Upsilon^2\cos^2\vartheta}\sqrt{r^2-\Upsilon^2}\sin\vartheta,
\end{equation}
such that
\begin{equation}
	q\partial_n p J_S = \bigoh\left(r^{-3}\right)\quad\text{whenever}\quad q=\bigoh\left(r^{-3}\right)\quad\text{and}\quad p = \bigoh\left(r^{-1}\right).
\end{equation}
That is, for the Petrov--Galerkin formulations
\begin{equation}
	\lim_{\gamma\to\infty}\int_{S^\gamma} q\partial_n p\idiff\Gamma = \lim_{\gamma\to\infty}\int_0^{2\PI}\int_0^\PI q\partial_n p J_S\idiff\vartheta\idiff\varphi = 0.
\end{equation}

\subsection{Bilinear form for unconjugated Petrov--Galerkin formulation}
The bilinear form (in the domain outside the artificial boundary) in~\Cref{Eq2:B_uc_a} (in the unconjugated case) can in the Petrov--Galerkin formulations be simplified to
\begin{align}\label{Eq2:BilinearFormInserted}
	\begin{split}
	B_{\textsc{PGU}}(R_I\psi_n,R_J\phi_m) &= \lim_{\gamma\to\infty}\int_{\Omega_{\mathrm{a}}^\gamma} \left[\nabla(R_I\psi_n)\cdot \nabla (R_J\phi_m)- k^2 R_I\psi_n R_J\phi_m\right]\idiff\Omega\\
		&=\int_{\Omega_{\mathrm{a}}^+} \left[\nabla(R_I\psi_n)\cdot \nabla (R_J\phi_m)- k^2 R_I\psi_n R_J\phi_m\right]\idiff\Omega
	\end{split}
\end{align}
as the mentioned surface integral in the far field vanishes (this is however not the case for the Bubnov--Galerkin formulations). Recall that the radial shape functions are given by
\begin{align*}
	\phi_m(r) &= \euler^{\imag k (r-r_{\mathrm{a}})}Q_m\left(\frac{r_{\mathrm{a}}}{r}\right),\quad m = 1,\dots,N\\
	\psi_n(r) &= \euler^{\imag k (r-r_{\mathrm{a}})}\tilde{Q}_n\left(\frac{r_{\mathrm{a}}}{r}\right),\quad n = 1,\dots,N
\end{align*}
such that the derivative can be computed by
\begin{equation*}
	\deriv{\phi_m}{r} = \left[\imag kQ_m\left(\frac{r_{\mathrm{a}}}{r}\right) - \frac{r_{\mathrm{a}}}{r^2}Q_m'\left(\frac{r_{\mathrm{a}}}{r}\right)\right]\euler^{\imag k (r-r_{\mathrm{a}})}
\end{equation*}
and corresponding expression for $\psi_n$. Using the expression for the nabla operator found in \Cref{Eq2:GeneralNablaOperator}
\begin{align*}
	\nabla(R_I\psi_n)\cdot \nabla (R_J\phi_m) &= \frac{1}{h_{\mathrm{r}}^2}\pderiv{(R_I\psi_n)}{r}\pderiv{(R_J\phi_m)}{r} + \frac{1}{h_{\uptheta}^2}\pderiv{(R_I\psi_n)}{\vartheta}\pderiv{(R_J\phi_m)}{\vartheta} + \frac{1}{h_{\upvarphi}^2}\pderiv{(R_I\psi_n)}{\varphi}\pderiv{(R_J\phi_m)}{\varphi}\\
	 &= \frac{1}{h_{\mathrm{r}}^2}\pderiv{\psi_n}{r}\pderiv{\phi_m}{r}R_IR_J + \frac{1}{h_{\uptheta}^2}\psi_n\phi_m\pderiv{R_I}{\vartheta}\pderiv{R_J}{\vartheta}+ \frac{1}{h_{\upvarphi}^2}\psi_n\phi_m\pderiv{R_I}{\varphi}\pderiv{R_J}{\varphi}
\end{align*}
which multiplied with the Jacobian $J_1$ yields
\begin{align*}
	\nabla(R_I\psi_n)\cdot \nabla (R_J\phi_m) J_1&= \left[\left(r^2-\Upsilon^2\right)\pderiv{\psi_n}{r}\pderiv{\phi_m}{r}R_IR_J + \psi_n\phi_m\pderiv{R_I}{\vartheta}\pderiv{R_J}{\vartheta}\right. \\
	 &{\hskip8em\relax}\left.+ \frac{r^2-\Upsilon^2\cos^2\vartheta}{(r^2-\Upsilon^2)\sin^2\vartheta}\psi_n\phi_m\pderiv{R_I}{\varphi}\pderiv{R_J}{\varphi}\right]\sin\vartheta
\end{align*}
Combining all of this into \Cref{Eq2:BilinearFormInserted} yields
\begin{align}
	B_{\textsc{PGU}}(R_I\psi_n,R_J\phi_m) = &\int_0^{2\PI}\int_0^\PI K(\vartheta,\varphi)\sin\vartheta\idiff\vartheta\idiff\varphi
\end{align}
where
\begin{align*}
	K(\vartheta,\varphi) &= \int_{r_{\mathrm{a}}}^{\infty} \left\{\left(r^2-\Upsilon^2\right)\pderiv{\psi_n}{r}\pderiv{\phi_m}{r}R_IR_J + \psi_n\phi_m\pderiv{R_I}{\vartheta}\pderiv{R_J}{\vartheta}\right. \\
	 &{\hskip4em\relax}\left.+ \frac{r^2-\Upsilon^2\cos^2\vartheta}{(r^2-\Upsilon^2)\sin^2\vartheta}\psi_n\phi_m\pderiv{R_I}{\varphi}\pderiv{R_J}{\varphi}-k^2(r^2-\Upsilon^2\cos^2\vartheta)\psi_n\phi_m R_I R_J\right\}\idiff r.
\end{align*}
Inserting the expressions for the radial shape functions $\phi$ and $\psi$ (with Einstein's summation convention) with their corresponding derivatives one obtains the following expression using the substitution $\rho = \frac{r}{r_{\mathrm{a}}}$ and the notation $\varrho_1 = \Upsilon/r_{\mathrm{a}}$ (the eccentricity of the infinite-element spheroid), $\varrho_2=kr_{\mathrm{a}}$ and $\varrho_3=k\Upsilon$
\begin{align*} 
	 K(\vartheta,\varphi) &= \left\{R_IR_J\left[-2\varrho_2^2B_{\tilde{n}+\tilde{m}}^{(1)} - \imag \varrho_2(\tilde{n}+\tilde{m}+2)B_{\tilde{n}+\tilde{m}+1}^{(1)} + \left[\tilde{m}(\tilde{n}+2) +\varrho_3^2\right]B_{\tilde{n}+\tilde{m}+2}^{(1)} \phantom{\left(\pderiv{r_{\mathrm{a}}}{\varphi}\right)^2} \right.\right.\\
	 &\left.\left.{\hskip5em\relax}+ \imag \varrho_1^2\varrho_2(\tilde{n}+\tilde{m}+2)B_{\tilde{n}+\tilde{m}+3}^{(1)} - \tilde{m}(\tilde{n}+2)\varrho_1^2 B_{\tilde{n}+\tilde{m}+4}^{(1)}+\varrho_3^2\cos^2\vartheta B_{\tilde{n}+\tilde{m}+2}^{(1)}\right] \right.\\
	 &\left.{\hskip2em\relax}+ \pderiv{R_I}{\vartheta}\pderiv{R_J}{\vartheta}B_{\tilde{n}+\tilde{m}+2}^{(1)}+\pderiv{R_I}{\varphi}\pderiv{R_J}{\varphi}\frac{1}{\sin^2\vartheta}\left(B_{\tilde{n}+\tilde{m}+1}^{(2)} - \varrho_1^2\cos^2\vartheta B_{\tilde{n}+\tilde{m}+3}^{(2)}\right)\right\}r_{\mathrm{a}}\euler^{-2\imag \varrho_2}\tilde{D}_{n\tilde{n}}D_{m\tilde{m}}
\end{align*}
where the radial integrals
\begin{equation*}
	B_{n}^{(1)} = \int_{1}^\infty \frac{\euler^{2\imag \varrho_2 \rho}}{\rho^n}\idiff \rho \qquad B_{n}^{(2)} = \int_{1}^\infty \frac{\euler^{2\imag \varrho_2\rho}}{(\rho^2-\varrho_1^2)\rho^{n-1}}\idiff \rho,\qquad n\geq 1
\end{equation*}
can be evaluated according to formulas in \Cref{Sec2:radIntegrals}.

Assume that the artificial boundary $\Gamma_{\mathrm{a}}$ is parameterized by $\xi$ and $\eta$. As $\Gamma_{\mathrm{a}}$ is a surface with constant radius, $r=r_{\mathrm{a}}$, in the prolate spheroidal coordinate system, it may also be parameterized by $\vartheta$ and $\varphi$. Therefore,
\begin{equation}
	\diff \vartheta\diff\varphi = \begin{vmatrix}
		\pderiv{\vartheta}{\xi} & \pderiv{\varphi}{\xi}\\
		\pderiv{\vartheta}{\eta} & \pderiv{\varphi}{\eta}		
	\end{vmatrix}\diff\xi\diff\eta
\end{equation}
where
\begin{align*}
	\pderiv{\vartheta}{\xi} &= \pderiv{\vartheta}{x}\pderiv{x}{\xi} + \pderiv{\vartheta}{y}\pderiv{y}{\xi} + \pderiv{\vartheta}{z}\pderiv{z}{\xi},\qquad
	\pderiv{\vartheta}{\eta} = \pderiv{\vartheta}{x}\pderiv{x}{\eta} + \pderiv{\vartheta}{y}\pderiv{y}{\eta} + \pderiv{\vartheta}{z}\pderiv{z}{\eta}\\
	\pderiv{\varphi}{\xi} &= \pderiv{\varphi}{x}\pderiv{x}{\xi} + \pderiv{\varphi}{y}\pderiv{y}{\xi} + \pderiv{\varphi}{z}\pderiv{z}{\xi},\qquad
	\pderiv{\varphi}{\eta} = \pderiv{\varphi}{x}\pderiv{x}{\eta} + \pderiv{\varphi}{y}\pderiv{y}{\eta} + \pderiv{\varphi}{z}\pderiv{z}{\eta}
\end{align*}
and the inverse partial derivatives with respect to the coordinate transformation (from the prolate spheroidal coordinate system to the Cartesian coordinate system) is found in \Cref{Eq2:dProlateSphericalCoordinatesdX}. This Jacobian matrix may be evaluated by
\begin{equation}
	J_3 = \begin{bmatrix}
		\pderiv{\vartheta}{\xi} & \pderiv{\vartheta}{\eta}\\
		\pderiv{\varphi}{\xi}	 & \pderiv{\varphi}{\eta}
	\end{bmatrix} = \begin{bmatrix}
		\pderiv{\vartheta}{x} & \pderiv{\vartheta}{y} & \pderiv{\vartheta}{z}\\
		\pderiv{\varphi}{x} & \pderiv{\varphi}{y} & \pderiv{\varphi}{z}
	\end{bmatrix}\begin{bmatrix}
		\pderiv{x}{\xi} & \pderiv{x}{\eta}\\
		\pderiv{y}{\xi} & \pderiv{y}{\eta}\\
		\pderiv{z}{\xi} & \pderiv{z}{\eta}
	\end{bmatrix}
\end{equation}
and the derivatives of the basis functions may then be computed by
\begin{equation}
	\begin{bmatrix}
		\pderiv{R_I}{\vartheta}\\
		\pderiv{R_I}{\varphi}
	\end{bmatrix} = J_3^{-\transpose}\begin{bmatrix}
		\pderiv{R_I}{\xi}\\
		\pderiv{R_I}{\eta}	
	\end{bmatrix}.
\end{equation}
Defining the angular integrals
\begin{equation}\label{Eq2:infiniteElementsSurfaceIntegrals}
\begin{alignedat}{2}
	& A_{IJ}^{(1)} = \int_0^{2\PI}\int_0^\PI R_I R_J\sin\vartheta\idiff\vartheta\idiff\varphi,\qquad\qquad && A_{IJ}^{(2)} = \int_0^{2\PI}\int_0^\PI \pderiv{R_I}{\vartheta} \pderiv{R_J}{\vartheta}\sin\vartheta\idiff\vartheta\idiff\varphi\\
	& A_{IJ}^{(3)} = \int_0^{2\PI}\int_0^\PI R_I R_J\cos^2\vartheta\sin\vartheta\idiff\vartheta\idiff\varphi, \quad && A_{IJ}^{(4)} = \int_0^{2\PI}\int_0^\PI \pderiv{R_I}{\varphi} \pderiv{R_J}{\varphi}\frac{1}{\sin\vartheta}\idiff\vartheta\idiff\varphi\\
	& A_{IJ}^{(5)} = \int_0^{2\PI}\int_0^\PI \pderiv{R_I}{\varphi} \pderiv{R_J}{\varphi}\frac{\cos^2\vartheta}{\sin\vartheta}\idiff\vartheta\idiff\varphi
\end{alignedat}	
\end{equation}
the bilinear form may then finally be written as (Einstein's summation convention is used for the indices $\tilde{n}$ and $\tilde{m}$)
\begin{align}\label{Eq2:finalBilinearFormB_uc_a}
\begin{split}
	B_{\textsc{PGU}}(R_I\psi_n,R_J\phi_m) 
	 &= \left\{A_{IJ}^{(1)}\left[-2\varrho_2^2 B_{\tilde{n}+\tilde{m}}^{(1)} - \imag\varrho_2(\tilde{n}+\tilde{m}+2)B_{\tilde{n}+\tilde{m}+1}^{(1)} + \left[(\tilde{n}+2)\tilde{m}+\varrho_3^2\right]B_{\tilde{n}+\tilde{m}+2}^{(1)}\right. \right. \\
	  &{\hskip4em\relax}\left. \left. + \imag\varrho_1\varrho_3(\tilde{n}+\tilde{m}+2)B_{\tilde{n}+\tilde{m}+3}^{(1)} - \varrho_1^2(\tilde{n}+2)\tilde{m}B_{\tilde{n}+\tilde{m}+4}^{(1)}\right] \right. \\
	 &{\hskip2em\relax}\left. +A_{IJ}^{(2)}B_{\tilde{n}+\tilde{m}+2}^{(1)} + \varrho_3^2A_{IJ}^{(3)} B_{\tilde{n}+\tilde{m}+2}^{(1)}\right.\\
	 &{\hskip2em\relax}\left. +A_{IJ}^{(4)}B_{\tilde{n}+\tilde{m}+1}^{(2)} -\varrho_1^2A_{IJ}^{(5)} B_{\tilde{n}+\tilde{m}+3}^{(2)}
	 \right\}r_{\mathrm{a}}\euler^{-2\imag \varrho_2}D_{m\tilde{m}}\tilde{D}_{n\tilde{n}}.
\end{split}
\end{align}
For completeness, the formulas for the other three formulations are included
\begin{equation}
\begin{aligned}
	B_{\textsc{BGU}}(R_I\psi_n,R_J\phi_m) 
	 &= \left\{A_{IJ}^{(1)}\left[-2\varrho_2^2 B_{\tilde{n}+\tilde{m}-2}^{(1)}(1-\delta_{\tilde{n}1}\delta_{\tilde{m}1}) - \imag\varrho_2(\tilde{n}+\tilde{m})B_{\tilde{n}+\tilde{m}-1}^{(1)} + \left(\tilde{n}\tilde{m}+\varrho_3^2\right)B_{\tilde{n}+\tilde{m}}^{(1)}\right. \right. \\
	  &{\hskip4em\relax}\left. \left. + \imag\varrho_1\varrho_3(\tilde{n}+\tilde{m})B_{\tilde{n}+\tilde{m}+1}^{(1)} - \varrho_1^2\tilde{n}\tilde{m}B_{\tilde{n}+\tilde{m}+2}^{(1)}\right] \right. \\
	 &{\hskip2em\relax}\left. +A_{IJ}^{(2)}B_{\tilde{n}+\tilde{m}}^{(1)} + \varrho_3^2A_{IJ}^{(3)} B_{\tilde{n}+\tilde{m}}^{(1)}\right.\\
	 &{\hskip2em\relax}\left. +A_{IJ}^{(4)}B_{\tilde{n}+\tilde{m}-1}^{(2)} -\varrho_1^2A_{IJ}^{(5)} B_{\tilde{n}+\tilde{m}+1}^{(2)}
	 \right\}r_{\mathrm{a}}\euler^{-2\imag \varrho_2}D_{m\tilde{m}}\tilde{D}_{n\tilde{n}}\\
	 &\quad-\imag\varrho_2 r_{\mathrm{a}} D_{m1}\tilde{D}_{n1}A_{IJ}^{(1)}
\end{aligned}
\end{equation}
\begin{equation}
\begin{aligned}
	B_{\textsc{PGC}}(R_I\psi_n,R_J\phi_m) 
	 &= \left\{A_{IJ}^{(1)}\left[-\imag\varrho_2(\tilde{n}-\tilde{m}+2)B_{\tilde{n}+\tilde{m}+1}^{(1)} + \left[(\tilde{n}+2)\tilde{m}-\varrho_3^2\right]B_{\tilde{n}+\tilde{m}+2}^{(1)}\right. \right. \\
	  &{\hskip4em\relax}\left. \left. + \imag\varrho_1\varrho_3(\tilde{n}-\tilde{m}+2)B_{\tilde{n}+\tilde{m}+3}^{(1)} - \varrho_1^2(\tilde{n}+2)\tilde{m}B_{\tilde{n}+\tilde{m}+4}^{(1)}\right] \right. \\
	 &{\hskip2em\relax}\left. +A_{IJ}^{(2)}B_{\tilde{n}+\tilde{m}+2}^{(1)} + \varrho_3^2A_{IJ}^{(3)} B_{\tilde{n}+\tilde{m}+2}^{(1)}\right.\\
	 &{\hskip2em\relax}\left. +A_{IJ}^{(4)}B_{\tilde{n}+\tilde{m}+1}^{(2)} -\varrho_1^2A_{IJ}^{(5)} B_{\tilde{n}+\tilde{m}+3}^{(2)}
	 \right\}r_{\mathrm{a}}D_{m\tilde{m}}\tilde{D}_{n\tilde{n}}
\end{aligned}
\end{equation}
\begin{equation}
\begin{aligned}
	B_{\textsc{BGC}}(R_I\psi_n,R_J\phi_m) 
	 &= \left\{A_{IJ}^{(1)}\left[-\imag\varrho_2(\tilde{n}-\tilde{m})B_{\tilde{n}+\tilde{m}-1}^{(1)} + \left(\tilde{n}\tilde{m}-\varrho_3^2\right)B_{\tilde{n}+\tilde{m}}^{(1)}\right. \right. \\
	  &{\hskip4em\relax}\left. \left. + \imag\varrho_1\varrho_3(\tilde{n}-\tilde{m})B_{\tilde{n}+\tilde{m}+1}^{(1)} - \varrho_1^2\tilde{n}\tilde{m}B_{\tilde{n}+\tilde{m}+2}^{(1)}\right] \right. \\
	 &{\hskip2em\relax}\left. +A_{IJ}^{(2)}B_{\tilde{n}+\tilde{m}}^{(1)} + \varrho_3^2A_{IJ}^{(3)} B_{\tilde{n}+\tilde{m}}^{(1)}\right.\\
	 &{\hskip2em\relax}\left. +A_{IJ}^{(4)}B_{\tilde{n}+\tilde{m}-1}^{(2)} -\varrho_1^2A_{IJ}^{(5)} B_{\tilde{n}+\tilde{m}+1}^{(2)}
	 \right\}r_{\mathrm{a}}D_{m\tilde{m}}\tilde{D}_{n\tilde{n}}\\
	 &\quad-\imag r_{\mathrm{a}}\varrho_2 D_{m1}\tilde{D}_{n1}A_{IJ}^{(1)} 
\end{aligned}
\end{equation}
where $\delta_{ij}$ is the Kronecker delta function in \Cref{Eq2:Kronecker}.

%% file: contents/radialIntegrals.tex
\section{Evaluation of radial integrals}
\label{Sec2:radIntegrals}
The exponential integral
\begin{equation}
	E_n(z) = \int_1^\infty \frac{\euler^{-z\rho}}{\rho^n}\idiff \rho,\qquad \Re(z) \geq 0
\end{equation}
is of great importance for the unconjugated formulations in the IEM. It is therefore important to be able to evaluate the integral accurately and efficiently, also for large (absolute) values of $z$ (which will correspond to high frequencies). In~\cite[p. 229, 5.1.12]{Abramowitz1965hom} the series representation for evaluation of these functions can be found\footnote{Here, $\upgamma$ is the Euler-Mascheroni constant which is defined by 
\begin{equation*}
\upgamma  = \lim_{n\to\infty}\left[-\ln(n)+\sum_{m=1}^n \frac{1}{m}\right]=0.577215664901532860606512090082\dots.
\end{equation*}}
\begin{equation}\label{Eq2:seriesRepresentation}
	E_n(z) = \frac{(-z)^{n-1}}{(n-1)!}\left[-\ln z-\upgamma+\sum_{m=1}^{n-1}\frac{1}{m}\right] -\sum_{\substack{m=0\\m\neq n-1}}^\infty \frac{(-z)^m}{(m-n+1)m!}
\end{equation}
with the empty sum interpreted to be zero. Moreover, using the continued fraction notation
\begin{equation}
	b_0+\frac{a_1}{b_1+}\frac{a_2}{b_2+}\frac{a_3}{b_3+}\cdots = b_0+\cfrac{a_1}{b_1+\cfrac{a_2}{b_2+\cfrac{a_3}{b_3+\cdots}}}
\end{equation}
the continued fraction representation of these functions are given by~\cite[p. 229, 5.1.22]{Abramowitz1965hom}
\begin{equation}\label{Eq2:contFracRepresentation}
	E_n(z) = \euler^{-z}\left(\frac{1}{z+}\frac{n}{1+}\frac{1}{z+}\frac{n+1}{1+}\frac{2}{z+}\frac{n+2}{1+}\frac{3}{z+}\cdots\right).
\end{equation}
In~\cite[p. 222]{Press1988nri} Press et al. present an even faster converging continued fraction given by
\begin{equation}\label{Eq2:contFracRepresentation2}
	E_n(z) = \euler^{-z}\left(\frac{1}{z+n-}\frac{1\cdot n}{z+n+2-}\frac{2(n+1)}{z+n+4-}\frac{3(n+2)}{z+n+6-}\cdots\right).
\end{equation}
It is here suggested to use \Cref{Eq2:seriesRepresentation} when $|z|\lesssim 1$ and \Cref{Eq2:contFracRepresentation} or \Cref{Eq2:contFracRepresentation2} when $|z|\gtrsim 1$. Press et al. then continue to present efficient algorithms for evaluation of these formulas.

Using series expansions at infinity
\begin{equation}
	\frac{1}{\rho^2-\varrho_1^2} = \frac{1}{\varrho_1^2}\sum_{j=1}^\infty \left(\frac{\varrho_1}{\rho}\right)^{2j},
\end{equation}
the radial integrals for 3D infinite elements may be computed by
\begin{align}
	\int_1^\infty \frac{1}{\rho^n}\idiff\rho &= \frac{1}{n-1}\\
	\int_1^\infty \frac{1}{(\rho^2-\varrho_1^2)\rho^{n-1}}\idiff\rho &= \sum_{j=0}^\infty \frac{\varrho_1^{2j}}{2j+n}
\end{align}
in the conjugated case and
\begin{align}
	\int_1^\infty \frac{\euler^{2\imag \varrho_2 \rho}}{\rho^n}\idiff\rho &= E_n(-2\imag \varrho_2)\\
	\int_1^\infty \frac{\euler^{2\imag \varrho_2\rho}}{(\rho^2-\varrho_1^2)\rho^{n-1}}\idiff\rho &= \sum_{j=0}^\infty \varrho_1^{2j}E_{2j+n+1}(-2\imag \varrho_2)
\end{align}
in the unconjugated case.

%% file: contents/betssi_description.tex
\section{The stripped BeTSSi submarine model}
\label{Sec2:BeTSSi_description}
\begin{figure}
	\centering
	\includegraphics[width=\textwidth]{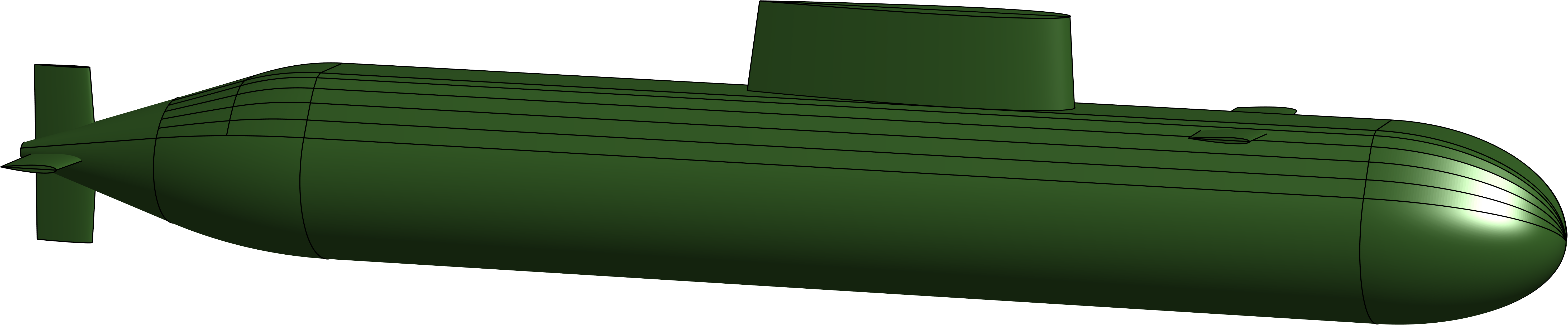}
	\caption{Outer pressure hull for BeTSSi submarine.}
	\label{Fig2:BeTSSi_BC}
\end{figure}

\begin{figure}
	\centering
	\includegraphics[width=\textwidth]{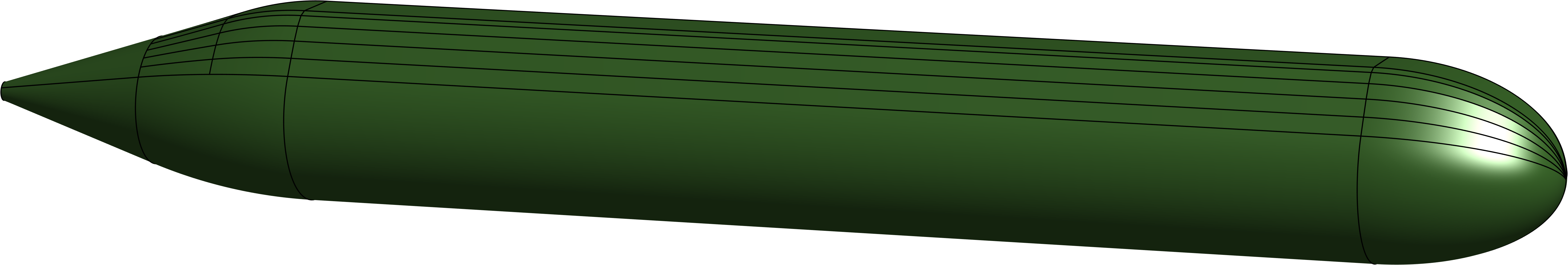}
	\caption{The stripped BeTSSi submarine model.}
	\label{Fig2:BeTSSi_BC_stripped}
\end{figure}

\begin{figure}
	\centering
	\includegraphics[scale=1]{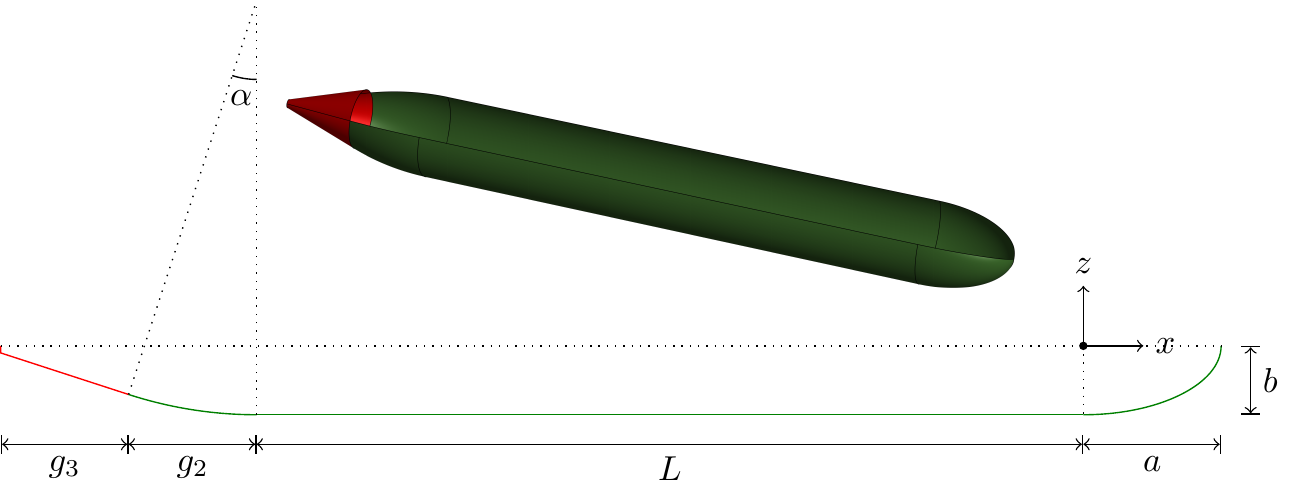}
	\caption{The side line of the lower part of the BeTSSi submarine. The side lines are formed (from the right) by an ellipse with semi-major axis $a$ and semi-minor axis $b$, followed by a straight line of length $L$, then an arc of angle $\alpha$ and finally two straight lines. The latter two straight lines (in red) are rotated about the $x$-axis and the remaining part (in green) are rotated an angle $\beta$ around the $x$-axis.}
	\label{Fig2:bettsi_bottom}
\end{figure}

\begin{figure}
	\centering
	\includegraphics[scale=1]{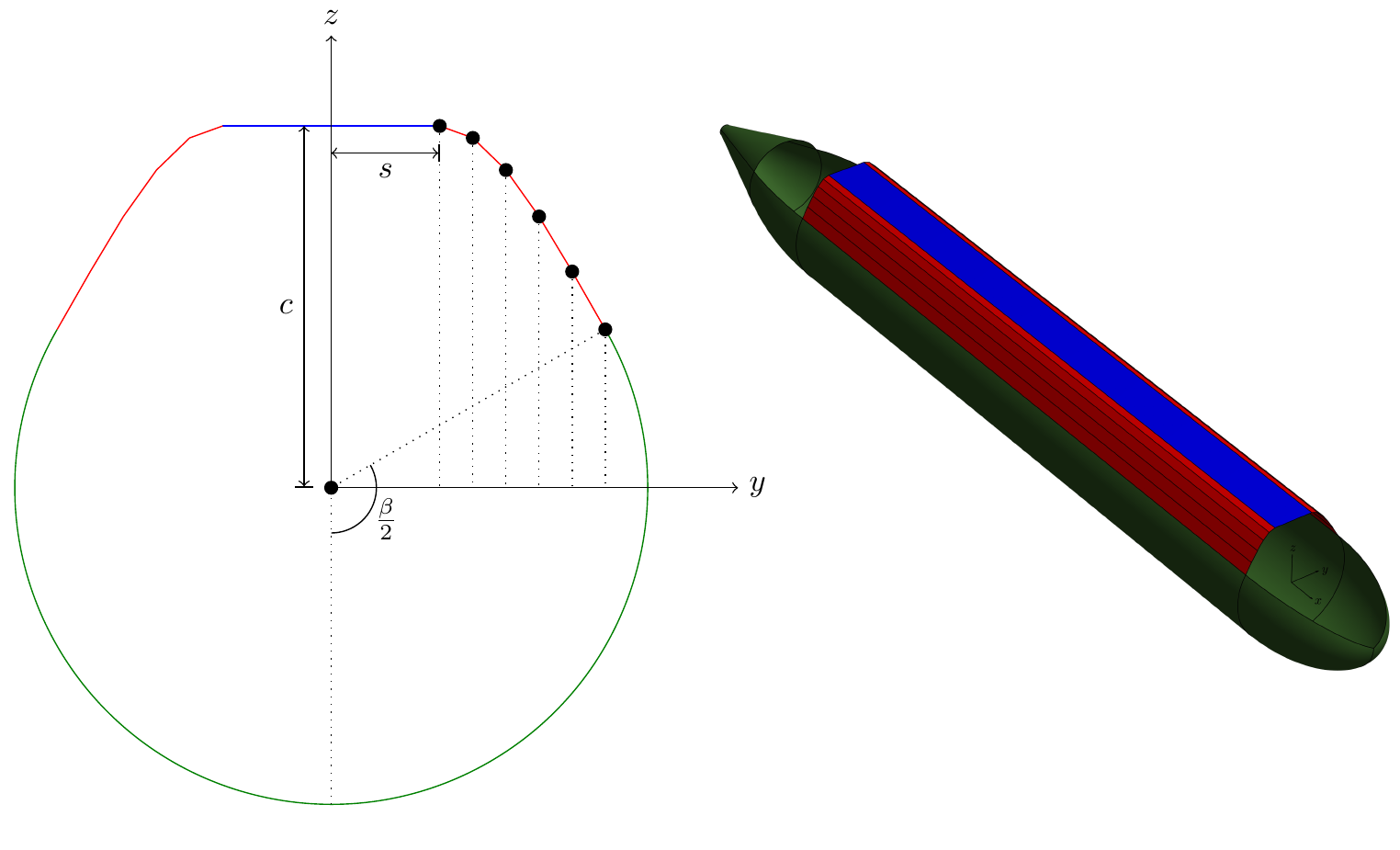}
	\caption{The transition (red line) from the axisymmetric hull (green line) to the deck (blue line) is given by sampling a cubic polynomial, $P(y)$, at 6 equidistant points in the $y$-direction and connecting the resulting points with straight lines (corresponding 6 points are found for negative values $y$-values, $(0,y,P(|y|))$).}
	\label{Fig2:bettsi_top}
\end{figure}

\begin{figure}
	\centering    
	\begin{subfigure}{0.49\textwidth}
		\centering
		\includegraphics[width=\textwidth]{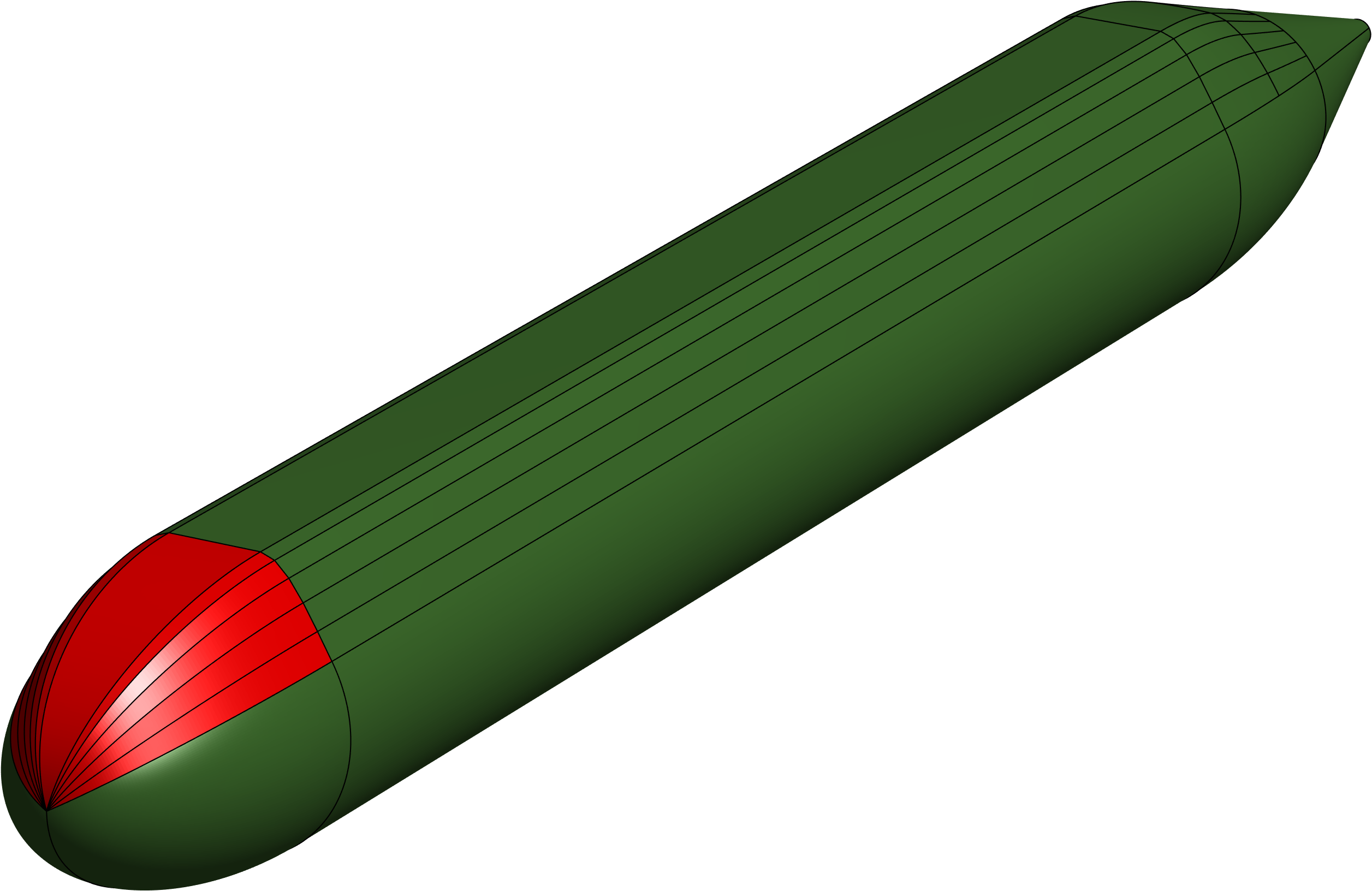}
		\caption{Illustration of upper bow part.}
		\label{Fig2:bettsi_upperBow}
	\end{subfigure}
	~    
	\begin{subfigure}{0.49\textwidth}
		\centering
		\includegraphics[width=\textwidth]{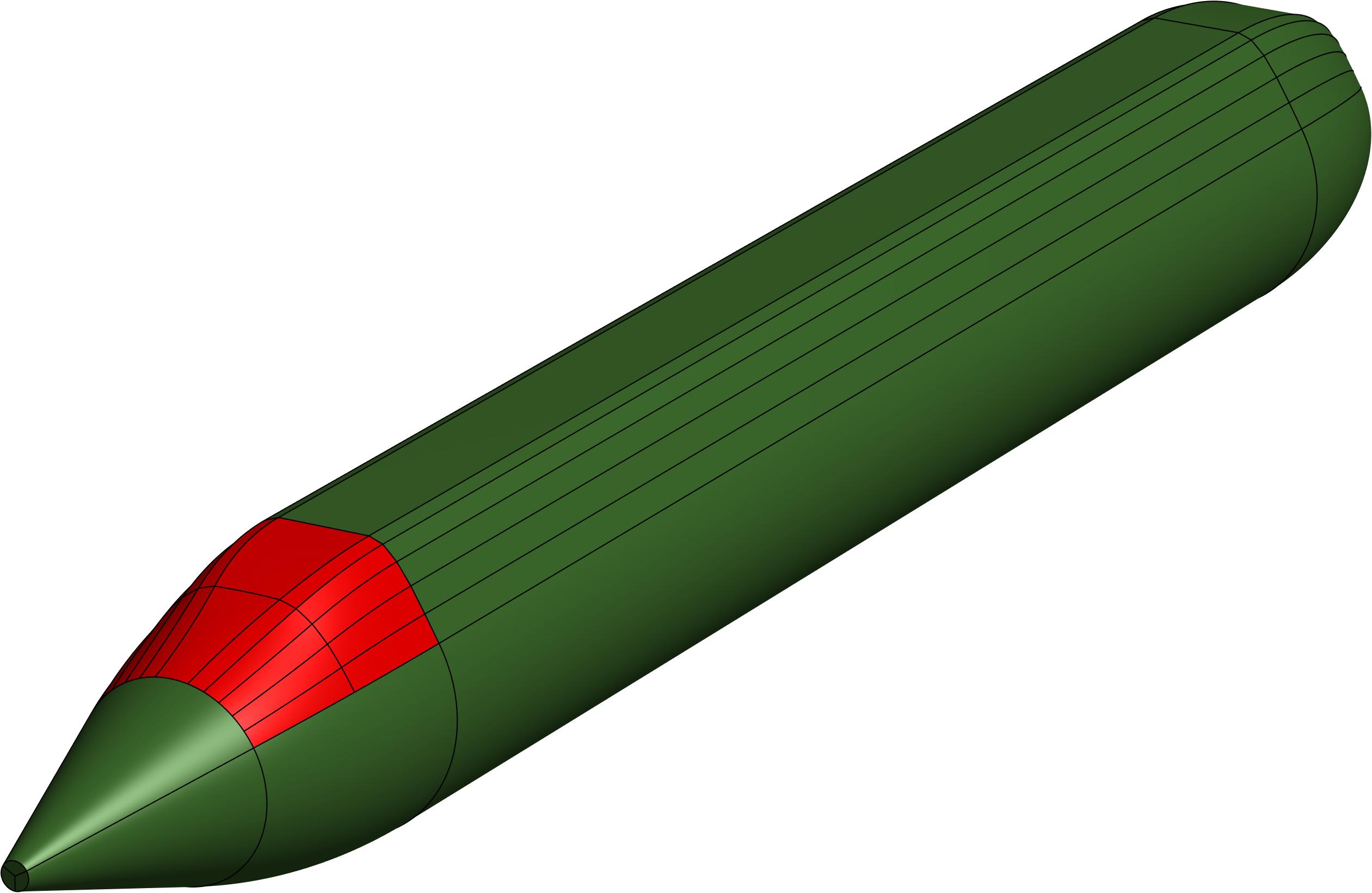}
		\caption{Illustration of upper transition part.}
		\label{Fig2:bettsi_upperPartOfTailSection}
	\end{subfigure}
	\caption{Final patches for the stripped BeTSSi submarine.}
\end{figure}

\begin{figure}
	\centering    
	\begin{subfigure}{0.49\textwidth}
		\centering
		\includegraphics[width=0.9\textwidth]{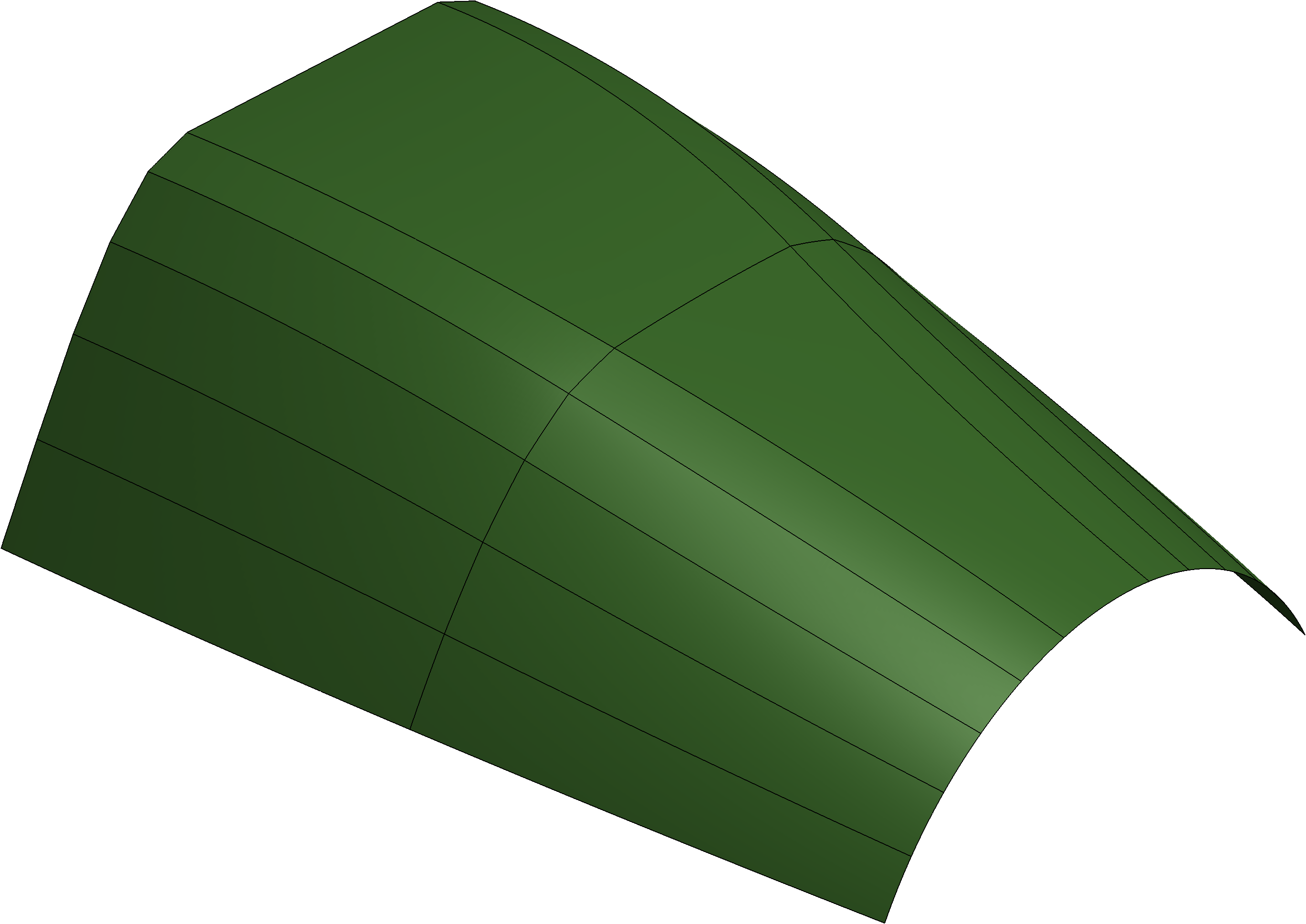}
		\caption{Illustration of the mesh.}
		\label{Fig2:BeTSSi_BC_tailSection}
	\end{subfigure}
	~    
	\begin{subfigure}{0.49\textwidth}
		\centering
		\includegraphics{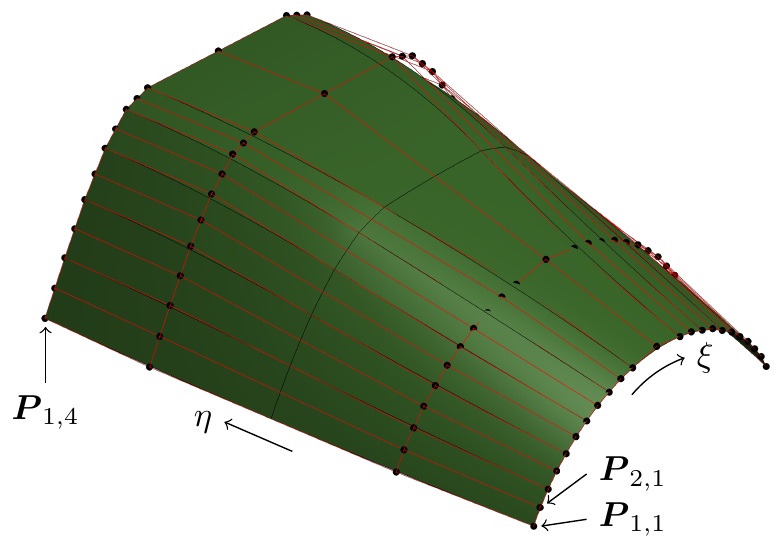}
		\caption{Illustration of the control polygon mesh.}
		\label{Fig2:BeTSSi_BC_tailSection_cp}
	\end{subfigure}
	\caption{Illustration of the upper transition part of the tail.}
\end{figure}

\begin{figure}
	\centering    
	\begin{subfigure}{0.49\textwidth}
		\centering
		\includegraphics{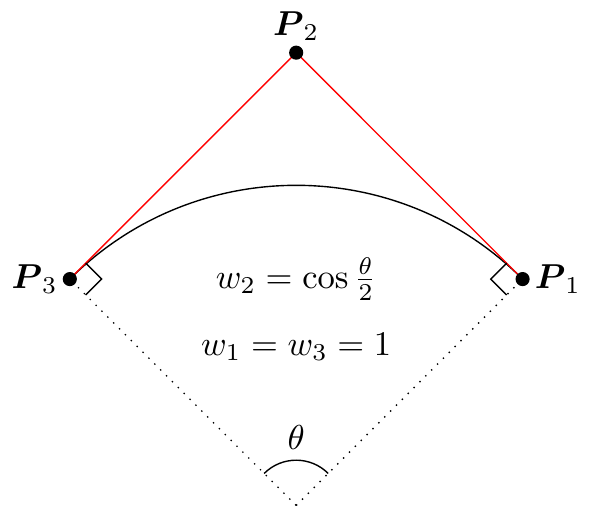}
		\caption{NURBS parametrization of arc of angle $\theta$ using three control points $\{\vec{P}_i\}_{i=1}^3$, the weights $\{w_i\}_{i=1}^3$ and the open knot vector $\Xi=\{0,0,0,1,1,1\}$.}
		\label{Fig2:arcParam1}
	\end{subfigure}
	~    
	\begin{subfigure}{0.49\textwidth}
		\centering
		\includegraphics{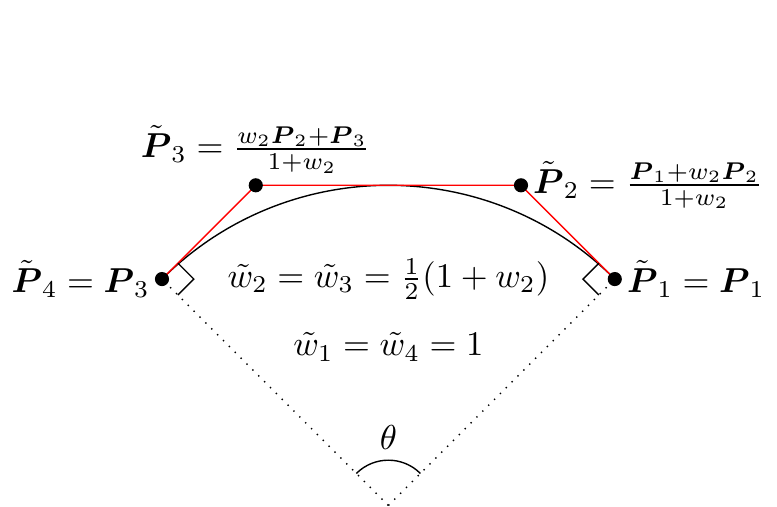}
		\caption{NURBS parametrization of arc of angle $\theta$ using four control points $\{\tilde{\vec{P}}_i\}_{i=1}^4$, the weights $\{\tilde{w}_i\}_{i=1}^4$ and the open knot vector $\tilde{\vec{t}}_\upxi=\{0,0,0,0.5,1,1,1\}$.}
		\label{Fig2:arcParam2}
	\end{subfigure}
	\caption{Two ways of parametrizing an arc using NURBS~\cite[p. 315]{Piegl1997tnb}.}
\end{figure}

In this section a simplified version of the BeTSSi submarine model (depicted in \Cref{Fig2:BeTSSi_BC}) will be presented. Namely a \textit{stripped BeTSSi submarine model} without sail and rudders as in \Cref{Fig2:BeTSSi_BC_stripped}.
The relevant BeTSSi parameters for the work presented herein are given in \Cref{Tab2:BeTSSiParameters}.
\begin{table}
	\centering
	\caption{\textbf{BeTSSi submarine:} Parameters for the BeTSSi submarine benchmark.}
	\label{Tab2:BeTSSiParameters}
	\begin{tabular}{l l}
		\toprule
		Parameter & Description\\
		\midrule
		$P_{\mathrm{inc}}=\SI{1}{Pa}$ & Amplitude of incident wave\\
		$E = \SI{2.10e11}{Pa}$ & Young's modulus\\
		$\nu = 0.3$ & Poisson's ratio\\
		$\rho_{\mathrm{s}} = \SI{7850}{kg.m^{-3}}$ & Density of solid\\
		$\rho_{\mathrm{f}} = \SI{1000}{kg.m^{-3}}$ & Density of water\\
		$c_{\mathrm{f}} = \SI{1500}{m.s^{-1}}$ & Speed of sound in water\\
		$t=\SI{0.01}{m}$ & Thickness of pressure hull\\
		$\alpha=\ang{18}$ & Arc angle of transition to the tail cone\\
		$\beta=\ang{240}$ & Rotational angle for the axisymmetric lower part of the pressure hull\\
		$g_2=\SI{6.5}{m}$ & Distance in $x$-direction of transition to the tail cone\\
		$g_3=\SI{6.5}{m}$ & Distance in $x$-direction of the tail cone\\
		$L=\SI{42}{m}$ & Length of the deck\\
		$a=\SI{7}{m}$ & Semi-major axis of bow\\
		$b=\SI{3.5}{m}$ & Semi-major axis of bow\\
		$c=\SI{4}{m}$ & Height from $x$-axis to the deck\\
		$s=\SI{1.2}{m}$ & Half of the width of the deck\\
		\bottomrule
	\end{tabular}
\end{table}
The model is symmetric about the $xz$-plane and has rotational symmetry for the lower part as described in \Cref{Fig2:bettsi_bottom}.
The transition from this axisymmetric part to the deck is described in \Cref{Fig2:bettsi_top}. This transition as well as the deck itself, contains a set of rectangular panels of length $L$.
The polynomial $P(y)$, is uniquely defined by the requirement that it defines a smooth transition between the hull and the deck. More precisely, the following requirement must be satisfied: 
\begin{alignat*}{3}
	P(s) &= c,\quad  &&P\left(b\sin\frac{\beta}{2}\right) = -b\cos\frac{\beta}{2}\\
	P'(s) &= 0,\quad &&P'\left(b\sin\frac{\beta}{2}\right) = \tan\frac{\beta}{2}
\end{alignat*}
which gives the polynomial
\begin{equation}
	P(y) = c+C_1(y-s)^2+C_2(y-s)^3
\end{equation}
where
\begin{equation*}
	C_1 = -\frac{3C_4+C_3\tan\frac{\beta}{2}}{C_3^2}, \quad
	C_2 = \frac{2C_4+C_3\tan\frac{\beta}{2}}{C_3^3}, \quad
	C_3 = b\sin\frac{\beta}{2}-s, \quad
	C_4 = c+b\cos\frac{\beta}{2}.
\end{equation*}
The upper part of the bow (highlighted in \Cref{Fig2:bettsi_upperBow}) is obtained by linear lofting of elliptic curves from the 12 points described in \Cref{Fig2:bettsi_top} to the tip of the bow. The upper part of the tail section (highlighted in \Cref{Fig2:bettsi_upperPartOfTailSection}) is connected using a tensor NURBS surface of degree 2 such that it defines a smooth transition from the axisymmetric cone to the deck. More precisely, the upper part of the cone tail is divided into 12 arcs with angle $\frac{2\PI-\beta}{12}$, and the resulting points are connected to corresponding points on the transition to the deck from the axisymmetric hull. 
As illustrated in \Cref{Fig2:BeTSSi_BC_tailSection}, the NURBS patch is given by 22 elements. Thus, $4\cdot 23 = 92$ control points, $\vec{P}_{i,j}$, is needed as shown in \Cref{Fig2:BeTSSi_BC_tailSection_cp} (23 and 4 control points in the $\xi$ direction and $\eta$ direction, respectively).
The control points $\vec{P}_{1,j}$ and $\vec{P}_{23,j}$ for $j=1,2,3,4$ must be defined as in \Cref{Fig2:arcParam2}, while the control points $\vec{P}_{i,1}$ must be defined as in \Cref{Fig2:arcParam1} (with corresponding weights). For $2\leq i\leq 22$ the weights are defined by $w_{i,j}=w_{i,1}$ for $j=2,3,4$. That is,
\begin{equation}
	w_{i,j} = \begin{cases}
		1 & i\,\, \text{odd}\\
		\cos\left(\frac{2\PI-\beta}{24}\right) & i\,\, \text{even},\,i\neq 12\\
		\cos\left(\frac{2\PI-\beta}{12}\right) &  i = 12
		\end{cases}		
\end{equation}
The location of the control points $\vec{P}_{i,j}$, $j=2,3$ and $2\leq i\leq 22$, are determined by the requirement that the $x$ component is the same as $\vec{P}_{1,j}$ and the fact that the control polygon lines must be tangential to the surface both at the deck and the cone tail.

The inner surface of the BeTSSi submarine is generated by scaling a copy of the outer surface with the following change in the parameters $a\to a-t$, $b\to b-t$, $c\to c-t$, $s\to s-t/2$,  $g_2\to g_2-t/2$ and $g_3\to g_3-t/2$ ($\alpha$, $\beta$ and $l$ remain unchanged).

%% file: contents/transformNURBStoBspline.tex
\section{Approximating NURBS parametrizations with B-spline parametrizations}
\label{Sec2:NURBStransformation}
Starting with any NURBS parametrizations of a geometry where every internal knot has multiplicity $m=\check{p}_\upxi$ in the $\xi$-direction and correspondingly in the other two parameter directions, we want to transform the NURBS parametrization of the exact geometry, to a B-spline representation. This representation approximates the geometry by interpolating the geometry at $n_\upxi\cdot n_\upeta\cdot n_\upzeta$ (not necessarily unique) physical points resulting from a grid in the parametric space.

Let $\vec{X}$ be the NURBS parametrization of the geometry (with notation similar to \cite[p. 51]{Cottrell2009iat})
\begin{equation}
	\vec{X}(\xi,\eta,\zeta) = \sum_{i=1}^{n_\upxi}\sum_{j=1}^{n_\upeta}\sum_{l=1}^{n_\upzeta} R_{i,j,l}^{\check{p}_\upxi,\check{p}_\upeta,\check{p}_\upzeta}(\xi,\eta,\zeta)\vec{P}_{i,j,l},
\end{equation}
with knot vectors $\Xi$, $\Eta$ and $\Zeta$, polynomial order $\check{p}_\upxi$, $\check{p}_\upeta$ and $\check{p}_\upzeta$. For each control point $\vec{P}_{i,j,l}$ we will need a corresponding interpolating point $\vec{Q}_{i,j,l}$ which will be located at the grid point $\left(\tilde{\xi}_i,\tilde{\eta}_j,\tilde{\zeta}_l\right)$. These points in the parameter domain are chosen to be the Greville abscissae
\begin{align}
	\tilde{\xi}_i &= \frac{1}{\check{p}_\upxi}\sum_{\tilde{i} = i+1}^{i+\check{p}_\upxi} \xi_{\tilde{i}},\quad i = 1,\dots, n_\upxi\\
	\tilde{\eta}_j &= \frac{1}{\check{p}_\upeta}\sum_{\tilde{j} = j+1}^{j+\check{p}_\upeta} \eta_{\tilde{j}},\quad j = 1,\dots, n_\upeta\\
	\tilde{\zeta}_l &= \frac{1}{\check{p}_\upzeta}\sum_{\tilde{l} = l+1}^{l+\check{p}_\upzeta} \zeta_{\tilde{l}},\quad l = 1,\dots, n_\upzeta,
\end{align}
where $\xi_i$, $\eta_j$ and $\zeta_l$ are the knots of the knot vectors $\Xi$, $\Eta$ and $\Zeta$, respectively.

We can now compute the interpolation points $\vec{Q}_{i,j,l}$ by
\begin{equation}
	\vec{Q}_{i,j,l} = \vec{X}(\tilde{\xi}_i,\tilde{\eta}_j,\tilde{\zeta}_l).
\end{equation}
To find a B-spline approximation of the geometry which interpolates the points $\vec{Q}_{i,j,l}$, we want this new parametrization $\tilde{\vec{X}}$ to be based on $\vec{X}$ such that their order and knot vectors are equal. As all weights will be set to 1 (to get a B-spline parametrization), we are only left with dofs in the control points, $\tilde{\vec{P}}_{i,j,l}$, of the B-spline parametrization. To find these points we require
\begin{equation}
	\tilde{\vec{X}}(\tilde{\xi}_i,\tilde{\eta}_j,\tilde{\zeta}_l) = \sum_{\tilde{i}=1}^{n_\upxi}\sum_{\tilde{j}=1}^{n_\upeta}\sum_{\tilde{l}=1}^{n_\upzeta} B_{\tilde{i},\check{p}_\upxi,\Xi}(\tilde{\xi}_i)B_{\tilde{j},\check{p}_\upeta,\Eta}(\tilde{\eta}_j)B_{\tilde{l},\check{p}_\upzeta,\Zeta}(\tilde{\zeta}_l)\tilde{\vec{P}}_{\tilde{i},\tilde{j},\tilde{l}} = \vec{Q}_{i,j,l}
\end{equation}
for all $i=1,\dots,n_\upxi$, $j=1,\dots,n_\upeta$ and $l=1,\dots,n_\upzeta$. We may therefore find $\tilde{\vec{P}}_{i,j,l}$ by solving a system of $3n_\upxi\cdot n_\upeta\cdot n_\upzeta$ equations.

\begin{figure}
	\centering        
	\begin{subfigure}[t]{0.3\textwidth}
		\centering
		\includegraphics[width=0.9\textwidth]{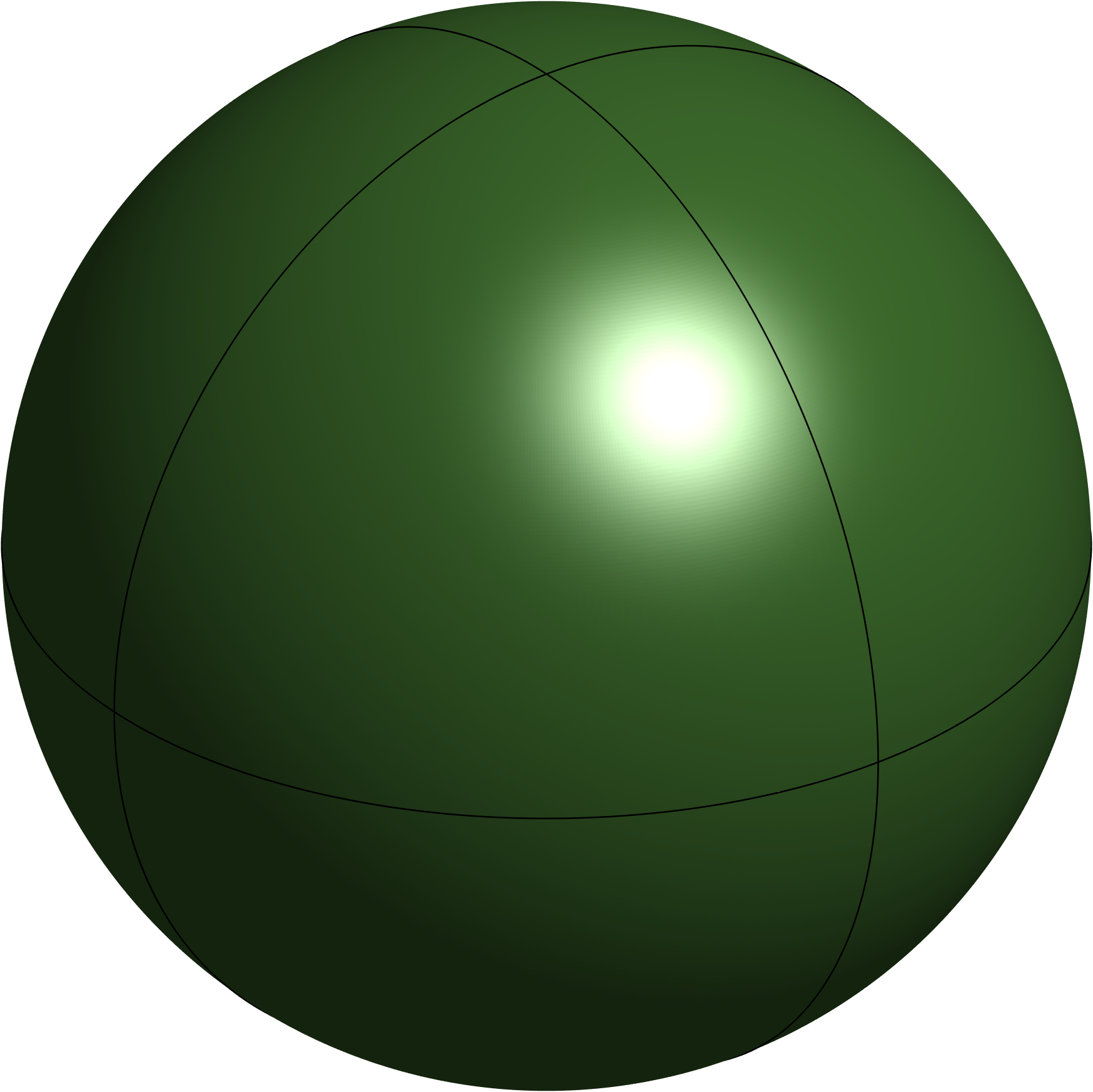}
		\caption{NURBS parametrization.}
	\end{subfigure}
	~
	\begin{subfigure}[t]{0.3\textwidth}
		\centering
		\includegraphics[width=0.9\textwidth]{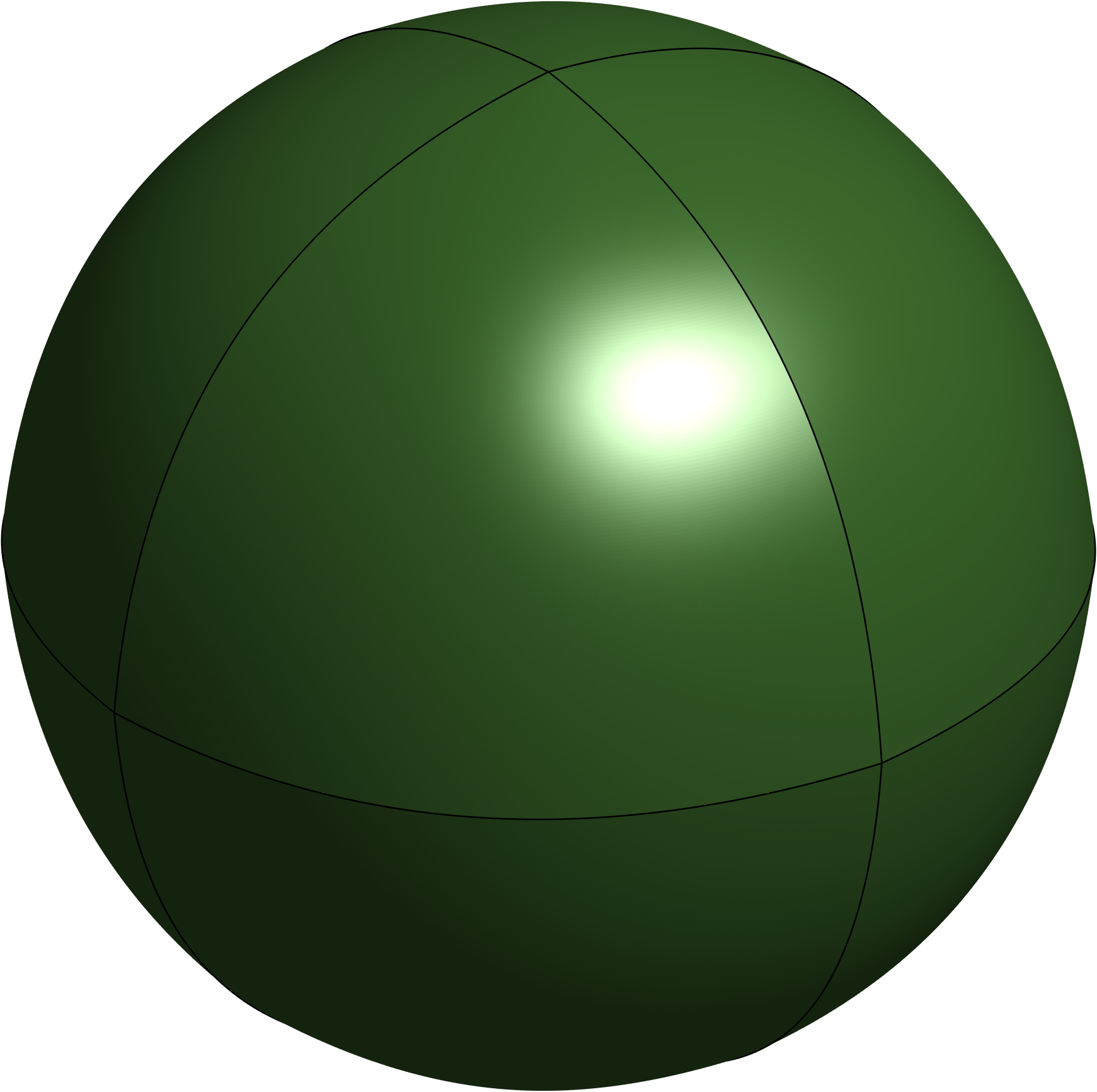}
		\caption{Approximation using $\check{p}_\upxi=\check{p}_\upeta=2$.}
	\end{subfigure} 
	~
	\begin{subfigure}[t]{0.3\textwidth}
		\centering
		\includegraphics[width=0.9\textwidth]{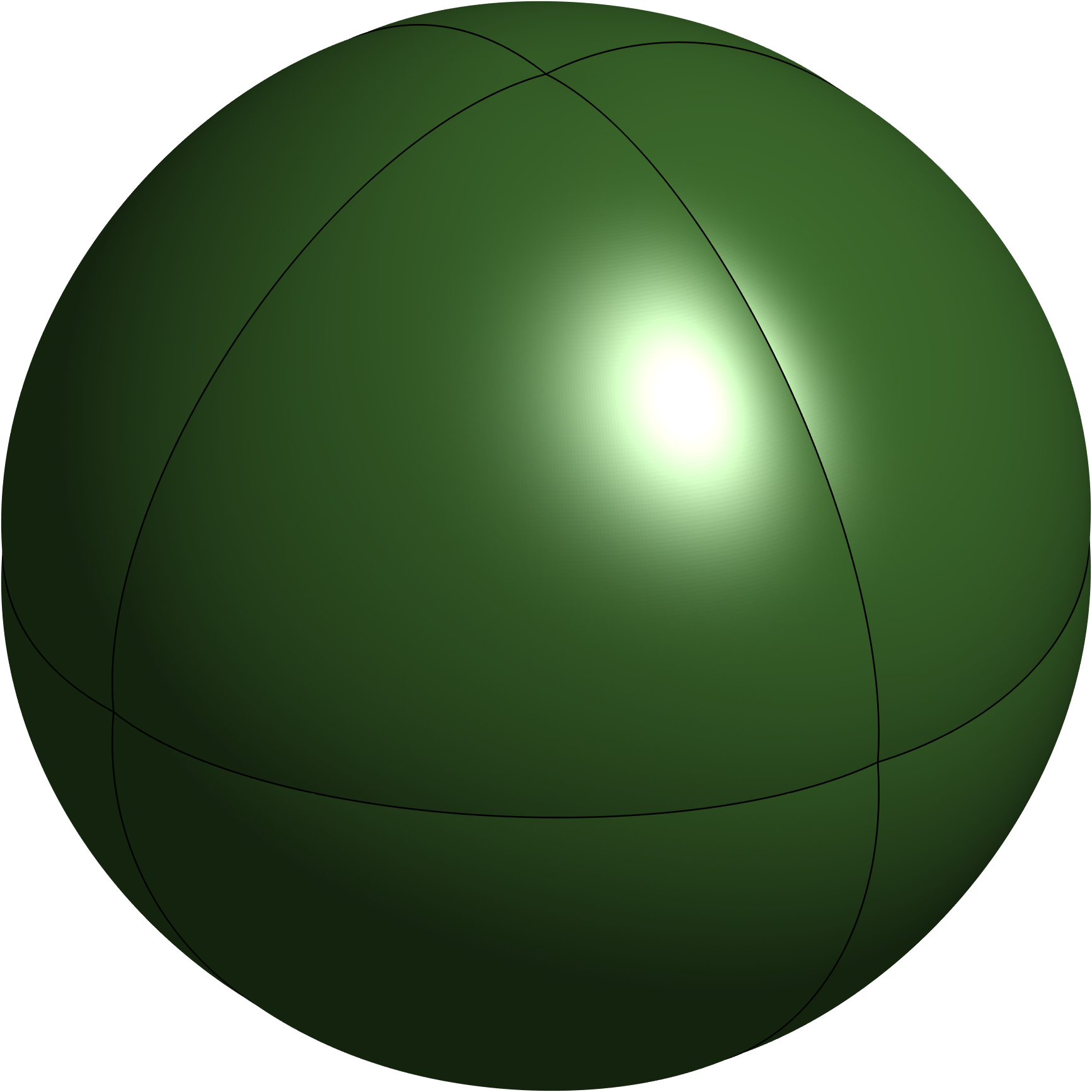}
		\caption{Approximation using $\check{p}_\upxi=\check{p}_\upeta=3$.}
	\end{subfigure}
	\caption{Transformation of an exact NURBS parametrization of a spherical shell to a B-spline approximation of the same geometry.}
	\label{Fig2:NURBStoBsplineTransVis}
\end{figure}

Application of this algorithm to the spherical shell parametrization using NURBS is illustrated in \Cref{Fig2:NURBStoBsplineTransVis}.
